\documentclass[a4paper,12pt]{article}

\usepackage[utf8]{inputenc}

\usepackage{amssymb,amsmath}
\usepackage{amsfonts}
\usepackage{epsfig,wrapfig}
\usepackage{mathptmx}
\usepackage{helvet}
\usepackage{courier}

\usepackage{type1cm}       

\usepackage{makeidx}

\numberwithin{equation}{section}  

\makeindex

\usepackage[bookmarks]{hyperref}

\makeatletter
\renewcommand{\@seccntformat}[1]{\@nameuse{the#1}.\quad}
\makeatother

\font\margofont=cmbx9

\gdef\Section#1{ 
\pagestyle{myheadings} \markboth{Peter Major}{\sc 
\hfill 
\margofont #1\hfill}
\section{#1} 
 }

\begin{document}

\title{Limit theorems for non-linear functionals of stationary 
Gaussian random fields}

\author{Lectures, held by P\'eter Major} 

\date{University of Bochum, 2017}

\maketitle

\Section{The subject of these lectures}

The main subject of this series of talks is the study of
limit theorems for normalized sums of elements in a 
stationary sequence of {\it dependent}\/ random variables in 
such cases when the central limit theorem does not hold for 
them. Because of lack of time I cannot 
give a detailed proof of all results I shall speak about. 
I shall concentrate instead on their content and the 
explanation of the picture behind them. I hope, this will 
be interesting in itself, and it can give considerable
help for those who are interested in a complete proof of 
the results. This can be found in my lecture note 
{\it Multiple Wiener--It\^o Integrals.\/} Lecture Notes 
in Mathematics~{\bf 849}, Revised (augmented) version, 
Springer Verlag, Berlin--Heidelberg--New York, (2014).

Let me describe this problem in more detail. Before 
discussing limit theorems for sums of dependent random 
variables let us recall some facts about limit theorems 
for i.i.d. random variables. There is a natural approach 
to the investigation of these limit theorems where first
we try to find the possible limits by means of the study 
of an appropriately formulated fixed point problem in the 
space of distribution functions. This fixed point problem
can be solved, and it shows that the possible limits are 
the normal and the so-called stable distributions whose
Fourier transforms can be described explicitly.

We want to formulate a natural analogue of this fixed point 
problem which helps to find the possible limit in the case 
of limit theorems for normalized sums of dependent random 
variables. In this case we have to look for the distribution 
of an appropriate random process as the possible limit, 
because the one-dimensional distributions --- unlike in
the case of independent random variables --- do not 
give sufficient information about the behaviour of the 
limit. This leads to the introduction of the notion of
{\it renormalization}\/ and {\it self-similar}\/ fields. 
To define them first we have to introduce some additional 
notions.

I shall consider $\nu$-dimensional stationary random fields.
A $\nu$-dimensional random field is a set of random variables 
$\xi_n$, $n\in {\mathbb Z}_\nu$, where ${\mathbb Z}_\nu$ 
denotes the $\nu$-dimensional integer lattice. I shall call it a
stationary random field if it satisfies the following definition.

\medskip\noindent
{\bf Definition of Discrete (Strictly) Stationary Random Fields.} 
{\it A set of random variables $\xi_n$, $n\in\mathbb Z_\nu$, 
is called a (strictly) stationary discrete random field if
$$
(\xi_{n_1},\dots,\xi_{n_k})\stackrel{\Delta}{=}
(\xi_{n_1+m},\dots,\xi_{n_k+m})
$$
for all $k=1,2,\dots$ and \ $n_1,\dots,n_k,\;m\in\mathbb Z_\nu$,
where $\stackrel{\Delta}{=}$ denotes equality in distribution.}

\medskip
Next I formulate the general limit problem we are interested in.

Given a discrete stationary random field $\xi_n$, 
$n\in\mathbb Z_\nu$, we define for all  
parameters $N=1,2,\dots$  the new stationary random field
\begin{equation}
Z_n^N=A_N^{-1}\sum_{j\in B_n^N}\xi_j, \qquad N=1,2,\dots, \quad
n\in{\mathbb Z}_\nu,
\label{(1.1)}
\end{equation}
where
$$
B_n^N=\{j\colon\; j\in \mathbb Z_\nu,\quad n^{(i)}N\le
j^{(i)}<(n^{(i)}+1)N,\;i=1,2,\dots, \nu\},
$$
and $A_N$, $A_N>0$, is an appropriate norming constant. The 
superscript $i$ denotes the $i$-the coordinate of a vector 
in this formula. We defined in formula~(\ref{(1.1)}) a new
stationary random field $Z_n^N$, $n\in{\mathbb Z}_\nu$, for
all indices~$N$. We are interested in the question when the 
finite dimensional distributions of these random fields 
$Z_n^N$, $n\in{\mathbb Z}_\nu$, called the renormalizations of 
the original field $\xi_n$, $n\in\mathbb Z_\nu$, have a limit 
as  $N\to\infty$, and we want to describe this limit if it 
exists. In particular, we would like to describe those random 
fields $Z_n^*$, $n\in\mathbb Z_\nu$, which appear as the 
limit of such random fields~$Z_n^N$ as $N\to\infty$. This 
problem, which is the natural counterpart of the fixed point 
problem leading to the description of the possible limits in 
the independent case suggested the introduction of the 
following notion.

\medskip\noindent
{\bf Definition of Self-similar (Discrete) Random Fields.} 
{\it A (discrete) random field 
$\xi_n$, $n\in\mathbb Z_\nu$, is called self-similar 
with self-similarity parameter~$\alpha$ if the random
fields $Z^N_n$ defined in~(\ref{(1.1)}) with their help 
and the choice $A_N=N^\alpha$ satisfy the relation
\begin{equation}
(\xi_{n_1},\dots,\xi_{n_k})\stackrel{\Delta}{=}
(Z^N_{n_1},\dots,Z_{n_k}^N) \label{(1.2)}
\end{equation}
for all $N=1,2,\dots$ and $n_1,\dots,n_k\in\mathbb Z_\nu$.}

\medskip
It is natural to expect that the self-similar random fields appear
as the limit fields in the limit problem we are interested in now.
We chose the norming constant $A_N=N^\alpha$, because under some
natural conditions we can satisfy formula~(\ref{(1.2)}) only with
such a choice. We shall consider only such random fields in 
which the random variables $\xi_n$ have a finite second moment. 
This excludes the random fields consisting of independent
random variables with (non-normal) stable distribution from the
classes of self-similar random fields we are interested in.
With this restriction a self-similar random fields with parameter 
$\alpha\neq\frac\nu2$ must consist of strongly dependent random
variables. 

The description of (stationary) self-similar random fields is a very 
hard problem, and we have only partial results. The description 
of the (stationary) Gaussian self-similar random fields and of their 
(Gaussian) domain of attraction is a relatively simple problem, 
because in this case only the correlation function of the 
elements of the random fields has to be studied. This problem
is essentially solved. We want to find non-Gaussian self-similar
random fields and to present such interesting, non-trivial limit 
theorems where they appear as the limit. To find such random fields 
we shall introduce the notion of random fields subordinated to a 
stationary Gaussian random field. We shall work out a method to work 
with such subordinated random fields, and we shall be able to 
construct non-trivial self-similar random  fields and to prove 
interesting limit theorems. To introduce the notion of random 
fields subordinated to a Gaussian random field first we have to 
define the shift transformation determined by a stationary Gaussian 
random field. This will be done in the next consideration.

Let $X_n$, $n\in\mathbb Z_\nu$, be a stationary Gaussian random
field. First we define the shift transformations $T_m$, 
$m\in\mathbb Z_\nu$, over this field by the formula 
$T_mX_n=X_{n+m}$, for all $n,\,m\in\mathbb Z_\nu$. Then we can
extend this shift transformation by means of some results in
measure theory for all such random variables $\xi(\omega)$,
which are measurable with respect to the $\sigma$-algebra
${\cal B}(X_n(\omega),\,n\in\mathbb Z_\nu)$. Indeed, by
some results of measure theory we can write such a random
variable in the form 
$\xi(\omega)=f(X_n(\omega),\,n\in\mathbb Z_\nu)$ with 
some Borel measurable function
$f(x_n,\,n\in\mathbb Z_\nu)$ on the product space 
$R^{\mathbb Z_\nu}$, and we can define with the help
of this representation the shift $T_m$ of the random variable
$\xi(\omega)$ by the formula
$T_m\xi(\omega)=f(X_{n+m}(\omega),\,n\in\mathbb Z_\nu)$. 
It must be still explained that although the function~$f$ is 
not unique in the representation of the random 
variable~$\xi(\omega)$, nevertheless the above definition of 
$T_m\xi(\omega)$ is meaningful. We shall identify two random
variables if they are equal with probability~1. To see that 
we gave a correct definition of the shift transformation we 
have to check that if $f_1(X_n(\omega),\,n\in\mathbb Z_\nu)
=f_2(X_n(\omega),\,n\in\mathbb Z_\nu)$ 
for two functions $f_1$ and $f_2$ with probability~1, then 
also $f_1(X_{n+m}(\omega),\,n\in\mathbb Z_\nu)
=f_2(X_{n+m}(\omega),\,n\in\mathbb Z_\nu)$
with probability~1 because of the stationarity of the random
field~$X_n(\omega)$, $n\in\mathbb Z_\nu$. 

In such a way we have defined the shift transformation
for all random variables measurable with respect to the 
$\sigma$-algebra generated by the random field 
$X_n,\;n\in\mathbb Z_\nu$. But since we shall work
only with random variables of finite second moment we shall 
consider the action of the shift transformation only for
a smaller class of random variables, for the elements of 
the Hilbert space $\cal H$ introduced below.

Let ${\cal H}$ denote the {\it real} Hilbert space consisting of 
the square integrable random variables measurable with respect 
to the $\sigma$-algebra
${\cal B}={\cal B}(X_n,\;n\in\mathbb Z_\nu)$. The scalar 
product in ${\cal H}$ is defined as $(\xi,\eta)=E\xi\eta$, 
$\xi,\,\eta\in{\cal H}$. We define the shift transformations $T_m$, 
$m\in\mathbb Z_\nu$, for the elements of $\cal H$ in the
way as we have done before in the general case. It is not 
difficult to check that the shift transformations $T_m$, 
$m\in\mathbb Z_\nu$, map the elements of $\cal H$ to
another element of $\cal H$, they are unitary, i.e. norm
preserving, invertible linear transformations. Moreover, they 
constitute a unitary group in ${\cal H}$, i.e. $T_{n+m}=T_nT_m$ 
for all $n,\,m\in\mathbb Z_\nu$, and $T_0=\textrm{Id}$. 
Now we introduce the following

\medskip\noindent
{\bf Definition of Subordinated Random Fields.} 
{\it Given a stationary Gaussian random field 
$X_n$, $n\in\mathbb Z_\nu$, we define the Hilbert space
${\cal H}$ and the shift transformations $T_m$, 
$m\in\mathbb Z_\nu$, over ${\cal H}$ as before. A 
discrete stationary field $\xi_n$ is called a random field 
subordinated to $X_n$ if $\xi_n\in{\cal H}$, and
$T_n\xi_m=\xi_{n+m}$ for all $n,\,m\in\mathbb Z_\nu$.}

\medskip
One of the main task of this series of talks is to work out
a good method that enables us to study the fields subordinated
to a stationary Gaussian random field together with the shift
transformation acting on it. This enables us both to find
non-trivial self-similar fields and to prove interesting
limit theorems. This program will consist of several steps.
First we study the underlying Gaussian fields. There is a
classical result in analysis that enables us to describe 
the correlation function of this Gaussian field as the
Fourier transform of a so-called spectral measure. We show
that a so-called random spectral measure can be constructed,
and a natural random integral can be defined with respect to 
it in such a way that the elements of the Gaussian random
field can be expressed in a form that can be interpreted
as the random Fourier transform of the random spectral measure.
Then we introduce a multiple random integral, called 
multiple Wiener--It\^o random integral with respect to the 
random spectral measure that enables us to express all
elements of the above defined Hilbert space $\cal H$ as
the sum of multiple Wiener--It\^o integrals of different
multiplicity. Moreover, this representation is unique,
and the action of the shift transformation on $\cal H$
can be calculated in a simple way with its help. Then 
we make a most important step in our investigation, we 
prove the so-called diagram formula that enables us to 
rewrite the product of Wiener--It\^o integrals as the sum 
of Wiener--it\^o integrals of different multiplicity. These 
results together with some basic facts about Hermite 
polynomials make possible to work out a technique that 
enables us to construct non-trivial self-similar fields,
and to prove non-trivial (non-central) limit theorems.

In the above discussion I dealt with discrete time stationary 
random fields. However it is useful to handle discrete time 
Gaussian self-similar random fields together with their 
continuous time versions, since they yield --- because of 
their stronger symmetry properties --- an essential help 
also in the study of discrete time stationary random fields.
However, in the study of continuous time stationary fields 
serious additional technical difficulties appear, because, 
as it turned out, it is more useful to work with generalized 
and not with classical continuous time random fields. The 
study of  generalized random fields demands some additional 
work. In particular, we have to present the appropriate notions 
and results needed in their study. Besides,we have to explain 
why we want to work with generalized random fields, why are 
the classical random fields inappropriate for us.

\medskip
To give some feeling why the elaboration of the above theory
is useful for us let us consider some limit problems which
we can handle with its help.

Let $X_n$, $n\in\mathbb Z_\nu$, be a stationary Gaussian 
field with expectation $EX_n=0$ and some correlation function 
$r(n)=EX_mX_{n+m}$, $n,m\in\mathbb Z_\nu$, and take the
limit problem introduced in formula (\ref{(1.1)}) with 
appropriate norming constants $A_N$ if the random variables
$X_j$ play the role of the random variables $\xi_j$. This
is a relatively simple problem, because it is enough to
check whether the correlation functions $r_N(n)=EZ_m^NZ_{m+n}^N$
have a limit as $N\to\infty$ with an appropriate norming
constants $A_N$. Next we are considering the following
harder problem. Let us consider the previous Gaussian random
field $X_n$, $n\in\mathbb Z_\nu$, and define with the
help of some real valued function $f(x)$, $x\in R$ such that
$Ef(X_n)=0$ and $Ef^2(X_n)<\infty$ a new stationary random
field $\xi_n=f(X_n)$, $n\in\mathbb Z_\nu$. Now we are
interested in whether the random field $Z_n^N$, defined
in formula (\ref{(1.1)}) with the help of this new random 
field $\xi_n$ have a limit with an appropriate norming 
constants $A_N$ as $N\to\infty$. This is a considerably 
harder problem, and at the first sight we may have no 
idea how to handle it.

But let us observe that this random field $\xi_n$, 
$n\in\mathbb Z_\nu$, is subordinated to the
Gaussian random field $X_n$, $n\in\mathbb Z_\nu$,
hence we can apply the theory worked out to study this
field. This theory enables us to present the normalized
random sums $Z_N^n$ in such a form that indicates what
kind of limit these random variables may have, and how to
chose the norming constants $A_N$ to get a limit. Some 
(natural) calculation shows that the limit theorems 
suggested by the formulas obtained by our theory really 
hold.

\Section{Random spectral measures}

In the first step of our study we give a useful representation
of the elements of the stationary Gaussian random field we are 
working with by means of a random integral. To get it first we 
apply the corollary of a classical result of analysis, --- 
called Bochner's theorem ---, about a good representation of 
the so-called positive definite functions by means of Fourier 
transform. This enables us to describe the correlation function 
of a stationary Gaussian random field as the Fourier transform 
of a finite measure which is called the spectral measure of this 
field in the literature. Because of lack of time I omit the 
discussion of the underlying Bochner theorem, I will formulate 
only the result about the description of the correlation 
function of a stationary random field as the Fourier transform 
of the spectral measure of this field, because we need only this 
result. Then I show that a so-called random spectral measure can 
be constructed with the help of this result, and the elements of 
our stationary Gaussian random field can be represented as a 
random integral with respect to this random spectral measure. 
In an informal way this statement can be interpreted so that 
while the correlation function of a stationary Gaussian random 
field can be represented as the Fourier transform of its 
spectral measure, the random variables of this field can be 
represented as random Fourier transforms of the random 
spectral measure.

The introduction of the random spectral measure turned out to 
be very useful. The representation of our random variables by 
means of a random integral with respect to it enables us to 
work well with the shift transformations of our random field. 
Moreover, multiple Wiener--It\^o integrals can be defined 
with respect to the random spectral measure, and they have an 
important role in our considerations. All elements of the 
Hilbert space $\cal H$ consisting of those random variables 
with finite second moment which are measurable with respect 
to the $\sigma$-algebra generated by the random variables of 
the stationary Gaussian random field we are working with can 
be written as the sum of multiple Wiener--It\^o integrals of
different multiplicity. This representation is very useful, 
and it will be the main tool in our investigation.

First I formulate the result about the description of the
correlation function of discrete time stationary random
fields.

\medskip\noindent
{\bf Theorem 2A about the spectral representation of the 
correlation function of a discrete stationary random field.}
{\it Let $X_n$, $n\in{\mathbb Z}_\nu$, be a discrete 
(Gaussian) stationary random field with expectation $EX_n=0$,
$n\in {\mathbb Z}_n$. There exists a unique
finite measure $G$ on $[-\pi,\pi)^\nu$ such that the 
correlation function $r(n)=EX_0X_n=EX_kX_{k+n}$, 
$n\in{\mathbb Z}_\nu$, $k\in {\mathbb Z}_\nu$, can be 
written in the form
\begin{equation}
r(n)=\int e^{i(n,x)}G(\,dx), \label{(3.1)}
\end{equation}
where $(\cdot,\cdot)$ denotes scalar product. Further, 
$G(A)=G(-A)$ for all $A\in [-\pi,\pi)^\nu$.}

\medskip\noindent
We can identify $[-\pi,\pi)^\nu$ with the torus
$R^\nu/2\pi{\mathbb Z}_\nu$. Thus
e.g. $-(-\pi,\dots,-\pi)=(-\pi,\dots,-\pi)$.

\medskip
We want to formulate the continuous time version of the 
above result about the representation of the correlation
function as the Fourier transform of a so-called 
spectral measure. There exists such a theorem, and it 
is very similar to the discrete time result. The only 
difference is that in the new case we have to consider 
a spectral measure $G$ on $R^\nu$ and not on
$[-\pi,\pi)^\nu$ as in formula (\ref{(3.1)}), and we 
have to assume that the correlation function 
$r(t)=EX_0X_t$ is continuous. We could work with such 
a result. Nevertheless, we shall follow a different
approach. We shall work instead of classical continuous 
time stationary random fields with (stationary) 
generalized random fields, and formulate the results 
for them. This demands the introduction of some new 
notions and some extra explanation.

First we have to understand why the classical continuous 
time stationary random fields are not good for us, why 
do we want to work with generalized random fields 
instead. R.~L.~Dobrushin gave the following informal 
explanation for this. The trajectory of a continuous
time random field can be very bad. It can be so bad, 
that it simply does not exist. In such cases we cannot 
consider our random field as a really existing field, 
but there may be a possibility to consider it as a 
generalized random field. As we shall later see, there 
are very important generalized random fields that cannot 
be obtained by means of a classical random field. Moreover, 
they play a useful role also in the study of discrete 
time random fields.

The following heuristic argument may explain the definition of
 generalized random fields. Let us have a classical continuous 
time random field $X(t)$ with parameters $t\in{R}^\nu$
in the $\nu$-dimensional Euclidean space, and a linear 
topological space $\cal F$ of functions on ${R}^\nu$ with 
some nice properties. In nice cases we can define the integral 
$X(\varphi)=\int_{R^\nu} \varphi(s)X(s)\,ds$ for all functions
$\varphi\in{\cal F}$, and if the space of functions $\cal F$ 
is sufficiently rich, then the random variables $X(\varphi)$
determine the values of the random variables $X(t)$ from which
we obtained it. The random variables $X(\varphi)$ have the property 
$X(a\varphi+b\psi)=aX(\varphi)+bX(\psi)$, and if $\varphi$
and $\psi$ are two functions of $\cal F$ which are close
to each other, then it is natural to expect that $X(\varphi)$
and $X(\psi)$ are also close to each other in some sense.

This means that we can correspond to a classical random
field $X(t)$ a class of random linear functionals 
$X(\varphi)$, $\varphi\in{\cal F}$, indexed by the 
elements of $\cal F$, which have some nice properties. 
We will call a class of continuous random functionals  
indexed by the elements of a nice linear topological
space $\cal F$ a generalized random field. We can 
correspond to each classical random field a generalized 
random field which determines it, but not all 
generalized random fields can be obtained in such a way. 
Now I introduce the precise definition of generalized 
random fields together with some additional notions we 
shall apply in our considerations.

\subsection{Generalized random fields and some related notions}

First I introduce the linear topological space
we shall be working with in the definition of
generalized random fields. There are several good choices
for it. I shall use the so-called Schwartz space 
$\cal S$, because we can work with it very well.

We define the Schwartz space ${\cal S}$ together
with its version ${\cal S}^c$ consisting of complex
valued function. The space ${\cal S}^c=({\cal S}_\nu)^c$ 
consists of those complex valued functions of $\nu$ 
variables which decrease at infinity, together with 
their derivatives, faster than any polynomial degree. 
More explicitly, $\varphi\in{\cal S}^c$ for a complex 
valued function $\varphi$ defined on $R^\nu$ if
$$
\left|x_1^{k_1}\cdots x_\nu^{k_\nu}\frac{\partial^{q_1+\cdots+q_\nu}}
{\partial x_1^{q_1}\dots \partial x_\nu^{q_\nu}}
\varphi(x_1,\dots,x_\nu)\right|
\le C(k_1,\dots,k_\nu,q_1,\dots,q_\nu)
$$
for all point $x=(x_1,\dots,x_\nu)\in R^\nu$ and  vectors
$(k_1,\dots,k_\nu)$, $(q_1,\dots,q_\nu)$ with non-negative
integer coordinates with some constant
$C(k_1,\dots,k_\nu,q_1,\dots,q_\nu)$ which may depend on the
function~$\varphi$. The elements of the space ${\cal S}$ 
are defined similarly, with the only difference that they 
are real valued functions.

To complete the definition of the the spaces ${\cal S}$ 
and ${\cal S}^c$ we still have to define the topology
in them. We introduce the following topology in these
spaces. 

Let a basis of neighbourhoods of the origin consist
of the sets
$$
U(k,q,\varepsilon)=\left\{\varphi\colon\;
\max_x(1+|x|^2)^k |D^q\varphi(x)|<\varepsilon\right\}
$$
with $k=0,1,2,\dots$, $q=(q_1,\dots,q_\nu)$ with non-negative 
integer coordinates and $\varepsilon>0$, where 
$|x|^2=x_1^2+\cdots+x_\nu^2$, and 
$D^q=\frac{\partial^{q_1+\cdots+q_\nu}}
{\partial x_1^{q_1}\dots \partial x_\nu^{q_\nu}}$.
A basis of neighbourhoods of an arbitrary function 
$\varphi\in{\cal S}^c$ (or $\varphi\in{\cal S}$) consists of 
sets of the form $\varphi+U(k,q,\varepsilon)$, where the 
class of sets $U(k,q,\varepsilon)$ is a basis of 
neighbourhood of the origin. Let me remark that a sequence of 
functions $\varphi_n\in{\cal S}^c$ (or $\varphi_n\in{\cal S}$) 
converges to a function $\varphi$ in this topology if and 
only if
$$
\lim_{n\to\infty}\sup_{x\in R^\nu}
(1+|x|^2)^k|D^q\varphi_n(x)-D^q\varphi(x)|=0.
$$
for all $k=1,2,\dots$ and $q=(q_1,\dots,q_\nu)$.
The limit function $\varphi$ is also in the
space~${\cal S}^c$ (or in the space ${\cal S}$).

I shall define the generalized random fields and some
related notions with the help of the notion of Schwartz spaces.

\medskip\noindent
{\bf Definition of Generalized Random Fields.} 
{\it We say that the set of random variables $X(\varphi)$,  
$\varphi\in{\cal S}$, is a generalized random field over the 
Schwartz space ${\cal S}$ of rapidly decreasing, smooth functions if:

\medskip
\begin{description}
\item[(a)]  $X(a_1\varphi_1+a_2\varphi_2)=a_1X(\varphi_1)+a_2X(\varphi_2)$
with probability 1 for all real numbers $a_1$ and $a_2$ and
$\varphi_1\in{\cal S}$, $\varphi_2\in{\cal S}$. (The exceptional set of
probability~0 where this identity does not hold may depend on $a_1$,
$a_2$, $\varphi_1$ and $\varphi_2$.)
\item[(b)] $X(\varphi_n)\Rightarrow X(\varphi)$ stochastically if
$\varphi_n\to\varphi$ in the topology of ${\cal S}$.
\end{description}

}
\medskip
We also introduce the following definitions.

\medskip\noindent
{\bf Definition of Stationarity and Gaussian Property of a 
Generalized Random Field and the Notion of Convergence of 
Generalized Random Fields in Distribution.} 
{\it The generalized random field
$X=\{X(\varphi),\,\varphi\in {\cal S}\}$ is stationary if
$X(\varphi)\stackrel{\Delta}{=}X(T_t\varphi)$ for all 
$\varphi\in{\cal S}$ and $t\in R^\nu$, where 
$T_t\varphi(x)=\varphi(x-t)$. It is Gaussian
if $X(\varphi)$ is a Gaussian random variable for all
$\varphi\in{\cal S}$. The relation
$X_n\stackrel{{\cal D}}{\rightarrow} X_0$ as $n\to\infty$ 
holds for a sequence of generalized random fields $X_n$, 
$n=0,1,2,\dots$, if 
$X_n(\varphi)\stackrel{{\cal D}}{\rightarrow} X_0(\varphi)$ 
for all $\varphi\in{\cal S}$, where
$\stackrel{{\cal D}}{\rightarrow}$ denotes convergence in 
distribution.}

\medskip
Next I formulate the version of our limit problem for
generalized fields.

\medskip
Given a stationary generalized  random field $X$ and a function
$A(t)>0$, $t>0$, on the set of positive real numbers we define the
(stationary) random fields $X^A_t$ for all $t>0$ by the formula
\begin{equation}
X^A_t(\varphi)=X(\varphi_t^A), \quad \varphi\in{\cal S}, \qquad
\textrm{where } \varphi_t^A(x)=A(t)^{-1}\varphi\left(\frac xt\right).
\label{(1.3)}
\end{equation}

We are interested in the following

\medskip\noindent
{\bf Question.} {\it When does a generalized random field $X^*$ 
exist such that $X_t^A\stackrel{{\cal D}}{\rightarrow} X^*$ as 
$t\to\infty$ (or as $t\to0$)?}

\medskip
In relation to this question we introduce the following

\medskip\noindent
{\bf Definition of Self-similarity.} 
{\it The stationary generalized random field $X$ is self-similar 
with self-similarity parameter $\alpha$ if  
$X^A_t(\varphi)\stackrel{\Delta}{=} X(\varphi)$ for all
$\varphi\in{\cal S}$  and $t>0$ with the function 
$A(t)=t^\alpha$.}

\medskip
To answer the above question one should first describe the 
generalized self-similar random fields.

We define analogously to the case of discrete random fields
the notion of generalized subordinated random fields. 

Let $X(\varphi)$, $\varphi\in{\cal S}$, be a generalized 
stationary Gaussian random field. The formula
$T_tX(\varphi))=X(T_t\varphi)$, $T_t\varphi(x)=\varphi(x-t)$,
 defines the shift transformation for all $t\in R^\nu$. Let
${\cal H}$ denote the real Hilbert space consisting of the
${\cal B}={\cal B}(X(\varphi),\;\varphi\in{\cal S})$ measurable 
random variables with finite second moment. The shift 
transformation can be extended to a group of unitary 
transformations over ${\cal H}$ similarly to the discrete case.
(This definition has the following idea. If a random variable
$\xi\in{\cal H}$ has the form 
$\xi=F(X(\varphi_1),\dots,X(\varphi_s))$
with some functions $\varphi_1,\dots\varphi_s\in{\cal S}$,
and a measurable function $F$ of $s$ variables, then we define
$T_t\xi=F(X(T_t\varphi_1),\dots,X(T_t\varphi_s))$. A general
random variable $\xi\in{\cal H}$ has a somewhat more complicated,
but similar representation, and its shift can be defined in a
similar way.) With the help of the notion of shift transformations 
in generalized fields we can introduce the following definition.

\medskip\noindent
{\bf Definition of Generalized Random Fields Subordinated to a
Generalized Stationary Gaussian Random Field.} 
{\it Given a generalized stationary Gaussian random field 
$X(\varphi)$, $\varphi\in{\cal S}$, we define the Hilbert 
space ${\cal H}$ and the shift transformations $T_t$, 
$t\in R^\nu$, over ${\cal H}$ as above. A generalized 
stationary random field $\xi(\varphi)$, $\varphi\in{\cal S}$, 
is subordinated to the field $X(\varphi)$,
$\varphi\in{\cal S}$, if $\xi(\varphi)\in{\cal H}$ and
$T_t\xi(\varphi)=\xi(T_t\varphi)$ for all $\varphi\in{\cal S}$ 
and $t\in R^\nu$, and $E[\xi\varphi_n)-\xi(\varphi)]^2\to0$ if
$\varphi_n\to\varphi$ in the topology of ${\cal S}$.}

\medskip
Now we can formulate the analogue of Theorem 2A about the
Fourier representation of the correlation function of a 
generalized field. Before doing it I recall an important
property of the Fourier transform of the functions in the
Schwartz spaces $\cal S$ and ${\cal S}^c$. Actually this
property of the Schwartz spaces made useful their application
in the definition of generalized fields. 

The Fourier transform $f\to\tilde f$ is a bicontinuous map 
from ${\cal S}^c$ to~${\cal S}^c$. (This means that this 
transformation is invertible, and both the Fourier transform 
and its inverse are continuous maps from ${\cal S}^c$ to 
${\cal S}^c$.) (The restriction of the Fourier transform to 
the space ${\cal S}$ of real valued functions is a bicontinuous 
map from ${\cal S}$ to the subspace of ${\cal S}^c$ consisting 
of those functions $f\in{\cal S}^c$ for which
$f(-x)=\overline{f(x)}$ for all $x\in R^\nu$.) I omit the
proof of this statement, I only remark that the smoothness
properties of the functions in ${\cal S}^c$ imply the fast
decrease of their Fourier transform at infinity, and their 
fast decrease at infinity imply the smoothness properties
of their Fourier transform. 

In a thorough analysis one also studies the properties of
the elements of the space of generalized functions ${\cal S}'$ 
which are the continuous linear functionals over $\cal S$.
But since they are needed only in such proofs which  I omit
in this discussion, I do not discuss these problems. Next
I formulate the following result.

\medskip\noindent
{\bf Theorem 2B about the spectral representation of the correlation
function of a generalized stationary field.} 
{\it Let $X(\varphi)$, $\varphi\in{\cal S}$, $EX(\varphi)=0$ for
$\varphi\in{\cal S}$, be a generalized  Gaussian stationary 
random field over ${\cal S}={\cal S}_\nu$. There exists a 
unique $\sigma$-finite measure $G$ 
on $R^\nu$ such that
\begin{equation}
EX(\varphi)X(\psi)
=\int\tilde\varphi(x)\,\bar{\!\tilde\psi}(x)G(\,dx)
\qquad \textrm{for all } \varphi,\,\psi\in{\cal S}, \label{(3.2)}
\end{equation}
where $\,\tilde{}\,$ denotes Fourier transform and $\,\bar{}\,$
complex conjugate. The measure $G$ has
the properties $G(A)=G(-A)$ for all $A\in{\cal B}^\nu$, and
\begin{equation}
\int (1+|x|)^{-r}G(\,dx)<\infty \quad\textrm{with an appropriate } r>0.
\label{(3.3)}
\end{equation}
}

\medskip\noindent
Let me remark that while in Theorem 2A the spectral measure $G$
had to be finite, in Theorem~2B it had to satisfy a much weaker
condition~(\ref{(3.3)}). This indicates that there are such
generalized stationary random fields that cannot be obtained
from non-generalized random fields. This difference between
the properties of the spectral measures in Theorems~2A and~2B
also have other interesting and important consequences about
which I shall write at the end of this section. 

The proof that the correlation function of a generalized field
must satisfy Theorem~2B depends on some deep theorems about
generalized functions, hence I omit it. On the other hand, I 
briefly show that in Theorem~2B we defined really the correlation
function of a stationary generalized Gaussian random field. Before 
doing this I present a short calculation that indicates that 
formula~(\ref{(3.2)}) can be considered as the natural analogue
of formula~(\ref{(3.1)}) when we are working with generalized
field.

Let us consider a continuous time Gaussian stationary field
$X(t)$, $t\in R^\nu$, with correlation function 
$EX(s)X(t)=\int e^{i(s-t,x)}G(\,dx)$, and let us calculate the
correlation function $EX(\varphi)X(\psi)$, 
$\varphi,\psi\in{\cal S}$, where 
$X(\varphi)=\int \varphi(t)X(t)\,dt$,
and
$X(\psi)=\int \psi(t)X(t)\,dt$. We have

\begin{eqnarray*}
EX(\varphi)X(\psi)&=& E\int \varphi(s)X(s)\,ds\int \psi(t)X(t)\,dt \\
&=&\int\int \varphi(s) \psi(t)EX(s)X(t)\,ds\,dt \\
&=&\int\int \varphi(s) \psi(t) \left[\int e^{i(s-t,x)}G(\,dx)\right]\,ds\,dt\\
&=&\int\left[\int e^{(i(s,x)}\varphi(s)\,ds\right]
\left[\int e^{-i(t,x)}\psi(t)\,dt\right]G(\,dx) \\
&=&\int\tilde\varphi(x)\,\bar{\!\tilde\psi}(x)G(\,dx).
\end{eqnarray*}

Next we show that formulas (\ref{(3.2)}) and (\ref{(3.3)}) in Theorem~2.1
define the correlation function of a generalized stationary Gaussian
random field.

First we show that $EX(\varphi)X(\psi)$ defined in
(\ref{(3.2)}) is a real number for all $\varphi,\psi\in\cal S$.
To show this we apply the change of variables $x\to -x$ in this
formula, and we exploit that $G(A)=G(-A)$, and
$\bar{\tilde f}(x)=\tilde f_-(x)$ with $f_-(x)=f(-x)$ for a real
valued function~$f$. This implies that 
$EX(\varphi)X(\psi)=\overline{EX(\varphi)X(\psi)}$ 
i.e. $EX(\varphi)X(\psi)$ is a real number.

By Kolmogorov's existence theorem a random process with
prescribed finite dimensional distributions exists if 
these distributions are consistent. By this result to prove 
that there is a Gaussian random field $X(\varphi)$, 
$\varphi\in\cal S$, with expectation zero and correlation 
function $EX(\varphi)X(\psi)$ defined in (\ref{(3.3)}) it 
is enough to show that for arbitrary finite set of 
functions $\varphi_1,\dots\varphi_n\in\cal S$ the matrix 
$(d_{j,k})$, $1\le j,k\le n$,  
$di_{j,k}=EX(\varphi_j)X(\varphi_k)$ is positive semidefinite.
This is equivalent to the statement that for any function $\psi$
of the form $\psi(x)=c_1\varphi_1(x)+\cdots+c_n\varphi_n(x)$ 
with real numbers $c_1,\dots,c_n$ the expression 
$EX(\psi)X(\psi)$ defined in (\ref{(3.3)}) is non-negative.
This fact can be simply checked.

We also have to show that a random field with such a 
distribution is a generalized field, i.e. it satisfies 
properties~(a) and~(b) given in the definition of 
generalized fields.

Property~(a) holds, because, as it is not difficult 
to check with the help of formula~(\ref{(3.2)}),
\begin{eqnarray*}
&&E[a_1X(\varphi_1)+a_2X(\varphi_2)
-X(\varphi(a_1\varphi_1+a_2\varphi_2)]^2\\
&&\qquad =\int\left|a_1\tilde\varphi_1(x)+a_2\tilde\varphi_2(x)
-(\widetilde{a_1\varphi_1+a_2\varphi_2})(x)\right|^2G(\,dx)=0.
\end{eqnarray*}

It is not difficult to show that if $\varphi_n\to\varphi$ 
in the topology of the space ${\cal S}$, then 
$E[X(\varphi_n)-X(\varphi)]^2
=\int|\tilde\varphi_n(x)-\tilde\varphi(x)|^2 G(\,dx)\to0$ 
as $n\to\infty$, hence property~(b) also holds. (Here we exploit 
that the transformation $\varphi\to\tilde\varphi$ is 
bicontinuous in the space~${\cal S}$.) 

It is clear that the Gaussian random field constructed in 
such a way is stationary.

Finally, I remark that some additional investigation
shows that the correlation function $EX(\varphi)X(\psi)$
uniquely determines the spectral measure $G$ in formula
(\ref{(3.3)}), since the class of functions $\cal S$
is sufficiently rich. 

\subsection{Construction of random spectral measures}

We shall construct generalized spectral measures both for
discrete valued and generalized stationary Gaussian random 
fields. The construction in the two cases is similar, but there
is some difference between them. In both cases we construct
an appropriate unitary operator $I$, and we define the
random spectral measure with its help. Let me remark
that I shall speak also about unitary operators between two
different Hilbert spaces. Given two Hilbert spaces
${\cal H}_0$ and ${\cal H}_1$ I call a linear transformation 
$I\colon\; {\cal H}_0\to {\cal H}_1$ unitary if it is 
norm preserving and invertible. We shall define the 
operator $I$ in the case of discrete and generalized
fields in a similar way. First we define them on a dense 
subspace of ${\cal H}_0$, and then we extend it to the 
whole space in a natural way.

First I define the Hilbert spaces ${\cal H}_0$ and
${\cal H}_1$, (more precisely its complexification
${\cal H}^c_1$ we shall work with) both in the discrete
and generalized random field case.

Let us consider a stationary Gaussian random field 
(discrete or generalized one) with spectral measure $G$. We 
shall denote the space $L_2([-\pi,\pi)^\nu,{\cal B}^\nu,G)$ 
or $L_2(R^\nu,{\cal B}^\nu,G)$ simply by $L^2_G$. This will
play the role of the Hilbert space ${\cal H}_0$. (The space
${\cal H}_0$ contains also complex valued functions.) 

Given a stationary Gaussian random field, either a discrete field,
 $X_n$, $n\in{\mathbb Z}_n$, or a generalized one $X(\varphi)$,
$\varphi\in\cal S$, first we define a real Hilbert
space ${\cal H}_1$ and then its complexification ${\cal H}_1^c$. 
The real Hilbert space ${\cal H}_1$ is that subspace of the 
Hilbert space of square integrable random  variables (in the 
probability space we are working with) which is generated by 
the finite linear combination of the random variables 
$X_n$, $n\in{\mathbb Z}_n$, in the discrete field case, and by
$X(\varphi)$, $\varphi\in{\cal S}$, in the generalized field case.
We define its complexification ${\cal H}^c_1$ in the following 
way. The elements of ${\cal H}^c_1$ are of the form $X+iY$, \ 
$X,\,Y\in{\cal H}_1$, and the scalar product is defined in it 
as $(X_1+iY_1,X_2+iY_2)=EX_1X_2+EY_1Y_2+i(EY_1X_2-EX_1Y_2)$.
We are going to construct a unitary transformation $I$ 
from $L^2_G$ to ${\cal H}^c_1$. We shall define the random 
spectral measure with the help of this transformation.

I recall that we also introduced the Schwartz space 
${\cal S}^c$  consisting of  the rapidly decreasing, 
smooth, {\it complex valued} functions with the usual 
topology of the Schwartz space. 
It can be proved with the help of some results in analysis
that the set of finite trigonometrical polynomials
$\sum c_ne^{i(n,x)}$ are dense in $L_G^2$ in the discrete field, 
and the functions $\varphi\in{\cal S}^c$ are dense in $L_G^2$ in 
the generalized field case. We shall exploit this fact in our
construction, 

We define the mapping
\begin{equation}
I\left(\sum c_n e^{i(n,x)}\right)=\sum c_nX_n \label{(3.4)}
\end{equation}
in the discrete field case, where the sum is finite, and
\begin{equation}
I(\widetilde{\varphi+i\psi)}=X(\varphi)+iX(\psi),
\quad \varphi,\,\psi\in{\cal S} \label{($3.4'$)}
\end{equation}
in the generalized field case.

Simple calculation with the help of Theorems~2A and~2B shows that
\begin{eqnarray*}
\left\|\sum c_ne^{i(n,x)}\right\|_{L^2_G}^2
&=&\sum\sum c_n\bar c_m\int e^{i(n-m,x)}G(\,dx)\\
&=&\sum\sum c_n\bar c_m EX_nX_m=E\left|\sum c_nX_n\right|^2,
\end{eqnarray*}
and
\begin{eqnarray*}
\|\widetilde{\varphi+i\psi}\|^2_{L^2_G}
&=&\int[\tilde\varphi(x)\bar{\tilde\varphi}(x)
-i\tilde\varphi(x)\bar{\tilde\psi}(x)
+i\tilde\psi(x)\bar{\tilde\varphi}(x)
+\tilde\psi(x)\bar{\tilde\psi}(x)]G(\,dx)\\
&=&EX(\varphi)^2-iEX(\varphi)X(\psi)+iEX(\psi)X(\varphi)
+EX(\psi)^2 \\
&=&E\left(|X(\varphi)+iX(\psi)|\right)^2.
\end{eqnarray*}
This means that the mapping $I$ from a linear subspace of $L_G^2$
to ${\cal H}_1^c$ is norm preserving. Besides, the subspace 
where $I$ was defined is dense in $L^2_G$, Hence the mapping 
$I$ can be uniquely extended to a norm preserving transformation 
from $L^2_G$, to ${\cal H}^c_1$. Since the random variables 
$X_n$ or $X(\varphi)$ are obtained as the image of some 
element from $L_G^2$ under this transformation, the range of~$I$
is the whole space ${\cal H}_1^c$, and $I$ is a 
unitary transformation from $L^2_G$ i.e. from ${\cal H}_0$ to
${\cal H}^c_1$. A unitary transformation preserves not only the
norm, but also the scalar product. Hence
$\int f(x)\bar g(x)G(\,dx)=EI(f)\overline{I(g)}$ for all
$f,\,g\in L^2_G$.

We shall define the random spectral measure $Z_G(A)$ related to
our Gaussian stationary random field for those Borel measurable 
sets $A\in{\cal B}^\nu$ for which $G(A)<\infty$, and it is 
defined by the formula
$$
Z_G(A)=I(\chi_A),
$$
where $\chi_A$ denotes the indicator function of the set~$A$. It 
is not difficult to see that

\medskip
\begin{description}
\item [(i)] The random variables  $Z_G(A)$ are complex valued, 
jointly Gaussian random variables. (The random variables 
$\textrm{Re}\, Z_G(A)$ and $\textrm{Im}\, Z_G(A)$ with possibly 
different sets~$A$ are jointly Gaussian.)
\item [(ii)] $EZ_G(A)=0$,
\item [(iii)] $EZ_G(A)\overline {Z_G(B)}=G(A\cap B)$,
\item [(iv)] $\sum\limits_{j=1}^nZ_G(A_j)
=Z_G\left(\bigcup\limits_{j=1}^n A_j\right)$ if
$A_1,\dots,A_n$ are disjoint sets.

Also the following relation holds.

\item [(v)] $Z_G(A)=\overline{Z_G(-A)}$.

This follows from the relation
\item [(v$'$)] $I(f)=\overline{I(f_-)}$ for all $f\in L^2_G$, where
$f_-(x)=\overline{f(-x)}$.
\end{description}

Relations (i)---(iv) simply follow from the properties of the 
operator~$I$. The proof of~(v), more precisely of its strengthened
form~(v$'$) demands some more work. This property is needed to
decide when a random integral with respect to the complex 
valued random measure~$Z_G$ is a real valued random variable. 
(We shall soon define the random spectral measure $Z_G$ and 
the (random) integral with respect to it.) I describe the 
proof of~(v$'$) in the generalized field case. The proof in the 
discrete parameter case is similar, but simpler. Then I shall
give the proof also of~(iii). 

Relation (v$'$) can be simply checked if $f$ is a finite
trigonometrical polynomial in the discrete field case, or if
$f=\tilde\varphi$, $\varphi\in{\cal S}^c$,  in the generalized field
case. (In the case $f=\tilde\varphi$, $\varphi\in{\cal S}^c$, the
following argument works. Put
$f(x)=\tilde\varphi_1(x)+i\tilde\varphi_2(x)$ with
$\varphi_1,\varphi_2\in{\cal S}$. Then $I(f)=X(\varphi_1)+iX(\varphi_2)$,
and $f_-(x)=\bar{\tilde\varphi}_1(-x)-i\bar{\tilde\varphi}_2(-x)
=\tilde\varphi_1(x)+i(\widetilde{-\varphi_2}(x))$, hence
$I(f_-)=X(\varphi_1)+iX(-\varphi_2)=X(\varphi_1)-iX(\varphi_2)
=\overline{I(f)}$.)
Then a simple limiting procedure implies~(v$'$) in the general
case. 

Relation~(iii) follows from the identity
$$
EZ_G(A)\overline {Z_G(B)}=EI(\chi_A)\overline{I(\chi_B)}
=\int \chi_A(x)\overline{\chi_B(x)}G(\,dx)=G(A\cap B).
$$

\medskip
We have constructed with the help of a stationary Gaussian
random field with spectral measure~$G$ a set of complex valued
random variables $Z_G(\cdot)$ which satisfy properties (i)---(v).
In the next definition we shall call any class of sets of
complex valued random variables with these properties (independently
of how we have obtained them) a  random spectral measure.

\medskip\noindent
{\bf Definition of Random Spectral Measures.} 
{\it Let $G$ be a spectral measure. A set of random variables 
$Z_G(A)$, $G(A)<\infty$, satisfying (i)--(v) is called a 
(Gaussian) random spectral measure corresponding to the 
spectral measure~$G$.}

\medskip
Given a Gaussian random spectral measure $Z_G$ corresponding to a
spectral measure $G$ we define the (one-fold) stochastic integral
$\int f(x)Z_G(\,dx)$ for an appropriate class of functions~$f$.

Let us first consider simple functions of the form
$f(x)=\sum c_i\chi_{A_i}(x)$, where the sum is finite, the sets
$A_i$ are disjoint, and $G(A_i)<\infty$ for all indices~$i$. 
In this case we define
$$
\int f(x)Z_G(\,dx)=\sum c_iZ_G(A_i).
$$
We have to justify that the above formula is meaningful.
The problem is that the representation 
$f(x)=\sum c_i\chi_{A_i}(x)$ of a simple function is not
unique. We can write a set $A_i$ as the partition of
finitely many disjoint sets $A_{i,j}$, and if we replace $A_i$ 
with these sets $A_{i,j}$, and define the function to be $c_i$ 
on all sets $A_{i,j}$, then we get a different
representation of the same function~$f$. Now the additivity
property~(iv) guarantees that the the integral defined
by the new representation has the same value. It is not
difficult to check with the help of this observation
that the value of the integral of a simple function with
respect to $Z_G$ does not depend on the representation
of the simple function. Later we meet a generalized version
of this problem in the definition of multiple integrals
with respect to a spectral random measure. We shall explain
there the above argument in more detail.

Then we have
\begin{equation}
E\left|\int f(x)Z_G(\,dx)\right|^2=\sum c_i\bar c_jG(A_i\cap A_j)
=\int |f(x)|^2G(\,dx) \label{(3.5)}
\end{equation}
for all elementary functions.

Since the simple functions are dense in $L^2_G$, 
relation~(\ref{(3.5)}) enables us to define 
$\int f(x)Z_G(\,dx)$ for all $f\in L^2_G$ via
$L_2$-continuity. It can be seen that this integral
satisfies the identity
\begin{equation}
E\int f(x)Z_G(\,dx)\overline{\int g(x)Z_G(\,dx)}
=\int f(x)\overline{g(x)}G(\,dx) \label{(3.6a)}  
\end{equation}
for all pairs of functions $f,g\in L_G^2$. Moreover, similar
approximation with simple functions yields that
\begin{equation}
\int f(x)Z_G(\,dx)=\overline{\int\overline{f(-x)}Z_G(\,dx)}
\label{(3.6b)}
\end{equation}
for a function $f\in L^2_G$. Here we exploit the identity 
$Z_G(A)=\overline{Z_G(-A)}$ formulated in property~(v) of 
the random spectral measure $Z_G$.

Formula (\ref{(3.6b)}) implies in particular that $\int f(x)Z_G(\,dx)$
is real valued if $f(x)=\overline{f(-x)}$.

The last two identities together with the relations~(\ref{(3.1)}) 
and~(\ref{(3.2)}) imply that if we define the set of random 
variables $X_n$ and $X(\varphi)$ by means of the formulas
\begin{equation}
X_n=\int e^{i(n,x)}Z_G(\,dx), 
\quad n\in{\mathbb Z}_\nu, 
\label{(3.6)}
\end{equation}
and
\begin{equation}
X(\varphi)=\int\tilde\varphi(x) Z_G(\,dx), 
\quad \varphi\in{\cal S},
\label{($3.6'$)}
\end{equation}
where we integrate with respect to the random spectral 
measure $Z_G$, then we get a Gaussian stationary 
random discrete and generalized field with spectral 
measure~$G$, i.e. with correlation function given in 
formulas~(\ref{(3.1)}) and~(\ref{(3.2)}). To check this 
statement first we have to show that the random 
variables $X_n$ and $X(\varphi)$ defined in (\ref{(3.6)}) 
and~(\ref{($3.6'$)}) are real valued, or equivalently 
saying the identities $X_n=\overline{X_n}$ and 
$X(\varphi)=\overline{X(\varphi)}$ hold with 
probability~1. This follows from relation~(\ref{(3.6b)}) 
and the identities $e^{i(n,x)}=\overline{e^{(i(n,-x)}}$ 
and $\tilde\varphi(x)=\overline{\tilde\varphi(-x)}$ 
for a (real valued) function $\varphi\in{\cal S}$. 
Then we can calculate the correlation functions 
$EX_nX_m=EX_n\overline{X}_m$ and 
$EX(\varphi)X(\psi)=EX(\varphi)\overline{X(\psi)}$
by means of formula~(\ref{(3.6a)}), (\ref{(3.6)})
and~(\ref{($3.6'$)}).

We also have
$$
\int f(x)Z_G(\,dx)=I(f) \quad \textrm{for all } f\in L_G^2
$$
if we consider the previously defined mapping $I(f)$ with the
stationary random fields defined in~(\ref{(3.6)}) 
and~(\ref{($3.6'$)}). In particular, if we have a discrete
or generalized stationary random field $X_n$, $n\in{\mathbb Z}_\nu$,
or $X(\varphi)$, $\varphi\in{\cal S}$, and we construct
the random spectral measure $Z_G$ with the help of the
operator~$I$ in the way as we have written down at the beginning
of this subsection, then we can write
$$
X_n=I(e^{i(n,x})=\int e^{i(n,x)}Z_G(\,dx) \qquad
\textrm{for all }n\in{\mathbb Z}_\nu,
$$
and
$$
X(\varphi)=I(\tilde\varphi)=\int \tilde\varphi(x)Z_G(\,dx) \qquad
\textrm{for all }\varphi\in{\cal S}.
$$

It is not difficult to prove the subsequent theorem with the help
of the above results. I omit the details.

\medskip\noindent
{\bf Theorem 2.1.} {\it For a stationary Gaussian random field (a
discrete or generalized one) with a spectral measure $G$ there exists
a unique Gaussian random spectral measure $Z_G$ corresponding to the
spectral measure~$G$ on the same probability space as the Gaussian
random field such that relation~(\ref{(3.6)}) 
or~(\ref{($3.6'$)}) holds in the
discrete or generalized field case respectively.

Furthermore
\begin{equation}
{\cal B}(Z_G(A),\; G(A)<\infty)=\left\{
\begin{array}{l}
{\cal B}(X_n,\;n\in{\mathbb Z}_\nu) 
\textrm{ in the discrete field case,}\\
{\cal B}(X(\varphi),\;\varphi\in{\cal S})
\textrm{ in the generalized field case.}
\end{array} \right. \label{(3.7)}
\end{equation}
}

\medskip
Given a Gaussian stationary random field, a discrete field 
$X_n$, $n\in{\mathbb Z}_\nu$ or a generalized one 
$X(\varphi)$, $\varphi\in{\cal S}$ with some spectral
measure $G$ we call a random spectral measure $Z_G$
adapted to it if it satisfies relation~(\ref{(3.6)}) 
or~(\ref{($3.6'$)}).

\medskip
We have given a good representation of the random
variables in the fields ${\cal H}_1$ by means of
random integrals with respect to the random spectral
measure. Later we shall show that one can also define 
multiple (Wiener--It\^o) integrals with respect to
the random spectral measure. In such a way we can get
a good representation of all random variables with
finite second moment which are measurable with respect
to the $\sigma$-algebra generated by the elements of
the original Gaussian stationary random field. Moreover,
this representation is useful in the study of the limit
theorem problems we are interested in.

To work out the theory of multiple Wiener--It\^o integrals
it is useful to have some additional knowledge about the
properties of random spectral measures. I list below
these properties, and I show how to prove them.

\medskip
\begin{description}
\item[(vi)] The random variables $\textrm{Re}\, Z_G(A)$ are 
independent of the random variables $\textrm{Im}\, Z_G(A)$.
\item[(vii)] The random variables of the form $Z_G(A\cup(-A))$ 
are real valued. If the sets $A_1\cup(-A_1)$,\dots, 
$A_n\cup(-A_n)$ are disjoint, then the random variables 
$Z_G(A_1)$,\dots, $Z_G(A_n)$ are independent.
\item[(viii)] The relations 
$\textrm{Re}\, Z_G(-A)=\textrm{Re}\, Z_G(A)$ and 
$\textrm{Im}\, Z_G(-A)=-\textrm{Im}\, Z_G(A)$ hold, and if 
$A\cap(-A)=\emptyset$, then the (Gaussian) random variables 
$\textrm{Re}\, Z_G(A)$ and $\textrm{Im}\, Z_G(A)$ are 
independent with expectation zero and variance $G(A)/2$.
\end{description}

\medskip
These properties easily follow from (i)--(v). Since $Z_G(\cdot)$
are complex valued Gaussian random variables, to prove the above
formulated independence it is enough to show that the real and 
imaginary parts are uncorrelated. We show, as an example, 
the proof of~(vi).
\begin{eqnarray*}
E\textrm{Re}\, Z_G(A)\textrm{Im}\, Z_G(B)&=&\frac1{4i} 
E(Z_G(A)+\overline{Z_G(A)})
(Z_G(B)-\overline{Z_G(B)})\\
&=&\frac1{4i}E(Z_G(A)+Z_G(-A))(\overline{Z_G(-B)}-\overline{Z_G(B)})\\
&=&\frac1{4i}G(A\cap(-B))-\frac1{4i}G(A\cap B)\\
&&\qquad+\frac1{4i}G((-A)\cap(-B))-\frac1{4i}G((-A)\cap B)=0
\end{eqnarray*}
for all pairs of sets $A$ and $B$ such that $G(A)<\infty$, $G(B)<\infty$,
since $G(D)=G(-D)$ for all $D\in{\cal B}^\nu$. The fact that
$Z_G(A\cup(-A))$ is real valued random variable, and the relations
$\textrm{Re}\, Z_G(-A)=\textrm{Re}\, Z_G(A)$, 
$\textrm{Im}\, Z_G(-A)=-\textrm{Im}\, Z_G(A)$ under the conditions
of~(viii) follow directly from~(v). The remaining statements of~(vii)
and~(viii) can be proved similarly to~(vi) only the calculations are
simpler in this case.

The properties of the random spectral measure $Z_G$ listed above imply
in particular that the spectral measure~$G$ determines the joint
distribution of the corresponding random variables $Z_G(B)$,
$B\in{\cal B}^\nu$.

In the definition of random spectral measure we have imposed 
conditions (i)--(v) in the definition of random spectral measures, 
properties (vi)--(viii) were their consequences. Actually, we could 
have omitted also condition~(iv) from the definition, because
it can be deduced from the remaining conditions. This can be done 
by showing that the absolute value of the difference of the 
expressions at the two sides of the identity~(iv) have zero 
second moment. (See the corresponding Remark at page~18 of my 
Lecture note.) Let me remark that if we want to define the 
distribution of a set of jointly Gaussian random variables,
then it is enough to give the expected value and correlation 
function of these random variables. We followed a similar 
approach in the formulation of properties (i)--(iii) in the 
definition of random spectral measures. But since we work 
here with complex valued random variables, we had to  
add property~(v) to these conditions to get a definition which
determines the distribution of the random spectral measures.

\subsection{An application of the results on random spectral measures}

We can define generalized stationary Gaussian field on the space 
$\cal S$ with such spectral measures $G$ which satisfy (\ref{(3.3)}), 
because this guarantees that the integral (\ref{(3.2)}) is meaningful 
for all pairs of test-functions $\varphi,\psi\in\cal S$. Let us enlarge 
the space of test-functions $\cal S$ to a larger linear linear space 
${\cal T}\supset{\cal S}$ which has the property that it is invariant 
under the shift transformations, i,e. if $\varphi(x)\in{\cal T}$, then 
$\varphi(x-t)\in{\cal T}$ for all arguments~$t$. If the integral 
(\ref{(3.2)}) is meaningful for all pairs of functions 
$\varphi,\psi\in{\cal T}$, then we can define a Gaussian stationary 
random field $X(\varphi)$, on a larger class of test functions 
$\varphi\in{\cal T}$, i.e. there is a Gaussian random field 
$X(\varphi)$ with test functions $\varphi\in{\cal T}$ which  
satisfies the identity
$X(a_1\varphi_1+a_2\varphi_2)=a_1X(\varphi_1)+a_2X(\varphi_2)$ with
probability 1 for all constants $a_1,a_2$ and functions 
$\varphi_1,\varphi_2\in\cal T$, and it has expectation 
$EX(\varphi)=0$ and correlation function satisfying~(\ref{(3.2)}).

An especially interesting case appears when the linear space
contains the indicator function of all rectangles of the form
$\prod_{j=1}^\nu[a_j,b_j)$. In this case we get, restricting the stationary
random field $X(\varphi)$, $\varphi\in{\cal T}$, to the set of cubes
$\prod_{j=1}^\nu[a_j-\frac12,a_j+\frac12)$, and identifying this cube
with its center point $(a_1,\dots,a_\nu)$ a discrete Gaussian
stationary random field which can be considered the discretization of
the original generalized field $X(\varphi)$, $\varphi\in\cal S$. 

The above construction is especially interesting in the case when the
generalized random field has the spectral measure $G$ with density function
$|X|^{-\alpha}$, with a parameter $\alpha$ is chosen so, that the relation
(\ref{(3.3)}) holds with it. Such a spectral measure defines a self-similar
stationary generalized Gaussian random field, and the above sketched 
method helps us to construct a discrete Gaussian stationary random field. 
I explain how we can construct the so-called fractional Brownian motions 
in such a way.

A fractional Brownian motion with Hurst parameter $H$, $0<H<1$, 
is defined as a Gaussian process $X(t)$, $t\ge0$, with 
continuous trajectories and zero expectation, i.e. $EX(t)=0$ 
for all $t\ge0$, and with correlation function 
$R_H(s,t)=EX(s)X(t)=\frac12(s^{2H}+t^{2H}-|t-s|^{2H})$ 
for all $0\le s,t<\infty$.
Naturally we must prove that such a process really exists.
 
I briefly explain that the  correlation function of a fractional 
Brownian motion  has a natural representation as the 
correlation function of the discretized version of an 
appropriately defined Gaussian stationary generalized 
self-similar field. 

To understand this approach observe that a fractional Brownian 
motion with Hurst parameter~$H$ has the self-similarity property 
$EX(as)X(at)=a^{2H}EX(s)X(t)$ for all $a>0$, and simple
calculation shows that it also has the following stationary
increments property: $EX(0)^2=0$, hence $X(0)=0$ with probability~1, 
and $E[X(s+u)-X(u)][X(t+u)-X(u)]=EX(s)X(t)$
for all $0\le s,t,u<\infty$. To construct a fractional Brownian 
motion $X(t)$ first we define an appropriate stationary, 
Gaussian generalized self-similar random field $\bar X(\varphi)$, 
$\varphi\in{\cal S}_1$ in the space of the real valued functions 
of the Schwartz space, and then we extend this field it to 
a larger parameter set (of functions), containing the indicator 
functions $\chi_{[0,t]}$ of the intervals $[0,t]$ for all $t\ge0$.
Finally we define the process $X(t)$ as $X(t)=\bar X(\chi_{[0,t]})$.

More explicitly, we can define for a parameter $\alpha$ a
stationary generalized Gaussian field $\bar X(\varphi)$, 
$\varphi\in{\cal S}^1$, with zero expectation and spectral 
density $|u|^{-2\alpha}$, i.e. put 
$E\bar X(\varphi)\bar X(\psi)=\int 
\tilde\varphi(u)\bar{\tilde\varphi}(u)|u|^{-2\alpha}\,du$.
Then we introduce its natural extension to a function 
space containing the functions $\chi_{[0,t]}$ for all $t>0$. 
Then we have
$$
E\bar X(\chi_{[0,s]})\bar X(\chi_{[0,t]})=
\int \tilde\chi_{[0,s]}(u)\bar{\tilde \chi}_{[0,t]}(u)
|u|^{-2\alpha}\,du
=\int \frac{e^{isu}-1}{iu}\frac{e^{-itu}-1}{-iu}
|u|^{-2\alpha}\,du,
$$
provided that these integrals are convergent.  

The above defined generalized fields exist if $2\alpha>-1$, 
and their extension to an appropriate space $\cal T$
containing the indicator functions $\chi_{[0,t]}$ exists if 
$-1<2\alpha<1$. The first condition is needed to guarantee 
that the singularity of the integrand in the formula 
expressing the correlation function is not too large 
in the origin, and the second condition is needed to 
guarantee that the singularity of this integrand is 
not too large at the infinity even if we work with
the Fourier transform of the indicator functions
$\chi_{[0,t]}$ in the discretized case.

Simple calculation shows that the correlation function of 
the above defined random field satisfies the identity 
$E\bar X(\varphi_a)\bar X(\psi_a)=a^{-(1+2\alpha)}
E\bar X(\varphi)\bar X(\psi)$, with the functions
$\varphi_a(x)=\varphi(ax)$, $\psi_a(x)=\psi(ax)$, and 
similarly, we have $EX(as)X(at)=a^{(1+2\alpha)}EX(s)X(t)$ 
for  all $a>0$. Besides, the Gaussian stochastic process 
$X(t)$, $t>0$, has stationary increments, i.e. 
$E[X(s+u)-X(u)][X(t+u)-X(u)]=EX(s)X(t)$ for all 
$0\le s,t,u<\infty$, and $EX(0)^2=0$. This follows from 
its construction with the help of a stationary Gaussian 
random field. 

The above calculations imply that with the choice 
$\alpha=H-1/2$ we get the correlation function of a 
fractional Brownian motion with Hurst parameter $H$ 
for all $0<H<1$, more precisely the correlation 
function of this process multiplied by an appropriate 
constant. Indeed, it follows from the stationary 
increments property of the process that 
$E(X(t)-X(s))^2=EX(t-s)^2$, if $t\ge s$, and the 
self-similarity property of this process implies that 
$EX(s)X(t)=\frac12[EX(s)^2+EX(t)^2-E(X(t)-X(s))^2]
=\frac12 EX(1)^2[s^{2H}+t^{2H}-|t-s|^{2H}]$.

We can get a representation of this process by means 
of a random  integral with respect to a random 
spectral measure. This representation has the form
$$
X(t)=\int\frac {e^{itu}-1}{iu}|u|^{-H+1/2}Z(\,du), \quad t>0,
$$
with the random spectral measure $Z(\cdot)$ corresponding to
the Lebesgue measure on the real line. Here I omit the proof 
that such a stochastic process also has a version with 
continuous trajectories.

The representation of the fractional Brownian processes
may be useful in the study of this process, but it seems to me
that this was not fully exploited in the research about this 
subject. Finally I remark that there is a rather complete 
description of the stationary, self-similar Gaussian processes 
(both discrete and generalized ones) together with their Gaussian 
domain of attraction. (See P. Major: On renormalizing Gaussian 
fields. {\it Z.~Wahrscheinlichkeitstheorie verw. Gebiete\/} 
{\bf 59} (1982), 515--533.) 

Here I explain the content of this paper in a very informal way. 
The self-similar stationary Gaussian random fields (with 
random elements of zero expectation) are those, whose 
spectral measure are homogeneous functions. The stationary 
Gaussian random fields in their domain of attraction are 
those random fields, whose spectral measures are close (in a 
natural sense) to the spectral measure of the limiting field.

\Section{Multiple Wiener--It\^o integrals}

Let us take a Gaussian stationary random field (either a
discrete field $X_n$, $n\in{\mathbb Z}_\nu$, or a generalized
field $X(\varphi)$, $\varphi\in{\cal S}$). We considered in both
cases the (real) Hilbert space $\cal H$ consisting of the 
square integrable random variables, measurable with respect
to the $\sigma$-algebra generated by the random variables
of the stationary random field with the usual scalar product
$(\xi,\eta)=E\xi\eta$, and we defined on it the group
of shift transformations $T_n$, $n\in{\mathbb Z}_\nu$ and $T_t$,  
$t\in R^\nu$, for the discrete and generalized random fields 
respectively. We want to get a good representation of this
Hilbert space together with the shift transformations on it.

First we decompose the Hilbert space $\cal H$ into the direct
sum of orthogonal subspaces which are invariant with respect
to all shift transformations.

We construct the invariant subspaces of the Hilbert space 
$\cal H$ with the help of its subspace ${\cal H}_1$ which is 
the closure of the finite linear combinations 
$\sum_{j=1}^kc_j X_{n_j}$ of the elements $X_n$,
$n\in{\mathbb Z}_\nu$, in the case of discrete stationary
random fields and of the random variables $X(\varphi)$, 
$\varphi\in\cal S$, in the case of generalized stationary
random fields, where the closure is taken in the Hilbert space
$\cal H$. First we define for all $n=1,2,\dots$ the  Hilbert 
subspace ${\cal H}_{\le n}\subset{\cal H}$, $n=1,2,\dots$, as 
the subspace which is the closure of the linear space 
consisting of the elements $P_n(X_{t_1},\dots,X_{t_m})$, where 
$P_n$ runs through all polynomials of degree less than or 
equal to~$n$, the integer~$m$ is arbitrary, and 
$X_{t_1},\dots,X_{t_m}$ are elements of ${\cal H}_1$. Let 
${\cal H}_0={\cal H}_{\le 0}$ consist of the constant functions, 
and put ${\cal H}_n={\cal H}_{\le n}\ominus{\cal H}_{\le n-1}$, 
$n=1,2,\dots$, where $\ominus$ denotes orthogonal completion. It 
is clear that the Hilbert space ${\cal H}_1$ given in this 
definition agrees with the previously defined Hilbert space 
${\cal H}_1$. 

The following theorem holds.

\medskip\noindent
{\bf Theorem 3.1 on the decomposition of the Hilbert space
${\cal H}$ consisting of the square integrable random variables
measurable with respect to a stationary random field.} {\it The 
Hilbert space ${\cal H}$ has the following decomposition with the 
help of the above defined Hilbert spaces ${\cal H}_n$, $n=0,1,2,\dots$.
\begin{equation}
{\cal H}={\cal H}_0+{\cal H}_1+{\cal H}_2+\cdots,  \label{(2.1)}
\end{equation}
where $+$ denotes direct sum. Besides, all subspaces
${\cal H}_n$ of this decomposition are invariant subspaces
of all shift transformations of the Hilbert space ${\cal H}$.}

\medskip
I explain the main ideas of the proof together with the 
formulation of some basic, classical results of the analysis 
that we need in the proof. First I recall the definition of 
Hermite polynomials, which play an important role both in 
this proof and in some further consideration. Then I also 
formulate some results about their properties.  

\medskip\noindent
{\bf Definition of Hermite polynomials.} {\it The $n$-th 
Hermite polynomial $H_n(x)$ with leading coefficient~1 is 
the polynomial of order~$n$ defined by the formula
$H_n(x)=(-1)^ne^{x^2/2}\frac{d^n}{dx^n}(e^{-x^2/2})$.} 

\medskip
The Hermite polynomials have the following property.

\medskip\noindent
{\bf Theorem 3A.} {\it The Hermite polynomials $H_n(x)$, $n=0,1,2,\dots$,
form a complete orthogonal system in
$L_2\left(R,{\cal B},\frac1{\sqrt{2\pi}}e^{-x^2/2}\,dx\right)$.
(Here ${\cal B}$ denotes the Borel $\sigma$-algebra on the real line.)}

\medskip
We also need the following measure theoretical results in
the proof of Theorem~3.1.

Let  $(X_j,{\cal X}_j,\mu_j)$, $j=1,2,\dots$, be countably many 
independent copies of a probability space $(X,{\cal X},\mu)$. 
Let $(X^\infty,{\cal X}^\infty,\mu^\infty)
=\prod\limits_{j=1}^\infty(X_j,{\cal X}_j,\mu_j)$. With such a 
notation the following result holds.

\medskip\noindent
{\bf Theorem 3B.} {\it Let $\varphi_0,\varphi_1,\dots$,
$\varphi_0(x)\equiv1$, be a complete orthonormal system in 
a Hilbert space $L_2(X,{\cal X},\mu)$. Then the functions
$\prod\limits_{j=1}^\infty\varphi_{k_j}(x_j)$, where only 
finitely many indices $k_j$ differ from~0, form a complete 
orthonormal basis in the product space 
$L_2(X^\infty,{\cal X}^\infty,\mu^\infty)$.}

\medskip\noindent
{\bf Theorem 3C.} {\it Let $Y_1,Y_2,\dots$ be random variables 
on a probability space $(\Omega,{\cal A},P)$ taking values in 
a measurable space  $(X,{\cal X})$. Let $\xi$ be a real valued 
random variable measurable with respect to the $\sigma$-algebra 
${\cal B}(Y_1,Y_2,\dots)$, and let $(X^\infty,{\cal X}^\infty)$ 
denote the infinite product
$(X\times X\times\cdots,{\cal X}\times {\cal X}\times\cdots)$ 
of the space $(X,{\cal X})$ with itself. Then there exists a 
real valued, measurable function~$f$ on the space 
$(X^\infty,{\cal X}^\infty)$ such that $\xi=f(Y_1,Y_2,\dots)$.}

\medskip\noindent
Theorem 3.1 can be proved with the help of Theorems~3.A, 3B,
and~3C in a natural way. One can show with the help of Theorems~3A
and~3B that if we introduce the infinite product $\mu_\infty$ of 
the standard Gaussian probability distribution with itself on the 
infinite product space $(R^\infty,{\cal B}^\infty)$, then 
the set of all finite products of the form 
$H_{j_1}(x_{k_1})\cdots H_{j_l}(x_{k_l})$ provide an
orthogonal basis in the Hilbert space 
$L_2(R^\infty,{\cal B}^\infty,\mu_\infty)$. Then we can choose
an orthonormal basis of standard Gaussian random variables 
$X_1,X_2,\dots$ in ${\cal H}_1$, and we can construct with
the help of this basis and Theorem~3C an appropriate embedding 
that implies Theorem~3.1. I omit the details of the proof. 

I also formulate a result in the next Corollary~3.2
which we can get as a by-product of this proof. This result
will be useful in our later considerations.

\medskip\noindent
{\bf Corollary 3.2.} {\it Let $\xi_1,\xi_2,\dots$ be an 
orthonormal basis in ${\cal H}_1$, and let $H_j(x)$ denote the
Hermite polynomial with  order~$j$ and leading coefficient~1.
Then the random variables $H_{j_1}(\xi_1)\cdots H_{j_k}(\xi_k)$, 
$k=1,2,\dots$, $j_1+\cdots+j_k=n$, and $j_k>0$ form a complete 
orthogonal basis in ${\cal H}_n$.}

\medskip
In a more detailed discussion we could have introduced the natural
multivariate version of the Hermite polynomials, the so-called
Wick polynomials. So I did in the lecture note which is the 
basis of this lecture. But here I omitted the discussion of Wick 
polynomials, because in the present text I do not work with them.

We constructed a good decomposition of the Hilbert space
$\cal H$ into orthogonal invariant subspaces. Next we shall
present the elements of these subspaces in the form of multiple 
Wiener--It\^o integrals, because this is useful in the study 
of the limit problems we are interested in.

\subsection{The construction of multiple Wiener--It\^o integrals}

The multiple random integral I shall discuss here is actually
different of the original Wiener--It\^o integral. This is a 
version of it which was introduced by R.~L.~Dobrushin. 
Nevertheless, I shall apply the original name. 

The original Wiener--It\^o integral is taken with respect to 
a Gaussian orthogonal random measure $Z_\mu$ corresponding 
to a measure~$\mu$. This Gaussian orthogonal random 
measure $Z_\mu$ is defined on a measure space $(M,{\cal M},\mu)$ 
with some $\sigma$-finite measure~$\mu$ on $\cal M$, and it 
consists of (jointly) Gaussian random variables  $Z_\mu(A)$ 
defined for all sets $A\in\cal M$ such that $\mu(A)<\infty$, 
and it has the properties that $EZ_\mu(A)=0$, 
$EZ_\mu(A)^2=\mu(A)$ for all $A\in{\cal M}$, the random 
variables $Z_\mu(A_j)$ are independent for disjoint sets 
$A_j$, and $Z_\mu(\cdot)$ is additive in the following sense. 
If $A_1,\dots,A_k$, are disjoint sets, then 
$Z_\mu\left(\bigcup\limits_{j=1}^k A_j\right)=\sum\limits_{j=1}^k Z_\mu(A_j)$.
It\^o defined the $k$-fold random integral integral with respect 
to a Gaussian random measure $Z_\mu$ for all multiplicities
$k=0,1,2,\dots$. He defined the $k$-fold integral for such 
functions $f(x_1,\dots,x_k)$ which are in the Hilbert space
$L_2((M^k,{\cal M}^k,\mu^k))$. He could express all random
variables which have finite second moment and are measurable
with respect to the $\sigma$-algebra generated by the
random variables $Z_\mu(\cdot)$ of the Gaussian orthogonal
random measure as a sum of random integrals of different 
multiplicity. Moreover, this representation is unique.
This result turned out to be useful in certain investigations.

Dobrushin worked out the theory of an analogous multiple
random integral, where we integrate with respect to a
random spectral measure instead of a Gaussian orthogonal
random measure. The proof of the results in this case is
somewhat more complicated. On the other hand, this integral
is more appropriate in our investigations of limit theorems
for non-linear functionals of Gaussian stationary random
fields. First I briefly explain why such an integral is 
useful for us, and what kind of technical difficulties  
have to be overcome in their study which do not appear in 
the theory of the original Wiener--It\^o integrals.

We want to study the Hilbert space $\cal H$ determined
by our stationary random field together with the shift
transformations on it. We can handle the shift transformation
better with the help of Fourier transforms. To understand
this let us take the following simple example. Let us
consider a function $f(x)$ on the real line together
with its shift $T_t f(x)=f(x-t)$. If we work with this 
shift transformation, then we can better calculate
with the Fourier transform $\tilde f(u)$ since 
$\widetilde{T_t f(u)}=e^{itu}\tilde f(u)$. (Those who are
familiar with the spectral theory of operators in Hilbert
spaces can give the following interpretation to this
example. If we take the Hilbert space of square
integrable functions with respect to the Lebesgue
measure, then the shift transformation 
$T_t\colon f(x)\to f(x-t)$ is a unitary operator on it.
In the above formula we gave the spectral representation
of this operator with the help of Fourier transforms.)
Integration with respect to the random spectral measure
plays a role similar to the Fourier transform on the
real line, and as a consequence we shall have a 
representation of the shift transformation on $\cal H$
which is similar to the above example.

In the definition of the original Wiener--It\^o integrals
it was exploited that the Gaussian orthogonal random 
measure $Z_\mu(A_j)$ of disjoint sets $A_j$ are independent. 
We want to apply a similar argument in the definition of 
random integrals with respect to random spectral measures. 
But here we have only a weaker independence property. We 
can state that $Z_G(A)$ and $Z_G(B)$ are independent
only if $A\cup(-A)$ and $B\cup(-B)$ are disjoint. (See
property~(vii) of the random spectral measures.) Another
point where we have to be careful is that we want to define 
the random integrals so that they are real valued, since 
$\cal H$ contains real valued random variables. Since the 
random spectral measure $Z_G(\cdot)$ is complex valued, we 
have to find the appropriate class of kernel functions to
guarantee this property of the random integrals.

\medskip
Now I turn to the definition of the multiple Wiener--It\^o
integrals. First I introduce the( {\it real}) Hilbert space 
$\bar{{\cal H}}_G^n$ and its symmetrization ${\cal H}_G^n$ 
whose elements will be the kernel functions of the 
$n$-fold Wiener--It\^o integrals with respect to a random
spectral measure~$Z_G$. As we shall later see, the 
Wiener--It\^o integrals of a function and of its
symmetrization agree. Hence it would be enough to work 
with Wiener--It\^o integrals whose kernel functions are 
in the symmetrized space ${\cal H}_G^n$. But for some 
technical reasons it will be better to work with 
Wiener--It\^o integrals with kernel functions from both spaces.

Let $G$ be the spectral measure of a stationary Gaussian 
random field (discrete or generalized one). We define the 
following {\it real}\/ Hilbert spaces 
$\bar{{\cal H}}_G^n$ and ${\cal H}_G^n$,
$n=1,2,\dots$. 

We have $f_n\in\bar{{\cal H}}_G^n$ if and only if
$f_n=f_n(x_1,\dots,x_n)$, \ $x_j\in R^\nu$, $j=1,2,\dots,n$, is a
complex valued function of $n$ variables, and

\medskip
\begin{description}
\item[(a)] $f_n(-x_1,\dots,-x_n)=\overline{f_n(x_1,\dots,x_n)}$,
\item[(b)]
$\|f_n\|^2=\int|f_n(x_1,\dots,x_n)|^2G(\,dx_1)\dots G(\,dx_n)<\infty$.
\end{description}

\medskip
Relation~(b) also defines the norm in $\bar{{\cal H}}^n_G$. The
subspace ${\cal H}^n_G\subset\bar{{\cal H}}_G^n$ contains those functions
$f_n\in\bar{{\cal H}}_G^n$ which are invariant under permutations of
their arguments, i.e.

\medskip
\begin{description}
\item[(c)] $f_n(x_{\pi(1)},\dots,x_{\pi(n)}))=f_n(x_1,\dots,x_n)$
for all $\pi\in\Pi_n$, where $\Pi_n$ denotes the group of all
permutations of the set $\{1,2,\dots,n\}$.
\end{description}

\medskip
The norm in ${\cal H}_G^n$ is defined in the same way as in
$\bar{{\cal H}}_G^n$. Moreover, the scalar product is also similarly
defined, namely if $f,\,g\in\bar{{\cal H}}_G^n$, then
\begin{eqnarray*}
(f,g)&=&\int f(x_1,\dots,x_n)\overline{g(x_1,\dots,x_n)}
G(\,dx_1)\dots G(\,dx_n)\\
&=&\int f(x_1,\dots,x_n)g(-x_1,\dots,-x_n)G(\,dx_1)\dots G(\,dx_n).
\end{eqnarray*}
Because of the symmetry $G(A)=G(-A)$ of the spectral measure
$(f,g)=\overline{(f,g)}$, i.e. the scalar product $(f,g)$ is a real
number for all $f,\,g\in\bar{{\cal H}}_G^n$. This means that
$\bar{{\cal H}}_G^n$ is a real Hilbert space.
We also define ${\cal H}_G^0=\bar{{\cal H}}_G^0$ as
the space of real constants with the norm $\|c\|=|c|$.  I remark
that $\bar{{\cal H}}_G^n$ is actually the $n$-fold direct product of
${\cal H}_G^1$, while ${\cal H}_G^n$ is the $n$-fold symmetrical direct
product of ${\cal H}^1_G$. Condition~(a) means heuristically that
$f_n$ is the Fourier transform of a real valued function.

We also define the so-called Fock space $\textrm{Exp\,}{\cal H}_G$
whose elements are sequences of functions $f=(f_0,f_1,\dots)$, 
$f_n\in{\cal H}_G^n$ for all $n=0,1,2,\dots$, such that
$$
\|f\|^2=\sum_{n=0}^\infty \frac1{n!}\|f_n\|^2<\infty,
$$
where $\|f_n\|$ denotes the norm of the function $f_n$ in the
Hilbert space ${\cal H}^n_G$.

Given a function $f\in\bar{{\cal H}}^n_G$ we define 
$\textrm{Sym}\, f$ as
$$
\textrm{Sym}\, f(x_1,\dots,x_n)=\frac1{n!}\sum_{\pi\in\Pi_n}
f(x_{\pi(1)},\dots,x_{\pi(n)}).
$$
Clearly, $\textrm{Sym}\, f\in{\cal H}_G^n$, and
\begin{equation}
\|\textrm{Sym}\, f\|\le \|f\|. \label{(4.1)}
\end{equation}

Let $Z_G$ be a Gaussian random spectral measure corresponding 
to the spectral measure~$G$ on a probability space 
$(\Omega,{\cal A},P)$. We shall define the $n$-fold 
Wiener--It\^o integrals
$$
I_G(f_n)=\frac1{n!}\int f_n(x_1,\dots,x_n)Z_G(\,dx_1)\dots Z_G(\,dx_n)
\quad \textrm{for } f_n\in\bar{{\cal H}}_G^n
$$
and
$$
I_G(f)=\sum_{n=0}^\infty I_G(f_n)\quad \textrm{for }
f=(f_0,f_1,\dots)\in\textrm{Exp}\,{\cal H}_G.
$$
We shall see that  $I_G(f_n)=I_G(\textrm{Sym}\, f_n)$ for all
$f_n\in\bar{{\cal H}}_G^n$. Therefore, it would have been 
sufficient to define the Wiener--It\^o integral only for 
functions in ${\cal H}_G^n$. Nevertheless, some arguments 
become simpler if we work in $\bar{{\cal H}}_G^n$. In the 
definition of Wiener--It\^o integrals we restrict 
ourselves to the case when the spectral measure is
non-atomic, i.e. $G(\{x\})=0$ for all $x\in R^\nu$. This 
condition is satisfied in all interesting cases. We could
extend the definition of Wiener--it\^o integrals also 
to the case when $G$ may be non-atomic, but we do not
do that, because it seems so that we would not gain 
very much with such an extension.   

In the definition of multiply Wiener--It\^o integrals
we follow a similar approach as in the definition of
the one-fold integrals with respect to a random spectral
measure. First we define them to a class of appropriately
defined simple functions, and then we show that this
integral can be extended because of an $L_2$-contraction
property of this integral to the whole space 
$\bar{{\cal H}}_G^n$.

In the definition of the simple functions we have to
take into account that they are elements of the space
$\bar{{\cal H}}_G^n$. It will be natural to define them
together with the notion of regular systems which are 
collections of disjoint subsets of $R^\nu$ with some
additional properties. The simple functions of 
$n$-variables are those functions of $\bar{{\cal H}}_G^n$
which are adapted in an appropriate way to a regular system. 
Here I give the definition of these notions.

\medskip\noindent
{\bf Definition of Regular Systems and the Class of Simple 
Functions.} {\it Let 
$$
{\cal D}=\{\Delta_j,\;j=\pm1,\pm2,\dots,\pm N\}
$$ 
be a finite collection of bounded, measurable sets in $R^\nu$  
indexed by the integers $\pm1,\dots,\pm N$. We say that 
${\cal D}$ is a regular system if $\Delta_j=-\Delta_{-j}$, and 
$\Delta_j\cap\Delta_l=\emptyset$ if $j\neq l$ for all 
$j,l=\pm1,\pm2,\dots,\pm N$. 
A function $f\in\bar{{\cal H}}_G^n$ is adapted to this system 
${\cal D}$ if $f(x_1,\dots,x_n)$ is constant on the sets
$\Delta_{j_1}\times\Delta_{j_2}\times\cdots\times\Delta_{j_n}$, \
$j_l=\pm1,\dots,\pm N$, $l=1,2,\dots,n$,  it vanishes outside 
these sets and also on those sets of the form
$\Delta_{j_1}\times\Delta_{j_2}\times\cdots\times\Delta_{j_n}$, 
for which $j_l=\pm j_{l'}$ for some $l\neq l'$. 

A function $f\in\bar{{\cal H}}_G^n$ is in the class
$\hat{\bar{{\cal H}}}_G^n$ of simple functions, 
and a (symmetric) function $f\in{\cal H}_G^n$ is in 
the class $\hat{{\cal H}}_G^n$ of simple symmetric 
functions if it is adapted to some regular system 
${\cal D}=\{\Delta_j,\;j=\pm1,\dots,\pm N\}$.} 

\medskip
Next we define the Wiener--It\^o integrals of simple functions.

\medskip\noindent
{\bf Definition of Wiener--It\^o Integral of Simple Functions.}
{\it Let a simple function $f\in\hat{\bar{{\cal H}}}_G^n$ be adapted 
to some regular systems ${\cal D}=\{\Delta_j,\;j=\pm1,\dots,\pm N\}$. 
Its Wiener--It\^o integral with respect to the random spectral
measure $Z_G$ is defined as
\begin{eqnarray}
&&\int f(x_1,\dots,x_n)Z_G(\,dx_1)\dots Z_G(\,dx_n) 
\label{(4.2)} \\
&&\qquad =n!I_G(f)=\sum_{\substack{j_l=\pm1,\dots,\pm N\\ l=1,2,\dots,n}}
f(x_{j_1},\dots,x_{j_n})Z_G(\Delta_{j_1})\cdots Z_G(\Delta_{j_n})
\nonumber,
\end{eqnarray}
where $x_{j_l}\in\Delta_{j_l}$, $j_l=\pm1,\dots,\pm N$, $l=1,\dots,n$.}
\medskip

I remark that although the regular system ${\cal D}$ to which 
$f$ is adapted is not uniquely determined (e.g. the elements of 
${\cal D}$ can be divided to smaller sets in an appropriate way), 
the integral defined in~(\ref{(4.2)}) is meaningful, i.e. it 
does not depend on the choice of ${\cal D}$. This can be seen 
by observing that a refinement of a regular system ${\cal D}$ 
to which the function $f$ is adapted yields the same value 
for the sum defining $n!I_G(f)$ in formula~(\ref{(4.2)}) as 
the original one. (We also exploit that if a function~$f$ 
is adapted to two different regular systems ${\cal D}_1$
and ${\cal D}_2$, then there is a regular system ${\cal D}_3$
which is a refinement of both of them, and the function~$f$
is adapted to it.) This follows from the additivity of the 
random spectral measure $Z_G$ formulated in its 
property~(iv), since this implies that each term
$f(x_{j_1},\dots,x_{j_n})Z_G(\Delta_{j_1})\cdots Z_G(\Delta_{j_n})$
in the sum at the right-hand side of formula~(\ref{(4.2)}) 
corresponding to the original regular system equals the sum 
of all such terms $f(x_{j_1},\dots,x_{j_n})
Z_G(\Delta'_{j'_1})\cdots Z_G(\Delta'_{j'_n})$ in the sum
corresponding to the refined partition for which
$\Delta'_{j'_1}\times\cdots\times\Delta'_{j'_n}\subset
\Delta_{j_1}\times\cdots\times\Delta_{j_n}$.

By property~(vii) of the random spectral measures all products
$$
Z_G(\Delta_{j_1})\cdots Z_G(\Delta_{j_n})
$$ 
with non-zero coefficient in~(\ref{(4.2)}) are products of 
independent random variables. We had this property in mind 
when requiring the condition that the function $f$ vanishes on 
a product $\Delta_{j_1}\times\cdots\times\Delta_{j_n}$ if 
$j_l=\pm j_{l'}$ for some $l\neq l'$. This condition is 
interpreted in the literature as discarding the hyperplanes 
$x_l=x_{l'}$ and $x_l=-x_{l'}$, \ $l,l'=1,2,\dots,n$, 
$l\neq l'$, from the domain of integration. (Let me remark
that here we also omitted  the hyperplanes $x_l=-x_{l'}$ and 
not only the hyperplanes $x_l=x_{l'}$, $l\neq l'$, from the domain of 
integration. In particular, we omitted the points $x_l=-x_l$, i.e.
$x_l=0$ for all $1\le l\le d$. This is a difference from the 
definition of the original Wiener--It\^o integrals with respect 
to a Gaussian orthogonal random measure, where we omit only the
hyperplanes $x_l=x_{l'}$, $l\neq l'$, from the domain of 
integration.)
Property~(a) of the functions in $\bar{{\cal H}}_G^n$ 
and property~(v) of the random spectral measures imply that
$I_G(f)=\overline{I_G(f)}$, i.e. $I_G(f)$ is a real valued 
random variable for all $f\in\hat{\bar{{\cal H}}}_G^n$. 
The relation
\begin{equation}
EI_G(f)=0, \quad \textrm{for }f\in\hat{\bar{{\cal H}}}_G^n,
\quad n=1,2,\dots \label{(4.3)}
\end{equation}
also holds. Let 
$\hat{{\cal H}}_G^n={\cal H}_G^n\cap\hat{\bar{{\cal H}}}_G^n$.
If $f\in\hat{\bar{{\cal H}}}_G^n$, then  
$\textrm{Sym}\, f\in\hat{{\cal H}}_G^n$, and
\begin{equation}
I_G(f)=I_G(\textrm{Sym}\, f). \label{(4.4)}
\end{equation}
Relation~(\ref{(4.4)}) holds, since 
$Z_G(\Delta_{j_1})\cdots Z_G(\Delta_{j_n})=Z_G(\Delta_{\pi(j_1)})\cdots
Z_G(\Delta_{\pi(j_n)})$ for all permutations $\pi\in\Pi_n$.
I also claim that
\begin{equation}
EI_G(f)^2\le\frac1{n!}\|f\|^2 \quad\textrm {for \ }
f\in\hat{\bar{{\cal H}}}_G^n, \label{(4.5)}
\end{equation}
and
\begin{equation}
EI_G(f)^2=\frac1{n!}\|f\|^2 \quad\textrm {for \ } f\in\hat{{\cal H}}_G^n.
\label{($4.5'$)}
\end{equation}
More generally, I claim that
\begin{eqnarray}
EI_G(f)I_G(h)=\frac1{n!}(f,g)&=&\frac1{n!}\int f(x_1,\dots,x_n)
\overline{g(x_1,\dots,x_n)}G(\,dx_1)\dots G(\,dx_n) \nonumber \\ 
&&\qquad\qquad\qquad\qquad \textrm {for \ } f,g\in\hat{{\cal H}}_G^n.
\label{($4.5''$)}
\end{eqnarray}

Because of~(\ref{(4.1)}) and~(\ref{(4.4)}) it is enough to 
check~(\ref{($4.5''$)}).

Let ${\cal D}$ be a regular system of sets in $R^\nu$,
and choose some sets $\Delta_{j_l}\in{\cal D}$ and 
$\Delta_{k_l}\in{\cal D}$ with indices $j_1,\dots,j_n$ 
and $k_1,\dots,k_n$ such that
$j_l\neq\pm j_{l'}$, $k_l\neq\pm k_{l'}$ if $l\neq l'$. To 
prove~(\ref{($4.5''$)}) I show that
\begin{eqnarray*}
&&EZ_G(\Delta_{j_1})\cdots Z_G(\Delta_{j_n})
\overline{Z_G(\Delta_{k_1})\cdots Z_G(\Delta_{k_n})}\\
&&\qquad\qquad=\left\{
\begin{array}{l}
G(\Delta_{j_1})\cdots G(\Delta_{j_n}) \quad\textrm{ if \ }
\{j_1,\dots,j_n\}=\{k_1,\dots,k_n\}, \\
0 \quad \textrm{otherwise.}
\end{array} \right.
\end{eqnarray*}

To see the second relation in the last formula we will
decompose the product whose expectation is investigated
at the left-hand side of this formula to the product of
two independent components in such a way that one of them
has zero expectation.We shall find such a decomposition
with the help of property~(vii) of the random spectral 
measures.

The identity in the second relation of the last formula 
has to be proved under the condition 
$\{j_1,\dots,j_n\}\neq\{k_1,\dots,k_n\}$. Hence in this case 
there is an index~$l$ such that either $j_l\neq\pm k_{l'}$ 
for all $1\le l'\le n$, or there exists an index $l'$, 
$1\le l'\le n$, such that $j_l=-k_{l'}$. In the first case 
$Z_G(\Delta_{j_l})$ is independent of the remaining
coordinates of the vector
$(Z_G(\Delta_{j_1}),\dots,Z_G(\Delta_{j_n}),
\overline{Z_G(\Delta_{k_1})},\dots,\overline{Z_G(\Delta_{k_n})})$,
and $EZ_G(\Delta_{j_l})=0$. Hence the expectation of the investigated
product equals zero, as we claimed. If ${j_l}=-k_{l'}$ with some index
$l'$, then a different argument is needed, since $Z_G(\Delta_{j_l})$
and $Z_G(-\Delta_{j_l})$ are not independent. In this case we can
state that since $j_p\neq\pm j_l$ if $p\neq l$, and
$k_q\neq\pm j_l$ if $q\neq l'$, the vector
$(Z_G(\Delta_{j_l}),Z_G(-\Delta_{j_l}))$ is independent of the
remaining coordinates of the above random vector. On the other hand,
the product $Z_G(\Delta_{j_l})\overline{Z_G(-\Delta_{j_l}})$
has zero expectation, since
$EZ_G(\Delta_{j_l})\overline{Z_G(-\Delta_{j_l})}
=G(\Delta_{j_l}\cap(-\Delta_{j_l}))=0$ by property~(iii) of the
random spectral measures and the relation
$\Delta_{j_l}\cap(-\Delta_{j_l})=\emptyset$ for the elements of a
regular system. Hence the expectation
of the considered product equals zero also in this case. If
$\{j_1,\dots,j_n\}=\{k_1,\dots,k_n\}$, then
$$
EZ_G(\Delta_{j_1})\cdots Z_G(\Delta_{j_n})
\overline{Z_G(\Delta_{k_1})\cdots Z_G(\Delta_{k_n})} 
=\prod\limits_{l=1}^n EZ_G(\Delta_{j_l})
\overline{Z_G(\Delta_{j_l})}=\prod\limits_{l=1}^n G(\Delta_{j_l}).
$$

If we have two functions $f,g\in\hat{{\cal H}}_G^n$, then
we may assume that they are adapted to the same regular system
${\cal D}=\{\Delta_j,\,j=\pm1,\dots,\pm N\}$. Then, by exploiting
that $I_G(g)$ is real valued, i.e. $I_G(g)=\overline{I_G(g)}$ we can 
calculate the expectation of $I_G(f)I_G(g)$ in the following way.
\begin{eqnarray*}
EI_G(f)I_G(g)&=&EI_G(f)\overline{I_G(g)}
=\left(\frac1n\right)^2\sum\sum f(x_{j_1},\dots,x_{j_n})
\overline{g(x_{k_1},\dots,x_{k_n})} \\
&&\qquad\qquad\qquad\qquad EZ_G(\Delta_{j_1})\cdots Z_G(\Delta_{j_n})
\overline{Z_G(\Delta_{k_1})\cdots Z_G(\Delta_{k_n})}\\
&=&\left(\frac1{n!}\right)^2\sum f(x_{j_1},\dots,x_{j_n})
\overline{g(x_{j_1},\dots,x_{j_n})}
G(\Delta_{j_1})\cdots G(\Delta_{j_n}) \cdot n! \\
&=&\frac1{n!}\int f(x_1,\dots,x_n)
\overline{g(x_1,\dots,x_n)} G(\,dx_1)\cdots G(\,dx_n)
=\frac1{n!}(f,g),
\end{eqnarray*}
where we took summation in the first sum for such
pairs of indices $(j_1,\dots,j_n)$ and $(k_1,\dots,k_n)$
which are permutations of each other.

I claim that Wiener--It\^o integrals of different order are
uncorrelated. More explicitly, take two functions
$f\in\hat{\bar{{\cal H}}}^n_G$ and $f'\in\hat{\bar{{\cal H}}}^{n'}_G$
such that $n\neq n'$. Then we have
\begin{equation}
EI_G(f)I_G(f')=0 \quad \textrm{if \ }f\in \hat{\bar{{\cal H}}}^n_G, \;\;
f'\in\hat{\bar{{\cal H}}}^{n'}_G, \textrm{ and \ } n\neq n'. 
\label{(4.6)}
\end{equation}
To see this relation observe that a regular system ${\cal D}$ can be
chosen is such a way that both $f$ and $f'$ are adapted to it.
Then a similar, but simpler argument as the previous one shows that
$$
EZ_G(\Delta_{j_1})\cdots Z_G(\Delta_{j_n})
\overline{Z_G(\Delta_{k_1})\cdots Z_G(\Delta_{k_{n'}})}=0
$$
for all sets of indices $\{j_1,\dots,j_n\}$ and $\{k_1,\dots,k_{n'}\}$
if $n\neq n'$, hence the sum expressing $EI_G(f)I_G(f')$ in this case
equals zero.

\medskip
We extend the definition of Wiener--It\^o integrals to a 
more general class of kernel functions with the help of 
the following Lemma~3.3. This  is a simple result whose proof
is contained in Lemma~4.1 of my Lecture Note.  Unfortunately 
this proof has a rather complicated notation, it contains 
several unpleasant technical details, hence it is rather
unpleasant to read. I inserted an Appendix to this note, 
where I try to present a more accessible proof.

\medskip\noindent
{\bf Lemma 3.3.} {\it The class of simple functions
$\hat{\bar{{\cal H}}}_G^n$ is dense in the (real) Hilbert space 
$\bar{{\cal H}}_G^n$, and the class of symmetric simple function 
$\hat{{\cal H}}_G^n$ is dense in the (real) Hilbert space 
${\cal H}_G^n$.} 

\medskip
As the transformation $I_G(f)$ is a contraction from
$\hat{\bar{{\cal H}}}_G^n$ into $L_2(\Omega,{\cal A},P)$, 
it can uniquely be extended to the closure of 
$\hat{\bar{{\cal H}}}_G^n$, i.e. to $\bar{{\cal H}}_G^n$.
(Here $(\Omega,{\cal A},P)$ denotes the probability 
space where the random spectral measure $Z_G(\cdot)$ 
is defined.) At this point we exploit that if 
$f\in\hat{\bar{{\cal H}}}_G^n$, $N=1,2,\dots$, is a 
convergent sequence in the space $\bar{{\cal H}}_G^n$, 
then the sequence of random variables $I_G(f_N)$ is 
convergent in the space $L_2(\Omega,{\cal A},P)$, 
since it is a Cauchy sequence. With the help of this 
fact and Lemma~3.3 we can introduce the definition of 
Wiener--It\^o integrals in the general case when the 
integral of a function $f\in\bar{{\cal H}}_G^n$ is taken.

\medskip\noindent
{\bf Definition of Wiener--It\^o Integrals.} {\it Given 
a function $f\in\bar{{\cal H}}^n_G$ with a spectral 
measure $G$ choose a sequence of simple functions 
$f_N\in \hat{\bar{{\cal H}}}_G^n$, $N=1,2,\dots$, which 
converges to the function $f$ in the space 
$\bar{{\cal H}}^n_G$. Such a sequence exists by Lemma~3.3. 
The random variables $I_G(f_N)$ converge to a random 
variable in the $L_2$-norm of the probability space 
where these random variables are defined, and the limit 
does not depend on the choice of the sequence $f_N$ 
converging to $f$. This enables us to define the 
$n$-fold Wiener--It\^o integral with kernel function 
$f$ as
$$
\int f(x_1,\dots,x_n)Z_G(\,dx_1)\dots Z_G(\,dx_n)
=n!I_G(f)=\lim_{N\to\infty}n!I_G(f_N),
$$
where $f_N\in\hat{\bar{{\cal H}}}^n_G$, $N=1,2,\dots$, 
is a sequence of simple functions converging to the 
function $f$ in the space $\bar{{\cal H}}^n_G$.}

\medskip
We have defined the Wiener--It\^o integral of a function
$f(x_1,\dots,x_n)$ with some nice properties. To simplify
some later considerations I introduce the following convention.
I shall sometimes consider the same function~$f$ with
reindexed variables. We say that this reindexed function 
has the same Wiener--It\^o integral. More explicitly, if 
$f(x_1,\dots,x_n)=\bar f(x_{j_1},\dots,x_{j_n})$ with arbitrary
(different) indices $j_1,\dots,j_n$, then I shall say that
$$
\int f(x_1,\dots,x_n)Z_G(\,dx_1)\dots Z_G(\,dx_n)
=\int \bar f(x_{j_1},\dots,x_{j_n})Z_G(\,dx_{j_1})\dots Z_G(\,dx_{j_n}n).
$$

\subsection{Further properties of Wiener--It\^o integrals}

I have claimed that the elements of the Hilbert spaces 
${\cal H}_n$ defined before the formulation of Theorem~3.1
can be expressed by means of $n$-fold Wiener--It\^o integrals,
and the elements of Hilbert space ${\cal H}$ can be expressed
by means of the sum of such integrals. I formulate this
statement explicitly in the next Theorem~3.4. This theorem 
also states that this representation is unique if we allow
only kernel functions from the space of symmetric functions,
${\cal H}^n_G$. Moreover, the mapping 
$I_G\colon\; \textrm{\rm Exp}\,{{\cal H}}_G\to {\cal H}$
defined after formula (\ref{(4.1)}) is a unitary transformation
between the Fock space $\textrm{\rm Exp}\,{{\cal H}}_G$ and 
${\cal H}$. A similar statement holds also for the transformation
$(n!)^{1/2}I_G\colon\; {\cal H}_G^n\to{\cal H}_n$. 

I have also claimed that the shift transformations on ${\cal H}$
can be well expressed by means of Wiener--It\^o integrals. This
will be explained in Theorem~3.6.

One can see from the definition of the $n$-fold Wiener--It\^o
integrals that they are in ${\cal H}_{\le n}$. Moreover, their
orthogonality properties imply that they are in ${\cal H}_n$.
To see that all elements of ${\cal H}_n$ can be expressed as
an $n$-fold Wiener--It\^o integral we need more information.
This can be proved with the help of Corollary~3.2 and It\^o's
formula presented in Theorem~3.5 which has a central role
in the theory of multiple Wiener--it\^o integrals. It is a
very useful result, because it helps us to represent
each elements of ${\cal H}_n$ in the form of an $n$-fold
Wiener--It\^o integral. This result  also indicates the 
strong relation between multiple Wiener--It\^o integrals 
and Hermite polynomials. I postpone its proof to the next 
section. We will get it as a consequence of the diagram 
formula, another important result in the theory of
multiple Wiener--It\^o integrals.

First I formulate Theorem~3.4. 

\medskip\noindent
{\bf Theorem 3.4.} {\it Let a stationary Gaussian random field 
be given (discrete or generalized one), and let $Z_G$ denote 
the random spectral measure adapted to it. If we integrate 
with respect to this $Z_G$, then the transformation
$I_G\colon\; \textrm{\rm Exp}\,{{\cal H}}_G\to {\cal H}$, 
where ${\cal H}$ denotes the Hilbert space of those square 
integrable random variables which are measurable with respect 
to the $\sigma$-algebra generated by the random variables of the 
random spectral measure $Z_G$ is unitary. More explicitly,
formula (\ref{(4.7)}) provides a unitary transformation between
$\textrm{Exp}\,{\cal H}_G$ and ${\cal H}$. The transformation 
$(n!)^{1/2}I_G\colon\; {\cal H}_G^n\to{\cal H}_n$ is 
also unitary.  }

\medskip
Next I formulate It\^o's formula which plays a central role
in the proof of Theorem~3.4.

\medskip\noindent
{\bf Theorem 3.5. (It\^o's Formula.)} 
{\it Let $\varphi_1,\dots,\varphi_m$, \ 
$\varphi_j\in{\cal H}_G^1$, $1\le j\le m$, be an orthonormal 
system in $L_G^2$. Let some positive integers $j_1,\dots,j_m$ 
be given, and let $j_1+\cdots+j_m=N$. Define for all 
$i=1,\dots,N$ the function $g_i$ as $g_i=\varphi_s$ for 
$j_1+\cdots+j_{s-1}<i\le j_1+\cdots+j_s$, $1\le s\le m$. (In 
particular, $g_i=\varphi_1$ for $0<i\le j_1$.) Then
\begin{eqnarray}
&&H_{j_1}\left(\int\varphi_1(x)Z_G(\,dx)\right)\cdots
H_{j_m}\left(\int\varphi_m(x)Z_G(\,dx)\right) \nonumber  \\
&&\qquad=\int g_1(x_1)\cdots g_N(x_N)\,Z_G(\,dx_1)\cdots Z_G(\,dx_N)
\nonumber\\
&&\qquad=\int \textrm{Sym\,}[ g_1(x_1)\cdots 
g_N(x_N)]\,Z_G(\,dx_1)\cdots Z_G(\,dx_N). \label{(3.8)}
\end{eqnarray}
($H_j(x)$ denotes again the $j$-th Hermite polynomial with leading
coefficient~1.)

In particular, if $\varphi\in{\cal H}_G^1$, and 
$\int\varphi^2(x)G(\,dx)=1$, then
\begin{equation}
H_n\left(\int \varphi(x)Z_G(\,dx)\right)=
\int \varphi(x_1)\cdots\varphi(x_n)\,Z_G(\,dx_1)\cdots Z_G(\,dx_n).
\label{(3.9)}
\end{equation}
}
\medskip\noindent
{\it Proof of Theorem 3.4.}\/ We have already seen that $I_G$ is 
an isometry. So it remains to show that it is a one to one map 
from $\textrm{\rm Exp}\,{{\cal H}}_G$ to ${\cal H}$ and from 
${\cal H}_G^n$ to ${\cal H}_n$.

The one-fold integral $I_G(f)$, $f\in{\cal H}_G^1$, agrees 
with the stochastic integral $I(f)$ defined in Section~2. 
Hence $I_G(e^{i(n,x)})=X(n)$ in the discrete field case, 
and $I_G(\tilde\varphi)=X(\varphi)$, $\varphi\in{\cal S}$, 
in the generalized field case. This implies that 
$I_G\colon\;{\cal H}_G^1\to{\cal H}_1$ is a unitary 
transformation. Let $\varphi_1,\varphi_2,\dots$ be a
complete orthonormal basis in ${\cal H}_G^1$. Then
$\xi_j=\int\varphi_j(x)\,Z_G(\,dx)$, $j=1,2,\dots$, 
is a complete orthonormal basis in ${\cal H}_1$. 
It\^o's formula implies that for all sets of positive 
integers $(j_1,\dots,j_m)$ the random variable 
$H_{j_1}(\xi_1)\cdots H_{j_m}(\xi_m)$ can be written 
as a $j_1+\cdots+j_m$-fold Wiener--It\^o integral. 
Therefore Theorem~3.1 implies that the image of 
$\textrm{Exp}\,{\cal H}_G$ is the whole space 
${\cal H}$, and since $EH_{j_k}(\xi_k)^2=j_k!$ the operator
$I_G\colon\;\textrm{Exp}\,{\cal H}_G\to{\cal H}$ is
unitary.

The image of ${\cal H}_G^n$ contains ${\cal H}_n$ 
because of Corollary~3.2 and It\^o's formula. Since 
these images are orthogonal for different~$n$, 
formula~(\ref{(2.1)}) implies that the image of
${\cal H}_G^n$ coincides with ${\cal H}_n$. Hence
$(n!)^{1/2}I_G\colon\; {\cal H}_G^n\to{\cal H}_n$ 
is a unitary transformation. 

\medskip
In Theorem 3.6 I shall describe the action of the shift 
transformations in ${\cal H}$. To do this let us first
remark that by Theorem~3.4 all $\eta\in{\cal H}$ can be 
written in the form
\begin{equation}
\eta=f_0+\sum_{n=1}^\infty\frac1{n!}
\int f_n(x_1,\dots,x_n)Z_G(\,dx_1)
\dots Z_G(\,dx_n) \label{(4.7)}
\end{equation}
with $f=(f_0,f_1,\dots)\in\textrm{Exp}\,{\cal H}_G$ 
in a unique way, where $Z_G$ is the random measure 
adapted to the stationary Gaussian field. Now, we have

\medskip\noindent
{\bf Theorem 3.6.} {\it Let $\eta\in{\cal H}$ have the 
form~(\ref{(4.7)}). Then
$$
T_t\eta=f_0+\sum_{n=1}^\infty\frac1{n!}\int e^{i(t,x_1+\cdots+x_n)}
f_n(x_1,\dots,x_n)Z_G(\,dx_1)
\dots Z_G(\,dx_n)
$$
for all $t\in R^\nu$ in the generalized field and for all
$t\in{\mathbb Z}_\nu$ in the discrete field case.}

\medskip\noindent
{\it Proof of Theorem 3.6.} Because of formulas~(\ref{(3.6)}) 
and~(\ref{($3.6'$)}) and
the definition of the shift operator~$T_t$ we have
$$
T_t\left(\int e^{i(n,x)}Z_G(\,dx)\right)=T_tX_n=X_{n+t}
=\int e^{i(t,x)}e^{i(n,x)}Z_G(\,dx), \quad t\in{\mathbb Z}_\nu,
$$
and because of the identity
$\widetilde{T_t\varphi}(x)=\int e^{(i(u,x)}\varphi(u-t)\,du
=e^{i(t,x)}\tilde\varphi(x)$ for $\varphi\in{\cal S}$
\begin{eqnarray*}
T_t\left(\int\tilde\varphi(x)\,Z_G(\,dx)\right)
&=&T_tX(\varphi)=X(T_t\varphi) \\
&=&\int e^{i(t,x)}\tilde\varphi(x)Z_G(\,dx), \quad \varphi\in{\cal S},
\quad t\in R^\nu,
\end{eqnarray*}
in the discrete and generalized field cases respectively. Observe,
that the finite linear combinations of the functions $e^{i(n,x)}$,
$n\in{\mathbb Z}_\nu$, is dense in the space ${\cal H}_1$ in the 
case of discrete stationary random fields, and the functions 
$\varphi\in{\cal S}_\nu$ are dense in ${\cal H}_1$ in the
generalized stationary field case. Hence 
$$
T_t\left(\int f(x) Z_G(\,dx)\right)=\int e^{i(t,x)}f(x) Z_G(\,dx) \quad
\textrm{if } f\in {\cal H}^1_G
$$
for all $t\in{\mathbb Z}_\nu$ in the discrete field and for all
$t\in R^\nu$ in the generalized field case. This means that
Theorem~3.6 holds in the special case when $\eta$ is a one-fold
Wiener--It\^o integral. Let $f_1(x),\dots,f_m(x)$ be an
orthogonal system in ${\cal H}^1_G$. The set of functions
$e^{i(t,x)}f_1(x),\dots,e^{i(t,x)}f_m(x)$ is also an orthogonal
system in ${\cal H}^1_G$. (We choose $t\in{\mathbb Z}_\nu$ in 
the discrete and $t\in R^\nu$ in the generalized field case.) 
Hence It\^o's formula implies that Theorem~3.6 also holds for 
random variables of the form
$$
\eta=H_{j_1}\left(\int f_1(x)Z_G(\,dx)\right)\cdots
H_{j_m}\left(\int f_m(x)Z_G(\,dx)\right)
$$
and for their finite linear combinations. Since these linear
combinations are dense in ${\cal H}$ Theorem~3.6 holds true.

\medskip
I formulate at the end of this section a somewhat technical,
but rather natural result. It is a formula for the change 
of variables in Wiener--It\^o integrals. It can be 
interpreted so that we describe how to deal with the 
situation when we integrate with respect to the random 
spectral measure $Z_{G'}(\,dx)=g^{-1}(x)Z_G(\,dx)$ instead of 
$Z_G(\,dx)$ with some function $g(\cdot)$.
(We have to assume that $g(x)=\overline{g(-x)}$ to get a
random spectral measure again.) This new random
measure corresponds to the spectral measure 
$G'(\,dx)=|g^{-2}(x)|G(\,dx)$, and to preserve the value of 
the (sum of) integrals we are working with we have to multiply 
the kernel function $f_n(x_1,\dots,x_n)$ by 
$\prod\limits_{j=1}^ng(x_j)$
to compensate the multiplying factor $g^{-1}(x)$ in the definition
of $Z_{G'}(\,dx)$.

\medskip\noindent
{\bf Theorem~3.7.} {\it Let $G$ and $G'$ be two non-atomic spectral
measures such that $G$ is absolutely continuous with respect to $G'$,
and let $g(x)$ be a complex valued function such that
\begin{eqnarray*}
g(x)&=&\overline{g(-x)}, \\
|g^2(x)|&=&\frac{dG(x)}{dG'(x)}.
\end{eqnarray*}
For every $f=(f_0,f_1,\dots)\in\textrm{\rm Exp}\,{\cal H}_G$, we define
$$
f'_n(x_1,\dots,x_n)=f_n(x_1,\dots,x_n)g(x_1)\cdots g(x_n),
\quad n=1,2,\dots,\quad f'_0=f_0.
$$
Then $f'=(f'_0,f'_1,\dots)\in \textrm{\rm Exp}\,{\cal H}_{G'}^n$, and
\begin{eqnarray*}
&&f_0+\sum_{n=1}^\infty\int\frac1{n!}f_n(x_1,\dots,x_n)
Z_G(\,dx_1)\dots Z_G(\,dx_n) \\
&&\qquad \stackrel{\Delta}{=}f'_0+\sum_{n=1}^\infty\frac1{n!}
\int f_n'(x_1,\dots,x_n)Z_{G'}(\,dx_1) \dots Z_{G'}(\,dx_n),
\end{eqnarray*}
where $Z_G$ and $Z_{G'}$ are Gaussian random spectral measures
corresponding to $G$ and~$G'$.}

\medskip\noindent
{\it Proof of Theorem 3.7.} We have $\|f'_n\|_{G'}=\|f_n\|_G$, hence
$f'\in\textrm{Exp}\,{\cal H}_{G'}$. Let us choose a complete 
orthonormal system $\varphi_1,\varphi_2,\dots$ in ${\cal H}_G^1$. 
Then $\varphi'_1,\varphi_2',\dots$, with $\varphi'_j(x)=\varphi_j(x)g(x)$ 
for all $j=1,2,\dots$ is a complete orthonormal system in ${\cal H}^1_{G'}$.
All functions $f_n\in{\cal H}_G^n$ can be written in the form
$$
f(x_1,\dots,x_n)=\sum c_{j_1,\dots,j_n}
\textrm{Sym}\,(\varphi_{j_1}(x_1)\cdots
\varphi_{j_n}(x_n)).
$$ 
Then $f'(x_1,\dots,x_n)=\sum c_{j_1,\dots,j_n}
\textrm{Sym(}\,\varphi'_{j_1}(x_1)\cdots
\varphi'_{j_n}(x_n))$. Rewriting all terms
$$
\int\textrm{Sym}\,(\varphi_{j_1}(x_1)\cdots\varphi_{j_n}(x_n))
Z_G(\,dx_1)\dots Z_G(,dx_n)
$$
and
$$
\int\textrm{Sym}\,(\varphi'_{j_1}(x_1)\cdots\varphi'_{j_n}(x_n))
Z_{G'}(\,dx_1)\dots Z_{G'}(,dx_n)
$$
by means of It\^o's formula we
get that $f$ and $f'$ depend on a sequence of independent standard
normal random variables in the same way. Theorem~3.7 is proved.

\Section{The diagram formula and It\^o's formula }

The first main subject of this section is the diagram formula 
for the product of Wiener--It\^o integrals. This formula enables 
us to rewrite the product of Wiener--It\^o integrals in the 
form of a sum of Wiener--It\^o integrals of different 
multiplicity. A Wiener--It\^o integral is an element of the 
Hilbert space ${\cal H}$ consisting from square integrable 
random variables measurable with respect to the 
$\sigma$-algebra generated by the random variables of the 
underlying Gaussian stationary random field. Moreover, every 
element of the Hilbert space ${\cal H}$ can be expressed as 
the sum of Wiener--It\^o integrals of different multiplicity. 
Thus, by the diagram formula the product of  two Wiener--It\^o 
integrals also belongs to the Hilbert space ${\cal H}$. The 
measurability property of the product is obvious, but the fact 
that its second moment is finite requires some explanation. 
Besides, the diagram formula is not a simple existence result, 
it also gives an explicit formula about how to rewrite the 
product of Wiener--It\^o integrals as a sum of such integrals.

The other subject of this section which will be discussed in
a special subsection is the proof of It\^o's formula formulated 
in Theorem~3.5. The proof is made by induction which  is based 
on the similarity of a recursion formula for a product of 
special form of Wiener--It\^o integrals and a recursion 
formula about the relation of Hermite polynomials of different 
order. The recursion formula for Wiener--It\^o integrals we 
apply is a simple consequence of the diagram formula.

\medskip
To formulate the diagram formula first we have to introduce
some notations.

Let us fix some functions $h_1\in\bar{{\cal H}}_G^{n_1}$,\dots,
$h_m\in\bar{{\cal H}}_G^{n_m}$. In the  diagram formula, we 
shall express the product $n_1!I_G(h_{n_1})\cdots n_m!I_G(h_{n_m})$ 
as the sum of Wiener--It\^o integrals. There is no unique 
terminology for this result in the literature.
I shall follow the notation of Dobrushin's paper
{\it Gaussian and their subordinated fields.} Annals of
Probability~7, 1--28 (1979).

We introduce a class of diagrams $\gamma$ denoted by
$\Gamma(n_1,\dots,n_m)$, and define with the help of
each diagram $\gamma$ in this class a function $h_\gamma$ 
which will be the kernel function of one of the 
Wiener--It\^o integrals taking part in the sum expressing 
the product of the Wiener--It\^o integrals we investigate.
First we define the diagrams $\gamma$ and the functions 
$h_\gamma$ corresponding to them, and then we formulate 
the diagram formula with their help. After the formulation 
of this result I present an example together with some 
figures which may help in understanding better what the 
diagram formula is like.

We shall use the term diagram of order $(n_1,\dots,n_m)$ 
for an undirected graph of $n_1+\cdots+n_m$ vertices
which are indexed by the pairs of integers~$(j,l)$, 
$l=1,\dots,m$, \ $j=1,\dots,n_l$ if the second term in the
pair $(j,l)$ equals~$l$, and we shall call the 
set of vertices $(j,l)$, $1\le j\le n_l$ the $l$-th row 
of the diagram.  The diagrams of order $(n_1,\dots,n_m)$ 
are those undirected graphs with these vertices which 
have the properties that no more than one edge enters 
into  each vertex, and  edges can connect only pairs 
of vertices from different rows of a diagram, i.e. 
such vertices~$(j_1,l_1)$ and $(j_2,l_2)$ for which 
$l_1\neq l_2$. Let $\Gamma=\Gamma(n_1,\dots,n_m)$ denote 
the set of all diagrams of order $(n_1,\dots,n_m)$.

Given a diagram $\gamma\in\Gamma$ let $|\gamma|$ denote the
number of edges in~$\gamma$. Let there be given a set of 
functions $h_1\in\bar{{\cal H}}_G^{n_1}$,\dots, 
$h_m\in\bar{{\cal H}}_G^{n_m}$.  Let us denote the 
variables of the function $h_l$ by $x_{(j,l)}$ instead of 
$x_j$, i.e. let us write $h_l(x_{(1,l)},\dots,x_{(n_l,l)})$
instead of $h_l(x_1,\dots,x_{n_l})$.  Put $N=n_1+\cdots+n_m$. 
We introduce the function of $N$~variables corresponding 
to the vertices of the diagram by the formula
\begin{equation}
h(x_{(j,l)},\; l=1,\dots,m,\;j=1,\dots,n_l)=\prod_{l=1}^m
h_l(x_{(j,l)},\;j=1,\dots,n_l). \label{(5.1a)}
\end{equation}

For each diagram $\gamma\in\Gamma=\Gamma(n_1,\dots,n_m)$ 
we define the reenumeration of the indices of the function 
in formula~(\ref{(5.1a)}) in the following way. We enumerate the 
variables $x_{(j,l)}$ in such a way that the vertices 
into which no edge enters will have the indices 
$1,2,\dots,N-2|\gamma|$, and the vertices connected by an
edge will have the indices~$p$ and $p+|\gamma|$, where
$p=N-2|\gamma|+1,\dots,N-|\gamma|$. In such a way we have
defined a function $h(x_1,\dots,x_N)$ (with an enumeration 
of the indices of the variables depending on the diagram 
$\gamma$). After the definition of this function 
$h(x_1,\dots,x_N)$ we take that function of $N-|\gamma|$ 
variables which we get by replacing the arguments 
$x_{N-|\gamma|+p}$ by the arguments $-x_{N-2|\gamma|+p}$, 
$1\le p\le|\gamma|$, in the function $h(x_1,\dots,x_N)$. 
Then we define the function $h_\gamma$ appearing in the 
diagram formula by integrating this function by the 
product measure  
$\prod\limits_{p=1}^{|\gamma|}G(\,dx_{N-2|\gamma|+p})$. 

More explicitly, we write
\begin{eqnarray}
h_\gamma(x_1,\dots,x_{N-2|\gamma|})&=&\idotsint 
h(x_1,\dots,x_{N-|\gamma|},
-x_{N-2|\gamma|+1},\dots,-x_{N-|\gamma|}) \nonumber  \\
&&\qquad G(\,dx_{N-2|\gamma|+1})\dots G(\,dx_{N-|\gamma|}).
\label{(5.1)}
\end{eqnarray} 
The function $h_\gamma$ depends only on the variables
$x_1,\dots,x_{N-2|\gamma|}$, i.e. it is independent of how 
the vertices connected by edges are indexed. Indeed, it 
follows from the evenness of the spectral measure that 
by interchanging the indices $s$ and $s+\gamma$ of two 
vertices connected by an edge we do not change the value 
of the  integral $h_\gamma$. Let us now consider the
Wiener--It\^o integrals $(N-2|\gamma|)!I_G(h_\gamma)$. In the
diagram formula we shall show that the product of the
Wiener--It\^o integrals we considered  can be expressed
as the sum of these Wiener--It\^o integrals. To see that
the identity appearing in the diagram formula is 
meaningful observe that although the function $h_\gamma$ 
may depend on the numbering of those vertices 
of~$\gamma$ from which no edge starts, the function 
$\textrm{Sym}\,h_\gamma$, and therefore the Wiener--It\^o 
integral $I_G(h_\gamma)$ does not depend on it. 

Now I formulate the diagram formula. Then I 
make a remark about the definition of the function 
$h_\gamma$ in it, and discuss an example to show 
how to calculate the terms appearing in this result.

\medskip\noindent
{\bf Theorem 4.1. (Diagram Formula.)} 
{\it For all functions $h_1\in\bar{{\cal H}}_G^{n_1}$,\dots, 
$h_m\in\bar{{\cal H}}_G^{n_m}$, \ $n_1,\dots,n_m=1,2,\dots$, 
the following relations hold with $N=n_1+\cdots+n_m$:

\medskip
\begin{description}
\item[(A)] $h_\gamma\in\bar{{\cal H}}_G^{N-2|\gamma|}$, and
$\|h_\gamma\|\le\prod\limits_{j=1}^m\|h_j\|$ for all 
$\gamma\in\Gamma$.

\item[(B)] $n_1!I_G(h_1)\cdots n_m! I_G(h_m)
=\sum\limits_{\gamma\in\Gamma}(N-2|\gamma|)!I_G(h_\gamma)$.
\end{description}

\medskip
Here $\Gamma=\Gamma(n_1,\dots,n_m)$, and the functions $h_\gamma$
agree with the functions $h_\gamma$ defined before the formulation of 
Theorem~4.1. In particular, $h_\gamma$ was defined in~(\ref{(5.1)}).}

\medskip\noindent
{\it Remark.}\/ Observe that at the end of the definition
of the function $h_\gamma$ we replaced the variable
$x_{N-|\gamma|+p}$ by the variable $-x_{N-2|\gamma|+p}$
and not by $x_{N-2|\gamma|+p}$. This is related to the 
fact that in the Wiener--It\^o integral we integrate 
with respect a complex valued random measure $Z_G$ which 
has the property 
$EZ_G(\Delta)Z_G(-\Delta)=EZ_G(\Delta)\overline{Z_G(\Delta)}
=G(\Delta)$, while $EZ_G(\Delta)Z_G(\Delta)=0$ if 
$\Delta\cap(-\Delta)=\emptyset$. In the case of the
original Wiener--It\^o integral with respect to a Gaussian
orthogonal random measure
the situation is a bit different. In that case we 
integrate with respect to a real valued Gaussian 
orthonormal random measure $Z_\mu$ which has the property 
$EZ_\mu^2(\Delta)=\mu(\Delta)$. In that case a diagram 
formula also holds, but it has a slightly different form.
The main difference is that in that case we define the
function $h_\gamma$ (because of the above mentioned 
property of the random measure $Z_\mu$) by replacing
the variable $x_{N-|\gamma|+p}$ by the variable 
$x_{N-2|\gamma|+p}$.

\medskip
To make the notation in the diagram formula more 
understandable let us consider the following example. 

\medskip\noindent
{\it Example.}\/ Let us take four functions 
$h_1=h_1(x_1,x_2,x_3)\in\bar{\cal{H}}^3_G$,
$h_2=h_2(x_1,x_2)\in\bar{\cal{H}}^2_G$,
$h_3=h_3(x_1,x_2,x_3,x_4,x_5)\in\bar{\cal{H}}^5_G$ 
and $h_4=h_4(x_1,x_2,x_3,x_4)\in\bar{\cal{H}}^4_G$, 
and consider the product of Wiener--It\^o integrals 
$3!I_G(h_1)2!I_G(h_2)5!I_G(h_3)4!I_G(h_4)$. Let us
look how to calculate the kernel function $h_\gamma$ 
of a Wiener--It\^o integral 
$(14-2|\gamma|)!I_G(h_\gamma)$, corresponding to a diagram
$\gamma\in\Gamma(3,2,5,4)$, in the diagram formula.

\begin{figure}[ht]
\begin{center}
\epsfig{file=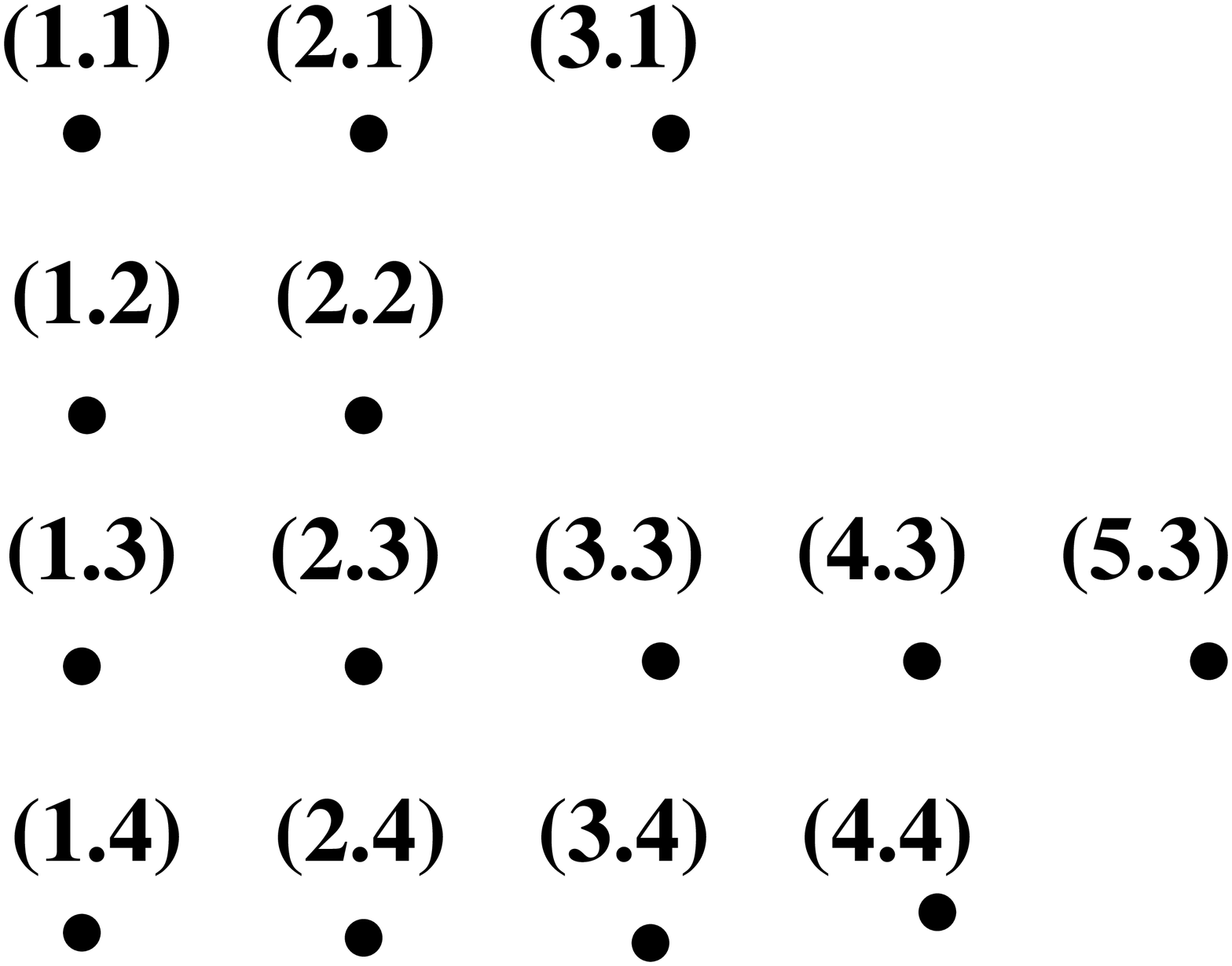, width=4cm}
\end{center}
\caption{The vertices of the diagrams $\gamma\in\Gamma(3,2,5,4)$}
\label{Diagram1}
\end{figure}

\medskip\noindent
We have to consider the class of diagrams $\Gamma(3,2,5,4)$,
i.e. the diagrams with vertices which are indexed in the first 
row as $(1,1)$, $(2,1)$, $(3,1)$, in the second row as 
$(1,2)$, $(2,2)$, in the third row in as $(1,3)$, $(2,3)$, 
$(3,3)$, $(4,3)$, $(5,3)$ and in the fourth row as $(1,4)$, 
$(2,4)$, $(3,4)$, $(4,4)$. (See Fig.~\ref{Diagram1}.)

Let us take a diagram $\gamma\in\Gamma(3,2,5,4)$, and let 
us see how we can calculate the kernel function $h_\gamma$ 
of the Wiener--It\^o integral corresponding to it. We 
also draw some pictures which may help in following this 
calculation. Let us consider for instance the diagram 
$\gamma\in\Gamma(3,2,5,4)$ with edges $((2,1),(4,3))$, 
$((3,1),(1,3))$, $((1,2),(2,4))$, $((2,2),(5,3))$, 
$((3,3),(3,4))$. Let us draw the diagram~$\gamma$ 
with its edges and with such a reenumeration of the 
vertices which helps in writing up the function 
$h(\cdot)$ (with $N=14$ variables) corresponding to 
this diagram $\gamma$ and introduced before the 
definition of the function $h_\gamma$.

\begin{figure}[ht]
\begin{center}
\epsfig{file=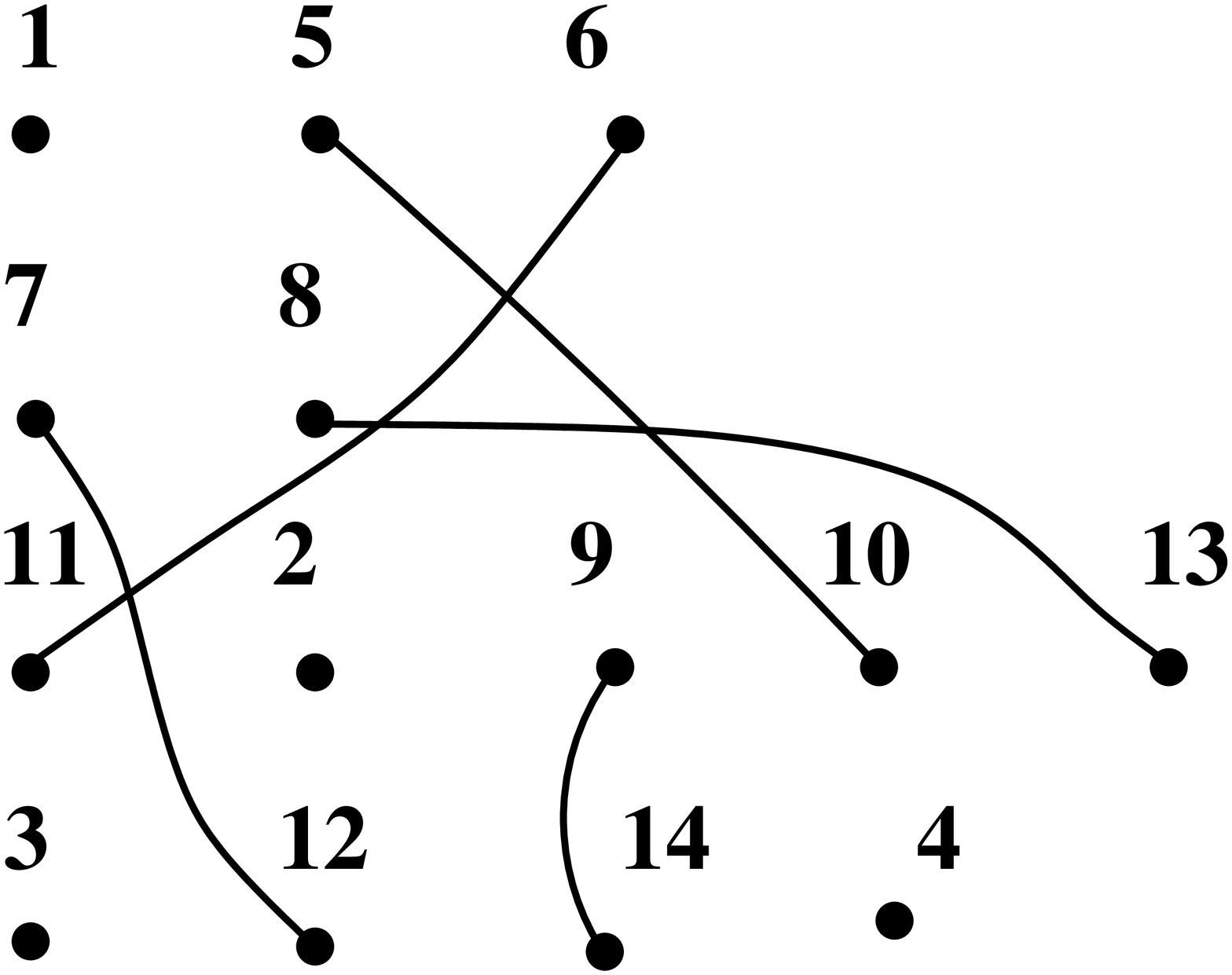, width=4cm}
\end{center}
\caption{The diagram $\gamma$ we are working with and the 
reenumeration of its vertices.}
\label{Diagram2}
\end{figure}

The function defined in~(\ref{(5.1a)}) equals in the present case
\begin{eqnarray*}
&&h_1(x_{(1,1)},x_{(2,1)},x_{(3,1)}) h_2(x_{(1,2)},x_{(2,2)})
h_3(x_{(1,3)},x_{(2,3)},x_{(3,3)},x_{(4,3)},x_{(5,3)}) \\ 
&&\qquad h_4(x_{(1,4)},x_{(2,4)},x_{(3,4)},x_{(4,4)}).
\end{eqnarray*}
The variables of this function are indexed by the labels
of the vertices of~$\gamma$. We made a relabelling of
the vertices of the diagram~$\gamma$ in such a way that
by changing the indices of the above function with the
help of this relabelling we get the function $h(\cdot)$ 
corresponding to the diagram~$\gamma$. In the next step 
we shall make such a new relabelling of the vertices 
of~$\gamma$ which helps to write up the functions
$h_\gamma$ we are interested in. (See Fig.~\ref{Diagram2}) 

The function $h(\cdot)$ (with $N=14$ variables) 
corresponding to the diagram $\gamma$ can be written 
(with the help of the labels of the vertices in the 
second diagram) as
\begin{eqnarray*}
&&h(x_1,x_2,\dots,x_{14}) \\
&&\qquad=h_1(x_1,x_5,x_6) h_2(x_7,x_8)
h_3(x_{11},x_2,x_9,x_{10},x_{13}) h_4(x_3,x_{12},x_{14},x_4).
\end{eqnarray*}

\begin{figure}[ht]
\begin{center}
\epsfig{file=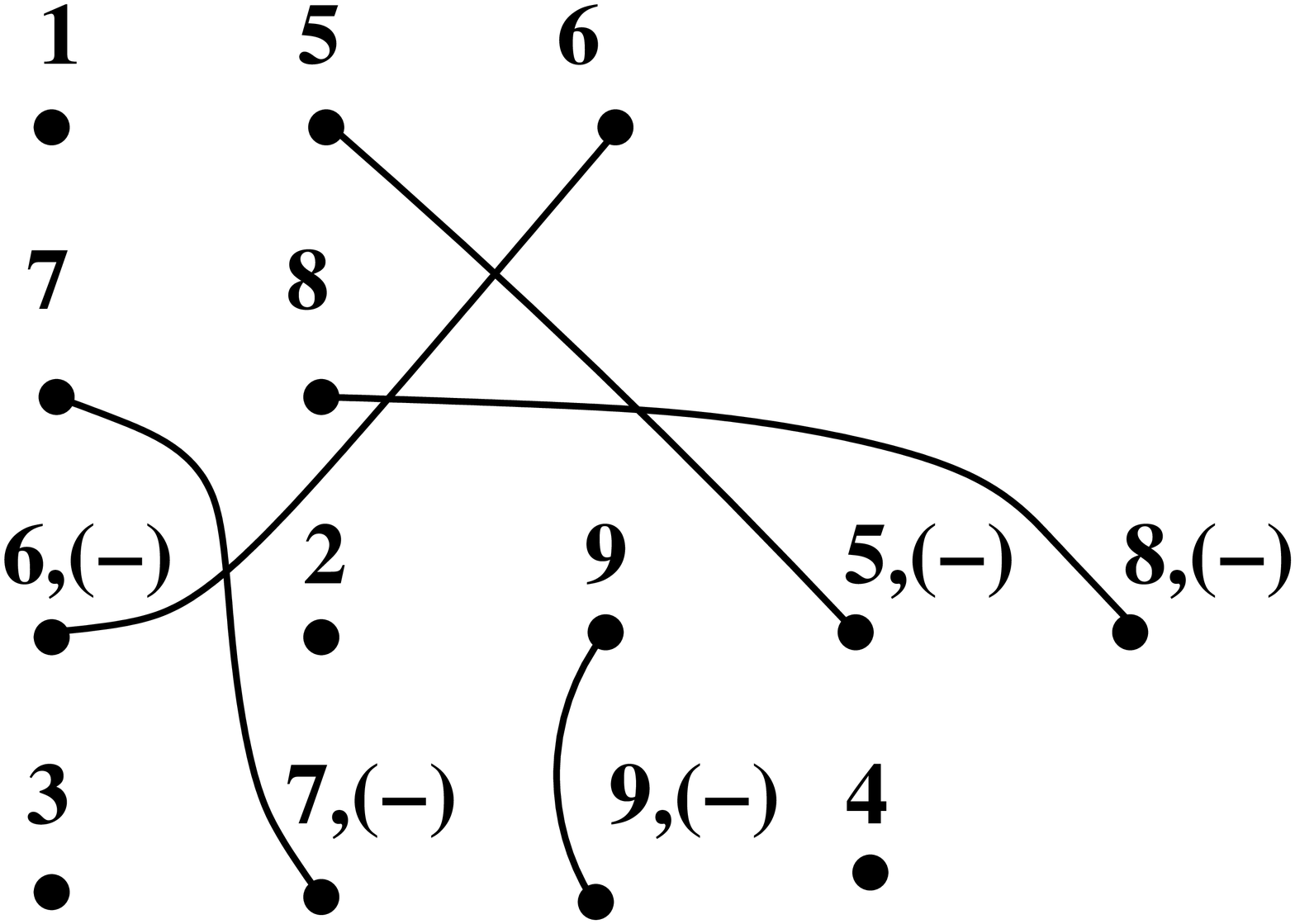, width=4cm}
\end{center}
\caption{The diagram applied for the calculation
of $h_\gamma$. The sign $-$ indicates that the corresponding
argument is multiplied by $-1$.}
\label{Diagram3}
\end{figure}

Let us change the enumeration of the vertices of the diagram
in a way that corresponds to the change of the arguments
$x_{N-|\gamma|+p}$ by the arguments $-x_{N-2|\gamma|+p}$. 
This is done in the next picture. (In this notation the 
sign $(-)$ denotes that the variable corresponding to this 
vertex is $-x_{N-2|\gamma|+p}$ and not $x_{N-2|\gamma|+p}$.
(See Fig~\ref{Diagram3}.)

With the help of the above diagram we can write up the function 
$$
h(x_1,\dots,x_{N-|\gamma|},-x_{N-2|\gamma|+1},\dots,-x_{N-|\gamma|})
$$ 
corresponding to the diagram~$\gamma$ in a simple way. This
yields that in the present case the function $h_\gamma$ 
defined in~(\ref{(5.1)}) can be written in the form
\begin{eqnarray*}
h_\gamma(x_1,x_2,x_3,x_4)\!\!\!\!\!\!  &&=\idotsint 
h_1(x_1,x_5,x_6) h_2(x_7,x_8) h_3(-x_6,x_2,x_9,-x_5,-x_8) \\ 
&&\quad h_4(x_3,-x_7,-x_9,x_4) 
G(\,dx_5)G(\,dx_6)G(\,dx_7)G(\,dx_8) G(\,dx_9).
\end{eqnarray*}

Here we integrate with respect to those variables $x_j$ whose
indices correspond to such a vertex of the last diagram from
which an edge starts.
Then the contribution of the diagram $\gamma$ to the sum at 
the right-hand side of diagram formula equals $4!I_G(h_\gamma)$ 
with this function~$h_\gamma$.

Let me remark that we had some freedom in choosing the 
enumeration of the vertices of the diagram~$\gamma$. Thus e.g. 
we could have enumerated the four vertices of the diagram from 
which no edge starts with the numbers 1, 2, 3 and 4 in an 
arbitrary order. A different indexation of these vertices would 
lead to a different function $h_\gamma$ whose Wiener--it\^o 
integral is the same. I have chosen that enumeration of the 
vertices which seemed to be the most natural for me.

\medskip
I shall omit the details of the proof of Theorem~4.1, because
it contains several complicated, unpleasant details. 
I only briefly explain the main ideas. 

\medskip\noindent
{\it The idea of the proof of Theorem~4.1.}\/ It suffices 
to prove Theorem~4.1 in the special case $m=2$. Then the 
case $m>2$ follows by induction.

We shall use the notation $n_1=n$, $n_2=m$, and we write
$x_1,\dots,x_{n+m}$ instead of
$x_{(1,1)},\dots,x_{(n,1)},x_{(1,2)}\dots,x_{(m,2)}$. It is 
clear that the function $h_\gamma$ satisfies Property~(a) 
of the classes $\bar{{\cal H}}_G^{n+m-2|\gamma|}$ defined 
in subsection~3.1. We show that Part~(A) of Theorem~4.1 is a 
consequence of the Schwartz inequality. I write down 
this part of the proof, because this is simple. The validity 
of this inequality means in particular that the functions 
$h_\gamma$ satisfy also Property~(b) of the class of 
functions $\bar{{\cal H}}_G^{n+m-2|\gamma|}$.
 
To prove this estimate on the norm of $h_\gamma$ it is 
enough to restrict ourselves to such diagrams~$\gamma$ 
in which the vertices $(n,1)$ and $(m,2)$, \ $(n-1,1)$ and 
$(m-1,2)$,\dots, $(n-k,1)$ and $(m-k,2)$ are connected by 
edges with some number $0\le k\le \min(n,m)$. 
In this case we can write
\begin{eqnarray*}
&&|h_\gamma(x_1,\dots,x_{n-k-1},x_{n+1},\dots,x_{n+m-k-1})|^2 \\
&&\qquad=\biggl|\int h_1(x_1,\dots,x_n)h_2(x_{n+1},
\dots,x_{n+m-k-1},-x_{n-k},
\dots,-x_n)\\
&&\qquad\qquad\qquad G(\,dx_{n-k})\dots G(\,dx_n)\biggr|^2 \\
&&\qquad\le \int |h_1(x_1,\dots,x_n)|^2 G(\,dx_{n-k})\dots G(\,dx_n) \\
&&\qquad\qquad\qquad
\int |h_2(x_{n+1},\dots,x_{n+m})|^2 G(\,dx_{n+m-k})\dots G(\,dx_{n+m})
\end{eqnarray*}
by the Schwartz inequality and the symmetry $G(-A)=G(A)$ of the
spectral measure~$G$. Integrating this inequality with respect to
the free variables we get Part~(A) of Theorem~4.1.

In the proof of Part~(B) first we restrict ourselves to 
the case when $h_1\in\hat{\bar{{\cal H}}}_G^n$ and 
$h_2\in\hat{\bar{{\cal H}}}_G^m$, i.e. to the case when 
they are simple functions. Moreover, we may assume that
they are adapted to such a regular system of subsets
$\Delta_j\in{\cal D}$ which satisfy the inequality 
$G(\Delta_j)<\varepsilon$ with a very small number 
$\varepsilon>0$. At this reduction we exploit that the
measure $G$ is non-atomic. This enables us to split up
the elements of a regular system to very small subsets.
By making a good approximation of the function $h_1$
and $h_2$ with such elementary functions and then
taking a limiting procedure with $\varepsilon\to0$
we get the proof of Part~(B). During the limiting
procedure we may exploit the already proven Part~(A)
of Theorem~4.1. 

To prove Part~(B) in this case let us consider a
regular system 
${\cal D}=\{\Delta_j,\;j=\pm1,\dots,\pm N\}$ of subsets 
of $R^n$ such that the functions $h_1$ and $h_2$
are adapted to it, and its elements satisfy the inequality 
$G(\Delta_j)<\varepsilon$ with a very small $\varepsilon>0$. 
Let us fix a point $u_j\in\Delta_j$ in all sets 
$\Delta_j\in{\cal D}$. 
We can express the product $n!I_G(h_1)m!I_G(h_2)$ as
\begin{eqnarray*}
I=n!I_G(h_1)m!I_G(h_2)&=&{\sum}' h_1(u_{j_1},\dots,u_{j_n})
h_2(u_{k_1},\dots,u_{k_m})\\
&&\qquad\qquad Z_G(\Delta_{j_1})\cdots Z_G(\Delta_{j_n})
Z_G(\Delta_{k_1})\cdots Z_G(\Delta_{k_m})
\end{eqnarray*}
with the numbers $u_{j_p}\in\Delta_{j_p}$ and
$u_{k_r}\in\Delta_{k_r}$ we have fixed. Here the summation in
$\sum'$ goes through all pairs $((j_1,\dots,j_n),(k_1,\dots,k_m))$,
$j_p,\,k_r\in\{\pm1,\dots,\pm N\}$, \ $p=1,\dots,n$, $r=1,\dots,m$,
such that $j_p\neq\pm j_{\bar p}$ and $k_r\neq\pm k_{\bar r}$ if
$p\neq\bar p$ or $r\neq\bar r$.

Write
\begin{eqnarray}
I&=&\sum_{\gamma\in\Gamma}{\sum}^\gamma \,
h_1(u_{j_1},\dots,u_{j_n}) h_2(u_{k_1},\dots,u_{k_m}) \nonumber \\
&&\qquad\qquad\qquad Z_G(\Delta_{j_1})\cdots Z_G(\Delta_{j_n})
Z_G(\Delta_{k_1})\cdots Z_G(\Delta_{k_m}),  \label{(5.13)}
\end{eqnarray}
where $\sum^\gamma$ contains those terms of $\sum'$ for which
$j_p=k_r$ or $j_p=-k_r$ if the vertices $(1,p)$ and $(2,r)$ are
connected in $\gamma$, and $j_p\neq \pm k_r$ if $(1,p)$ and $(2,r)$
are not connected.

Let us introduce the notation
$$
\Sigma^\gamma={\sum}^\gamma\, h_1(u_{j_1},\dots,u_{j_n}) h_2(u_{k_1},\dots,u_{k_m})
Z_G(\Delta_{j_1})\cdots Z_G(\Delta_{j_n})
$$
for all $\gamma\in\Gamma$.
We prove Theorem~4.1 if we show that the inner sum $\Sigma^\gamma$
in formula (\ref{(5.13)}) is very close to 
$(n+m-2|\gamma|)!I_G(h_\gamma)$ for all $\gamma\in\Gamma$ if 
$\varepsilon>0$ is chosen very small.
To explain why it is so we make a good approximation
of $\Sigma^\gamma$. For this goal we introduce the following
 notations. Put
\begin{eqnarray*}
A_1&=&A_1(\gamma)=\{p\colon\;p\in\{1,\dots,n\},\textrm{ and no
edge starts from }(p,1)\textrm{ in }\gamma\},\\
A_2&=&A_2(\gamma)=\{r\colon\;r\in\{1,\dots,m\},\textrm{ and no
edge starts from }(r,2)\textrm{ in }\gamma\}
\end{eqnarray*}
and
\begin{eqnarray*}
B=B(\gamma)&=&\{(p,r)\colon\; p\in\{1,\dots,n\},\; 
r\in\{1,\dots,m\},\\
&&\qquad\qquad (p,1) \textrm{ and } (r,2) 
\textrm{ are connected in }\gamma\}.
\end{eqnarray*}
We define with the help of this notation the expression
\begin{eqnarray*}
\Sigma_1^\gamma&=&{\sum}^\gamma\, h_1(u_{j_1},\dots,u_{j_n})
h_2(u_{k_1},\dots,u_{k_m})
\prod_{p\in A_1}Z_G(\Delta_{j_p})\prod_{r\in A_2} Z_G(\Delta_{k_r})\\
&&\qquad\qquad\qquad\cdot
\prod_{(p,r)\in B} E\left(Z_G(\Delta_{j_p})Z_G(\Delta_{k_r})\right).
\end{eqnarray*}
It can be shown that for one part $\Sigma_1^\gamma$ is very close
to $\Sigma^\gamma$, and on the other hand it approximates well 
$(n+m-2\gamma)!I_G(h_\gamma)$ if $\varepsilon>0$ is small. 

We can prove the first statement by showing that 
$E(\Sigma^\gamma-\Sigma^\gamma_1)^2$ is very small for small $\varepsilon$.
To prove the second statement we show that the
expression $\Sigma^\gamma_1$ is very similar to the integral defining
$(n+m-2\gamma)!I_G(h_\gamma)$. To see this observe that the terms
$E\left(Z_G(\Delta_{j_p})Z_G(\Delta_{k_r})\right)$ in the expression
$\Sigma^\gamma_1$ can be simplified. Indeed, since the terms 
of these products have indices $(p,r)\in B$, we have 
$j_p=\pm k_r$, and the products with such indices satisfy either 
the identity
$E\left(Z_G(\Delta_{j_p})Z_G(-\Delta_{j_p})\right)=G(\Delta_{j_p})$ 
or the identity $EZ_G(\Delta_{j_p})^2=0$. Writing these
relations in the expression $\Sigma^\gamma_1$, and by exploiting 
the properties of the functions $h_1$ and $h_2$ we get that the 
sum $\Sigma^\gamma_1$ can be written as an $(n+m-2|\gamma|)$-fold 
Wiener--It\^o integral of an elementary function which almost 
agrees with the function $h_\gamma$. (There is a small difference 
because this elementary function disappears on a small set 
where the function $h_\gamma$ may not disappear. This set 
contains such points $x$ which have two different coordinates 
$x_u$, $x_v$ with indices $u\neq v$ such that $x_u\in\Delta_j$
and $x_v\in\pm\Delta_j$ with the same element $\Delta_j$ of 
the regular system~$\cal D$.)

A careful analysis shows that both properties mentioned before
hold. Their proof is natural, but it requires the application 
of a rather complicated notation. Hence I omitted 
the explanation of the details. Next I turn to the proof of 
It\^o's formula.

\subsection{The proof of It\^o's formula}

We shall prove It\^o's formula with the help of two results.
Here is the first one.

\medskip\noindent
{\bf Proposition 4.2.} {\it Let $f\in\bar{{\cal H}}_G^n$ and
$h\in\bar{{\cal H}}_G^1$. Let us define the functions
$$
f\underset{k}{\times} h(x_1,\dots,x_{k-1},x_{k+1},\dots,x_n)
=\int f(x_1,\dots,x_n)\overline{h(x_k)} G(\,dx_k), \quad k=1,\dots,n,
$$
and
$$
fh(x_1,\dots,x_{n+1})=f(x_1,\dots,x_n)h(x_{n+1}).
$$
Then $f\underset{k}{\times} h$, $k=1,\dots,n$, and $fh$ are in
$\bar{{\cal H}}_G^{n-1}$ and $\bar{{\cal H}}_G^{n+1}$
respectively, and
their norms satisfy the inequality
$\|f\underset{ k}{\times} h\|\le\|f\|\cdot\|h\|$ and
$\|fh\|\le\|f\|\cdot\|h\|$. The relation
$$
n!I_G(f)I_G(h)=(n+1)!I_G(fh)+\sum_{k=1}^n (n-1)!
I_G(f\underset{k}{\times}h)
$$
holds true.}

\medskip
Proposition 4.2 is a simple consequence of the diagram formula
if we apply it for the product of the Wiener--It\^o integrals
$n!I_G(f)$ and $I_G(h)$ with the functions
 $f\in\bar{{\cal H}}_G^n$ and $h\in\bar{{\cal H}}_G^1$.
We have to observe that in this case such diagrams appear which
have two rows, the first row containing the vertices $(1,1)$, 
$(2,1)$,\dots, $(n,1)$, and the second row having one vertex
$(1,2)$. There are two kind of diagrams in this model. The first 
kind of diagrams contains no edge, and it gives the kernel 
function~$fh$. The other kind of diagrams contains one edge 
connecting the vertices $(k,1)$ an $(1,2)$, $1\le k\le n$, 
giving the kernel function $f\underset{k}{\times}h$. In the
last step we exploited that $h(-x)=\overline{h(x)}$, because 
$h\in{\cal H}_G^1$. These observations imply Proposition~4.2.

\medskip
The other result we need is the following (well-known) 
recursion formula for Hermite polynomials.

\medskip\noindent
{\bf Lemma 4.3.} {\it The identity
$$
H_n(x)=xH_{n-1}(x)-(n-1)H_{n-2}(x) \quad\textrm{for \ }n=1,2,\dots,
$$
holds with the notation $H_{-1}(x)\equiv0$.}

\medskip\noindent
{\it Proof of Lemma 4.3.}
\begin{eqnarray*}
H_n(x)=(-1)^ne^{x^2/2}\frac{d^n}{dx^n}(e^{-x^2/2})
&=&-e^{x^2/2}\frac d{dx}\left(H_{n-1}(x)e^{-x^2/2}\right)\\
&=&xH_{n-1}(x)-\frac d{dx}H_{n-1}(x).
\end{eqnarray*}
Since $\frac d{dx}H_{n-1}(x)$ is a polynomial of order $n-2$ 
with leading coefficient $n-1$ we can write
$$
\frac d{dx} H_{n-1}(x)=(n-1)H_{n-2}(x)+\sum_{j=0}^{n-3}c_jH_j(x).
$$
To complete the proof of Lemma~4.3 it remains to show that 
in the last expansion all coefficients~$c_j$ are zero. This 
follows from the fact that $e^{-x^2/2}H_{n-1}(x)$ is orthogonal
to any polynomial whose order is not larger than $n-2e^{-x^2/2}H_{n-1}(x)
$ and 
the calculation
\begin{eqnarray*}
\int e^{-x^2/2}H_j(x)\frac d{dx}H_{n-1}(x)\,dx
&=&-\int H_{n-1}(x)\frac d{dx}(e^{-x^2/2}H_j(x))\,dx\\
&=&\int 
e^{-x^2/2}H_{n-1}(x)
P_{j+1}(x)\,dx=0
\end{eqnarray*}
with the polynomial $P_{j+1}(x)=xH_j(x)-\frac d{dx}H_j(x)$ of
order~$j+1$ for $j\le n-3$. 

\medskip\noindent
{\it The proof of It\^o's formula.}\/ We prove 
It\^o's formula by induction with respect to~$N$. It holds 
for $N=1$. Assume that it holds for~$N-1$. Let us define 
the functions
\begin{eqnarray*}
f(x_1,\dots,x_{N-1})&=&g_1(x_1)\cdots g_{N-1}(x_{N-1})\\
h(x)&=&g_N(x).
\end{eqnarray*}
Then
\begin{eqnarray*}
J&=&\int g_1(x_1)\cdots g_N(x_N)Z_G(\,dx_1)\dots Z_G(\,dx_N)\\
&=&N!\,I_G(fh)=(N-1)!\,I_G(f)I_G(h)-\sum_{k=1}^{N-1} (N-2)!\,
I_G(f\underset{k}{\times} h)
\end{eqnarray*}
by Proposition~4.1. We can write because of our induction 
hypothesis that
\begin{eqnarray*}
J&=&H_{j_1}\left(\int\varphi_1(x)Z_G(\,dx)\right)\cdots
H_{j_{m-1}}\left(\int\varphi_{m-1}(x)Z_G(\,dx)\right) \\
&&\qquad\qquad\qquad H_{j_m-1}\left(\int\varphi_{m}(x)Z_G(\,dx)\right)
\int\varphi_m(x)Z_G(\,dx)\\
&&\qquad-(j_m-1)H_{j_1}\left(\int\varphi_1(x)Z_G(\,dx)\right)\cdots
H_{j_{m-1}}\left(\int\varphi_{m-1}(x)Z_G(\,dx)\right)\\
&&\qquad\qquad\qquad H_{j_m-2}\left(\int\varphi_m(x)Z_G(\,dx)\right),
\end{eqnarray*}
where $H_{j_m-2}(x)=H_{-1}(x)\equiv0$ if $j_m=1$. This relation
holds, since
\begin{eqnarray*}
&&
\!\!\!\!\!\!\! \!\!\!\!\!\! \!\!
f\underset{k}{\times} h(x_1,\dots,x_{k-1},x_{k+1},\dots,x_{N-1})
=\int g_1(x_1)\cdots g_{N-1}(x_{N-1})
\overline{\varphi_m(x_k)}G(\,dx_k)\\
&&
\!\!\!\!\!\!\! \!\!\!\!\!\! \!\!\! \!\!
\quad \,\,=\left\{   
\begin{array}{l}
0 \quad \textrm{if \ } k\le N-j_m\\
g_1(x_1)\cdots g_{k-1}(x_{k-1})g_{k+1}(x_{k+1})\cdots 
g_{N-1}(x_{N-1}) \quad \textrm{if } N-j_m<k\le N-1.
\end{array} \right.
\end{eqnarray*}
Hence Lemma~4.3 implies that
\begin{eqnarray*}
J&=&\prod_{s=1}^{m-1}H_{j_s}
\left(\int \varphi_s(x)Z_G(\,dx)\right)
\biggl[ H_{j_m-1}\left(\int\varphi_m(x)Z_G(\,dx)\right)
\int\varphi_m(x)Z_G(\,dx) \\
&&\qquad -(j_m-1)H_{j_m-2}\left(\int\varphi_m(x)
Z_G(\,dx)\right)\biggr]
=\prod_{s=1}^m H_{j_s}\left(\int\varphi_s(x)Z_G(\,dx)\right),
\end{eqnarray*}
as claimed. 

\Section{Some applications of the diagram formula}

In this section I discuss two kinds of results related 
to the diagram formula. The first of them is about
the description of subordinated random fields and the 
construction of non-Gaussian self-similar random fields. 
I shall formulate the results both for discrete and 
generalized random fields. I shall omit some proofs, in 
particular the proof of the results about generalized random 
fields which are related to the theory of generalized 
function, a circle of problems not discussed in this
note. The other problem I shall discuss is about the
estimation of high moments and the tail distribution 
of Wiener--It\^o integrals by means of the diagonal formula.

\subsection{Description of subordinated random fields, 
construction of self-similar random fields}

First I deal with the description of subordinated random 
fields.

Let $X_n$, $n\in{\mathbb Z}_\nu$, be a discrete stationary 
Gaussian random field with a non-atomic spectral measure, 
and let the random field $\xi_n$, $n\in{\mathbb Z}_\nu$, be 
subordinated to it. Let $Z_G$ denote the random spectral 
measure adapted to the random field $X_n$. By Theorem~3.4 
the random variable $\xi_0$ can be represented as
$$
\xi_0=f_0+\sum_{k=1}^\infty \frac1{k!}\int
f_k(x_1,\dots,x_k)Z_G(\,dx_1)\dots Z_G(\,dx_k)
$$
with an appropriate  function 
$f=(f_0,f_1,\dots)\in\,\textrm{Exp}\,{\cal H}_G$ in a unique way. 
This formula together with Theorem~4.4 yields the following

\medskip\noindent
{\bf Theorem 5.1.} {\it A random field $\xi_n$, $n\in{\mathbb Z}_\nu$,
subordinated to the stationary Gaussian random field $X_n$,
$n\in{\mathbb Z}_\nu$, with non-atomic spectral measure can be 
written in the form
\begin{equation}
\xi_n=f_0+\sum_{k=1}^\infty \frac1{k!}\int e^{i((n,x_1+\cdots+x_k)}
f_k(x_1,\dots,x_k)Z_G(\,dx_1)\dots Z_G(\,dx_k), 
\quad n\in{\mathbb Z}_\nu, \label{(6.1)}
\end{equation}
with some $f=(f_0,f_1,\dots)\in\,\textrm{\rm Exp}\,{\cal H}_G$,
where $Z_G$ is the random spectral measure adapted to the random 
field~$X_n$. This representation is unique. On the other hand, 
formula~(\ref{(6.1)}) defines a subordinated field for all 
$f\in\,\textrm{\rm Exp}\,{\cal H}_G$.}

\medskip\noindent
We rewrite formula (\ref{(6.1)}) in a slightly different form
that shows the similarity between Theorem~5.1 and its analogue,
Theorem~5.2 that gives a representation of subordinated generalized
fields.

Let $G$ denote the spectral measure of the underlying stationary
Gaussian random field. If it has the property $G(\{x\colon\;x_p=u\})=0$
for all $u\in R^1$ and $1\le p\le\nu$, where $x=(x_1,\dots,x_\nu)$
(this is a strengthened form of the non-atomic property of~$G$), 
then the functions
$$
\bar f_k(x_1,\dots,x_k)=f_k(x_1,\dots,x_k)
\tilde\chi_0^{-1}(x_1+\cdots+x_k),
\quad k=1,2,\dots,
$$
are meaningful, as functions in the measure space
$(R^{k\nu},{\cal B}^{k\nu},G^k)$, where
$\tilde\chi_n(x)=e^{i(n,x)}
\prod\limits_{p=0}^\nu\frac{e^{ix^{(p)}}-1}{ix^{(p)}}$,
$n\in{\mathbb Z}_\nu$, denotes the Fourier transform of 
the indicator function of the $\nu$-dimensional unit cube
$\prod\limits_{p=1}^\nu[n^{(p)},n^{(p)}+1]$. Then the random
variable $\xi_n$ in formula~(\ref{(6.1)}) can be rewritten 
in the form
\begin{equation}
\xi_n=f_0+\sum_{k=1}^\infty \frac1{k!}\int\tilde\chi_n(x_1+\cdots+x_k)
\bar f_k(x_1,\dots,x_k)Z_G(\,dx_1)\dots Z_G(\,dx_k),\quad n\in{\mathbb Z}_\nu.
\label{(6.1a)}
\end{equation}
Hence the following Theorem~5.2 can be considered as the continuous
time version of Theorem~5.1.

\medskip\noindent
{\bf Theorem~5.2.} {\it Let a generalized random field 
$\xi(\varphi)$, $\varphi\in{\cal S}$, be subordinated to a 
stationary Gaussian generalized random field $X(\varphi)$, 
$\varphi\in{\cal S}$. Let $G$ denote the spectral measure 
of the field $X(\varphi)$, and let $Z_G$ be the random 
spectral measure adapted to it. Let the spectral measure~$G$
be non-atomic. Then $\xi(\varphi)$ can be written in the form
\begin{equation}
\xi(\varphi)=f_0\cdot\tilde\varphi(0)+\sum_{k=1}^\infty\frac1{k!}
\int\tilde\varphi(x_1+\cdots+x_k)f_k(x_1,\dots,x_k)Z_G(\,dx_1)
\dots Z_G(\,dx_k), \label{($6.1'$)}
\end{equation}
where the functions $f_k$ are invariant under all permutations 
of their variables,
$$
f_k(-x_1,\dots,-x_k)=\overline{f_k(x_1,\dots,x_k)}, 
\quad k=1,2,\dots,
$$
and
\begin{equation}
\sum_{k=1}^\infty \frac1{k!}\int (1+|x_1+\cdots+x_k|^2)^{-p}
|f_k(x_1+\cdots+x_k)|^2G(\,dx_1)\dots G(\,dx_k)<\infty \label{(6.2)}
\end{equation}
with an appropriate number $p>0$. This representation is unique.

On the other side, all random fields $\xi(\varphi)$, $\varphi\in{\cal S}$,
defined by formulas~(\ref{($6.1'$)}) and (\ref{(6.2)}) 
are subordinated to the stationary, Gaussian random field 
$X(\varphi)$,~$\varphi\in{\cal S}$.}

\medskip
I shall omit the proof of Theorem~5.2. I only make some comments
on it. The proof depends heavily on the theory of generalized 
functions. Even the proof of the statement that formula 
(\ref{($6.1'$)}) defines a (generalized) stationary field is 
not simple. It is not enough to show that
$T_t\xi(\varphi))=\xi(T_t\varphi)$ in this case. We also have to
prove that $E[\xi(\varphi_n)-\xi(\varphi)]^2\to0$ if 
$\varphi_n\to\varphi$ in the topology of $\cal S$, and this demands
some special argument.

But the really hard part of Theorem~5.2 is to show that all 
subordinated fields can be represented in the form of (\ref{($6.1'$)}).
In particular we have to find the kernel functions $f_k$ in this
formula. To find them first we show that there is a function 
$\varphi_0\in{\cal S}$, whose Fourier transform nowhere disappears,
and the linear combinations made with the help of its shifts are 
everywhere dense in ${\cal S}$. Then writing the random variable 
$\xi(\varphi_0)$ in the form (\ref{(4.7)}), and writing the 
functions $f_n(x_1,\dots,x_n)$ in this representation as
$f_n(x_1,\dots,x_n)=\frac{f_n(x_1,\dots,x_n)}
{\tilde\varphi_0(x_1+\cdots+x_n)}\tilde\varphi_0(x_1+\cdots+x_n)$, 
we get that we have to choose the kernel functions
$\frac{f_n(x_1,\dots,x_n)}{\tilde\varphi_0(x_1+\cdots+x_n)}$
in formula~(\ref{($6.1'$)}). A detailed proof of Theorem~5.2
would demand much work, and I omit it.

\medskip
We shall call the representations given in Theorems~5.1 
and~5.2 the canonical representation of a subordinated field.
This notion will play an important role in our investigation 
about limit problems. We shall rewrite the random fields
$Z_n^N$ defined in formula~(\ref{(1.1)}) in the form of their 
canonical representation with the help of It\^o's formula, 
and this helps us to study their limit behaviour. 

From now on we restrict ourselves to the case $E\xi_n=0$ or 
$E\xi(\varphi)=0$ respectively, i.e. to the case when $f_0=0$ 
in the canonical representation. Next I construct self-similar
stationary random fields. To find such fields observe that if
$$
\xi(\varphi)=\sum_{k=1}^\infty\frac1{k!}\int\tilde\varphi(x_1+\cdots+x_k)
f_k(x_1,\dots,x_k)Z_G(\,dx_1)\dots Z_G(\,dx_k),
$$
then
$$
\xi(\varphi^A_t)=\sum_{k=1}^\infty\frac1{k!}\frac{t^\nu}{A(t)}
\int\tilde\varphi(t(x_1+\cdots+x_k))
f_k(x_1,\dots,x_k)Z_G(\,dx_1)\dots Z_G(\,dx_k)
$$
with the function $\varphi^A_t$ defined in~(\ref{(1.3)}), where 
we apply the function $A(t)>0$ appearing in that formula. 
Define the spectral measures $G_t$ by the formula 
$G_t(A)=G(tA)$ for all sets~$A$. Then it is not difficult to 
see looking first at the definition of  Wiener--It\^o 
integrals when it is applied for elementary functions and 
then taking limit in the general case that 
$$
\xi(\varphi^A_t)\stackrel{\Delta}{=}\sum_{k=1}^\infty\frac1{k!}
\frac{t^\nu}{A(t)}
\int\tilde\varphi(x_1+\cdots+x_k)
f_k\left(\frac{x_1}t,\dots,\frac{x_k}t\right)
Z_{G_t}(\,dx_1)\dots Z_{G_t}(\,dx_k).
$$

If the spectral measure $G$ and the kernel functions $f_k$
in the formula expressing $\xi(\varphi)$ have the homogeneity
properties that $G(tB)=t^{2\kappa}G(B)$ with some $\kappa>0$ 
for all $t>0$ and $B\in{\cal B}^\nu$, and the identity 
$f_k(\lambda x_1,\dots,\lambda x_k)=
\lambda^{\nu-\kappa k-\alpha}f_k(x_1,\dots,x_k)$ holds, and 
$A(t)$ is chosen as $A(t)=t^\alpha$, then Theorem~3.7 
(with the choice $G'(B)=G(tB)=t^{2\kappa}G(B)$) implies 
that $\xi(\varphi_t^A)\stackrel{\Delta}{=}\xi(\varphi)$.
Hence we obtain the following

\medskip\noindent
{\bf Theorem 5.3.} {\it Let a generalized random 
field $\xi(\varphi)$ be given by the formula
\begin{equation}
\xi(\varphi)=\sum_{k=1}^\infty\frac1{k!}\int \tilde\varphi
(x_1+\cdots+x_k)f_k(x_1,\dots,x_k)Z_G(\,dx_1)\dots Z_G(\,dx_k).
\label{(6.5)}
\end{equation}
If $f_k(\lambda x_1,\dots,\lambda x_k)=\lambda^{\nu-\kappa k-\alpha}
f_k(x_1,\dots,x_k)$ for all~$k$, $(x_1,\dots,x_k)\in R^{k\nu}$ and
$\lambda>0$, $G(\lambda A)=\lambda^{2\kappa}G(A)$ for all
$\lambda>0$ and $A\in{\cal B}^\nu$, then $\xi$ is a self-similar
random field with parameter $\alpha$.}

\medskip
The discrete time version of this result can be proved in the same
way. It states the following

\medskip\noindent
{\bf Theorem 5.4.} {\it If a discrete random field $\xi_n$,
$n\in{\mathbb Z}_\nu$, has the form
\begin{equation}
\xi_n=\sum_{k=1}^\infty \frac1{k!}\int\tilde\chi_n(x_1+\cdots+x_k)
f_k(x_1,\dots,x_k)Z_G(\,dx_1)\dots Z_G(\,dx_k), \quad
n\in{\mathbb Z}_\nu, \label{($6.5'$)}
\end{equation}
and $f_k(\lambda x_1,\dots,\lambda x_k)=\lambda^{\nu-\kappa k-\alpha}
f_k(x_1,\dots,x_k)$ for all~$k$, $G(\lambda A)=\lambda^{2\kappa}G(A)$,
then $\xi_n$ is a self-similar random field with parameter~$\alpha$.}

\medskip
Theorems~5.3 and~5.4 enable us to construct self-similar random
fields. Nevertheless, we have to check whether formulas~(\ref{(6.5)})
and~(\ref{($6.5'$)})  are meaningful. The hard part of this problem is
to check whether the inequality
$$
\sum \frac1{k!}\int |\tilde\chi_n(x_1+\cdots+x_k)|^2
|f_k(x_1,\dots,x_k)|^2 G(\,dx_1)\dots G(\,dx_k)<\infty
$$
holds in the discrete parameter case or whether the inequality
$$
\sum \frac1{k!}\int |\tilde\varphi(x_1+\cdots+x_k)|^2
|f_k(x_1,\dots,x_k)|^2 G(\,dx_1)\dots G(\,dx_k)<\infty
\quad \textrm{for all }\varphi\in{\cal S}
$$
holds in the generalized field case.

It is a rather hard problem to decide when these 
expressions are finite. This is a hard question even if 
we consider a single integral and not an infinite sum of 
integrals. One may consider the question whether an 
integral is convergent or divergent a technical problem, 
but one should not underestimate it. The question whether 
some integrals are convergent or divergent is closely
related to the problem whether in a certain model we have
a new type of limit theorem with a non-standard normalization
and a new non-Gaussian limit or the classical central
limit theorem holds in that model. In the next result I 
prove such a result about the finiteness of a certain 
integral which is needed to guarantee the existence of 
an important self-similar field. This self-similar field 
will appear in the limit theorems we shall prove.

Let us define the measure~$G$
\begin{equation}
G(A)=\int_A |x|^{2\kappa-\nu}a\left(\frac x{|x|}\right)\,dx,
\quad A\in{\cal B}^\nu,
\label{(6.6)}
\end{equation}
where $a(\cdot)$ is a non-negative, measurable and even 
function on the $\nu$-dimensional unit sphere $S_{\nu-1}$, 
and $\kappa>0$. (The condition $\kappa>0$ is imposed to 
guarantee the relation $G(A)<\infty$ for all bounded sets 
$A\in{\cal B}^\nu$.) We prove the following

\medskip\noindent
{\bf Proposition 5.5.} {\it Let the measure $G$ be defined 
in formula~(\ref{(6.6)}).

\medskip
If the function $a(\cdot)$ is bounded on the 
unit sphere $S_{\nu-1}$, and $\frac \nu k>2\kappa>0$, then
$$
D(n)=\int|\tilde\chi_n(x_1+\cdots+x_k)|^2 
G(\,dx_1)\dots G(\,dx_k)<\infty
\quad \textrm{for all \ } n\in{\mathbb Z}_\nu,
$$
and
\begin{eqnarray*}
D(\varphi)&=&\int|\tilde\varphi(x_1+\cdots+x_k)|^2
G(\,dx_1)\dots G(\,dx_k)\\
&\le& C\int(1+|x_1+\cdots+x_k)|^2)^{-p} 
G(\,dx_1)\dots G(\,dx_k)<\infty
\end{eqnarray*}
for all $\varphi\in{\cal S}$ and $p>\frac\nu2$ with some
$C=C(\varphi,p)<\infty$.}

\medskip\noindent
{\it Remark.}\/ In the lecture note which is the basis of these
lectures I also proved that this result is sharp. Namely, 
if the function $a(\cdot)$ in the definition of the spectral 
measure $G$ is always larger than some number $\varepsilon>0$, 
and $2\kappa\le0$ or $2\kappa\ge\frac\nu k$, then the integrals 
defining $D(n)$ and $D(\varphi)$ are divergent. This means
that the conditions imposed on $\kappa$ in Proposition~5.5
cannot be improved. Actually, the condition about the 
property I imposed on $a(\cdot)$ can be weakened in this statement.

\medskip\noindent
{\it Proof of Proposition 5.5.}\/  
We may assume that $a(x)=1$ for all $x\in S_{\nu-1}$. Define
$$
J_{\kappa,k}(x)=\int_{x_1+\cdots+x_k=x}
|x_1|^{2\kappa-\nu}\cdots|x_k|^{2\kappa-\nu}
\,dx_1\dots\,dx_k, \quad x\in R^{\nu},
$$
for $k\ge2$, where $\,dx_1\dots\,dx_k$ denotes the Lebesgue
measure on the hyperplane $x_1+\cdots+x_k=x$, and let
$J_{\kappa,1}(x)=|x|^{2\kappa-\nu}$. The identity
$$
J_{\kappa,k}(\lambda x)
=|\lambda|^{k(2\kappa-\nu)+(k-1)\nu}J_{\kappa,k}(x)
=|\lambda|^{2k\kappa-\nu}J_{\kappa,k}(x),
\quad x\in R^\nu,\;\; \lambda>0,
$$
holds because of the homogeneity of the integral
(provided that the
integral $J_{\kappa, k}(x)$ is finite). We can write, because
of (\ref{(6.6)}) with $a(x)\equiv1$
\begin{equation}
D(n)=\int_{R^\nu}|\tilde\chi_n(x)|^2 J_{\kappa,k}(x)\,dx, 
\label{(6.7)}
\end{equation}
and
$$
D(\varphi)=\int_{R^\nu}|\tilde\varphi(x)|^2J_{\kappa,k}(x)\,dx.
$$
We prove by induction on $k$ that
\begin{equation}
J_{\kappa,k}(x)\le C(\kappa,k)|x|^{2\kappa k-\nu} \label{(6.8)}
\end{equation}
with an appropriate constant $C(\kappa,k)<\infty$ if
$\frac\nu k>2\kappa>0$.

Inequality~(\ref{(6.8)}) holds for $k=1$, and we have
$$
J_{\kappa,k}(x)=\int J_{\kappa,k-1}(y)|x-y|^{2\kappa-\nu}\,dy
$$
for $k\ge2$. Hence
\begin{eqnarray*}
J_{\kappa,k}(x)\!\!&\le& \!\! C(\kappa,k-1)\int
|y|^{(2\kappa(k-1)-\nu}|x-y|^{2\kappa-\nu}\,dy\\
\!\!&=& \!\! C(\kappa,k-1)|x|^{2\kappa k-\nu}
\int |y|^{(2\kappa(k-1)-\nu}\left|\frac x{|x|}-y\right|^{2\kappa-\nu}\!\!\!dy
=C(\kappa,k)|x|^{2\kappa k-\nu},
\end{eqnarray*}
since $\int |y|^{(2\kappa(k-1)-\nu}
\left|\frac x{|x|}-y\right|^{2\kappa-\nu}\,dy<\infty$.

The last integral is finite, since its integrand behaves at zero
asymptotically as $C|y|^{2\kappa(k-1)-\nu}$, at the point
$e=\frac x{|x|}\in S_{\nu-1}$ as $C_2|y-e|^{2\kappa-\nu}$ and at
infinity as $C_3|y|^{2\kappa k-2\nu}$. Relations~(\ref{(6.7)}) 
and~(\ref{(6.8)}) imply that
\begin{eqnarray*}
D(n)&\le& C'\int |\tilde\chi_0(x)|^2|x|^{2\kappa k-\nu}\,dx
\le C''\int |x|^{2\kappa k-\nu}\prod_{l=1}^\nu\frac1{1+|x^{(l)}|^2}\,dx\\
&\le& C'''\int_{|x^{(1)}|=\max\limits_{1\le l\le\nu}|x^{(l}|}
|x^{(1)}|^{2\kappa k-\nu}\prod_{l=1}^\nu\frac1{1+|x^{(l)}|^2}\,dx\\
&=&\sum_{p=0}^\infty C'''
\int_{|x^{(1)}|=\max\limits_{1\le l\le\nu}|x^{(l}|,
\;2^p\le |x^{(1)}|<2^{p+1}}
+C'''\int_{|x^{(1)}|=\max\limits_{1\le l\le\nu}|x^{(l}|, |x^{(1)}|<1}.
\end{eqnarray*}
The second term in the last sum can be simply bounded by a constant, 
since $B=\left\{x\colon\; |x^{(1)}|=\max\limits_{1\le l\le\nu}|x^{(l}|, \;
|x^{(1)}|<1\right\}\subset\{x\colon\; |x|\le\sqrt \nu\}$,
and we have \hfill\break
$|x^{(1)}|^{2\kappa k-\nu}\prod\limits_{l=1}^\nu\frac1{1+|x^{(l)}|^2}
\le\textrm{const.}\, |x|^{2\kappa k-\nu}$ on the set~$B$. Hence
$$
D(n)\le C_1\sum_{p=0}^\infty 2^{p(2\kappa k-\nu)}
\left[\int_{-\infty}^\infty\frac 1{1+x^2}\,dx\right]^\nu+C_2<\infty.
$$
We have $|\varphi(x)|\le C(1+|x^2|)^{-p}$ with some $C>0$ and $D>0$
if $\varphi\in{\cal S}$.
The proof of the estimate $D(\varphi)<\infty$ for $\varphi\in{\cal S}$
is similar but simpler.

\medskip
We can prove some similar theorems, but they have smaller importance,
so I omit them. I discuss instead another useful application of
the diagram formula, the estimation of high moments of 
Wiener--it\^o integrals.

\subsection{Moment estimates on Wiener--It\^o integrals}

Next I show that the diagram formula, Theorem 4.1, enables
us to estimate the expectation of a product of 
Wiener--it\^o integrals, in particular the moments of 
a Wiener--It\^o integral. 

By  applying the diagram formula we can rewrite the 
product of Wiener--It\^o integrals as a sum of Wiener 
integrals of different multiplicity. The expected value 
of this sum equals the sum of the expected value of the 
individual terms. On the other hand, each Wiener--It\^o 
integral of multiplicity $n\ge1$ has zero expectation. 
Only the constant terms, i.e. Wiener--It\^o integrals 
of zero multiplicity can have non-zero expectation. 
The constant terms in the diagram formula correspond 
to those diagrams in which there starts an edge from 
each vertex. This makes natural to introduce the notion 
of complete diagrams, defined in the following way.
Let $\bar\Gamma\subset\Gamma$ denote the set of complete 
diagrams, i.e.\ let a diagram $\gamma\in\bar\Gamma$  
if an edge enters in each vertex of~$\gamma$. 

Clearly, we have $EI(h_\gamma)=0$ for all
$\gamma\in\Gamma\setminus\bar\Gamma$, since~(\ref{(4.3)}) 
holds for all $f\in\bar{{\cal H}}_G^n$, $n\ge1$, and if 
$\gamma\in\bar\Gamma$, then $I(h_\gamma)\in\bar{{\cal H}}_G^0$. 
Let $h_\gamma$ denote the value of $I(h_\gamma)$ in this case.
Let us also observe that in part~(A) of Theorem~4.1 we 
gave an upper bound for $|I(h_\gamma)|=\|h_\gamma\|$ if 
$\gamma$ is a closed diagram. These facts imply  
the following

\medskip\noindent
{\bf Proposition 5.6.} {\it For all 
$h_1\in\bar{{\cal H}}_G^{n_1}$,\dots, 
$h_n\in\bar{{\cal H}}_G^{n_m}$
$$
En_1!I_G(h_1)\cdots n_m!I_G(h_m)=
\sum_{\gamma\in\bar\Gamma}h_\gamma.
$$
(The sum on the right-hand side equals zero if $\bar\Gamma$ is
empty.) Besides, we have
$$
|h_\gamma|\le \prod_{j=1}^m\|h_j\| \quad \textrm{for all }\gamma\in\bar\Gamma.
$$
}

\medskip
Proposition 5.6 enables us to give a good estimate on the high
moments of a Wiener--It\^o integral. We may
assume that the kernel function of this integral is a symmetric
function. In the next Corollary I formulate an estimate on the
$2N$-th moment of an $m$-fold Wiener--It\^o integral. The interesting
case is when $N$ is large.

\medskip
{\bf Corollary 5.7.} {\it Let $h\in{\cal H}_G^m$. Then
\begin{eqnarray*}
E\left[(m!I_G(h))^{2N}\right]&\le& C(m,N)\|h\|^{2N}
= C(m,N)(E(m!I_G(h))^2)^N \\
&\le&(2mN-1)(2mN-3)\cdots3\cdot 1 (E(m!I_G(h))^2)^N,
\end{eqnarray*}
where $C(m,N)$ denotes the number of complete diagrams consisting of
$2N$ rows with $m$ elements in each row.}

\medskip\noindent
{\it Proof of Corollary~5.7.} The first inequality in Corollary~5.7
follows immediately from Proposition~5.6 if we apply it to the 
$2N$-fold product of the Wiener--It\^o integral $I_G(h)$ with
itself. To prove the next identity it is enough to observe that
$$
E (m!I_G(h))^2=m!\|h\|^2 \quad \textrm{if \ } h\in{\cal H}_G^m.
$$ 
Finally, we have to give an upper bound on the number of complete
diagrams $C(m,N)$. Let us calculate the number of those `generalized'
closed diagrams with the same number of rows $2N$ and $m$ vertices 
in each row, where also one edge starts from each vertex, but an edge 
also may connect vertices from the same row. Then it is not 
difficult to see that he number of such `generalized' closed 
diagrams is $(2mN-1)(2nN-3)\cdots3\cdot 1$, and this is an upper 
bound for $C(m,N)$.

\medskip
Next I formulate some results which can be considered as a consequence
of the above statements. I shall not work out the details of the
proofs. Finally I make some comments about the content of these
results. 

First I formulate the following

\medskip\noindent
{\bf Theorem 5.8.} {\it Let $(\xi_1,\dots,\xi_k)$ be a normal 
random vector, and $P(x_1,\dots,x_k)$ a polynomial of degree~$m$. 
Then
$$
E\left[P(\xi_1,\dots,\xi_k)^{2N}\right]\le Cm,N)(m+1)^N
\left(EP(\xi_1,\dots,\xi_k)^2\right)^N
$$
with the constant $C(m,N)$ introduced in Corollary~5.7.}

\medskip
I omit the proof of Theorem~5.8, I only explain its main idea.
The random variable $P(\xi_1,\dots,\xi_k)$ can be expressed as 
the sum of $j$-fold Wiener--It\^o integrals with $0\le j\le m$.
The moments of each integral can be bounded by means of 
Corollary~5.7. For $j<m$ we have a better estimate than for
$j=m$. A careful analysis provides the proof of Theorem~5.8.

The next result gives an interesting estimate on the tail-distribution
of Wiener--It\^o integrals.

\medskip\noindent
{\bf Theorem 5.9.} {\it Let $G$ be a non-atomic spectral measure 
and $Z_G$ a random spectral measure corresponding to~$G$. For all 
$h\in{\cal H}_G^m$ there exist some constants $K_1>K_2>0$ and 
$x_0>0$ depending on the function~$h$ such that
$$
e^{-K_1x^{2/m}}\le P(|I_G(h)|>x)\le e^{-K_2x^{2/m}}
$$
for all $x>x_0$.}

\medskip\noindent
{\it Remark 1.}\/ As the proof of Theorem~5.9 shows the constant~$K_2$
in the upper bound of the above estimate can be chosen as
$K_m=C_m (EI_G(h)^2)^{-1/m}$ with a constant~$C_m$ depending only on
the order~$m$ of the Wiener--It\^o integral of~$I_G(h)$. This means
that for a fixed number~$m$ the constant~$K_2$ in the above estimate
can be chosen as a constant depending only on the variance of the
random variable~$I_G(h)$. On the other hand, no simple
characterization of the constant~$K_1>0$ appearing in the lower
bound of this estimate is possible.

\medskip\noindent
{\it Remark 2.}\/ Theorem 5.9 has some interesting consequences.
For instance, we know that all bounded random variables in the 
space ${\cal H}$ can be written as a sum of Wiener--It\^o integrals.
Theorem~5.9 implies that this representation of a bounded random
variables must be an infinite sum, since the tail distribution of
a finite sum tends too slowly to zero at infinity.

\medskip
The proof of the lower bound in Theorem 5.9 requires a special
argument that I omit. On the other hand, the upper bound follows
from Corollary 5.7 and the Markov inequality.

\medskip\noindent
{\it Proof of the  upper estimate in Theorem 5.9.}\/ 
By the Markov inequality
$$
P(|I_G(h)|>x)\le x^{2N}E(I_G(h)|^{2N}).
$$
On the other hand, by  Corollary~5.7
$$
E(I_G(h)|^{2N})\le\frac{(2mN-1)(2mN-3)\cdots3\cdot 1}{(m!)^N}(EI_G(h)^2)^N.  
$$

We get, by multiplying the inequalities
$$
(2Nm-2j-1)(2Nm-2j-1-2N)\cdots (2Nm-2j-1-2N(m-1))\le (2N)^mm!,
$$
for all $j=1,\dots,N$ that
$$
\frac{(2mN-1)(2mN-3)\cdots3\cdot 1}{(m!)^N} \le (2N)^{mN}.
$$
(This inequality could be sharpened, but it is sufficient for our
purpose.) Choose a sufficiently small number $\alpha>0$, and define
$N=[\alpha x^{2/m}]$, where $[\cdot]$ denotes integer part. With
this choice we have
$$
P(|I_G(h)|>x)\le(x^{-2}(2\alpha)^mx^2)^N (EI_G(h)^2)^N
=[(2\alpha)^m EI_G(h)^2]^N
\le e^{-K_2x^{2/m}}
$$
if $\alpha$ is chosen in such a way that 
$(2\alpha)^m E(I_G(h)^2\le\frac1e$,
$K_2=\frac\alpha2$, and $x>x_0$ with an appropriate $x_0>0$.

\medskip
Observe that if $\xi$ is a standard normal variable, 
then $P(|\xi|^m>x)=P(|\xi|>x^{1/m})<e^{-x^{2/m}}$ for $x>1$, 
and this estimate is sharp. Thus Theorem~5.9 means 
that an $m$-fold Wiener--It\^o integral, i.e. a random 
variable $\eta\in{\cal H}_m$ has a similar tail 
distribution behaviour as the $m$-th power of a normal 
random variable with expectation zero. This shows a new 
property of the decomposition of the Hilbert space 
${\cal H}$ (consisting of the square integrable random 
variables measurable with respect to the $\sigma$-algebra 
generated by the underlying Gaussian random field.)

In such a way we got a different characterization of the
space of random variables which can be written down as an
$m$-fold Wiener--It\^o integral with a fixed number~$m$.
Let me also remark that Theorem~5.9 is closely related to
the previous results of this subsection, since the growth 
behaviour of the high moments of a random variable and the
behaviour of its tail distribution are closely related.

If we want to study the high moments (or the behaviour
of the tail distribution) of the random variables in 
${\cal H}$, then we can do this also by means of the theory 
of the original Wiener--it\^o integrals and the diagram 
formula about the calculation of their products. In the
theory of the original Wiener--It\^o integrals we 
integrate with respect to Gaussian orthogonal measures. 
This has some advantages. It is simpler, and moreover
it has some useful modifications. We can work out the
theory of multiple integrals with respect to such random
measures which are not Gaussian, but they preserve that
property that in these random measures the measure
of disjoint sets are independent random variables. 
Also a version of the diagram formula can be proved for 
such random integrals, and this result has some useful 
applications. 

Thus for instance I could prove good estimates on the
moments and on the behaviour of the tail distribution of  
$U$-statistics by applying the diagram formula for a
version of the diagram formula for an appropriately
defined random measure. These results played a very 
important role in my Springer Lecture Note {\it On the 
estimation of multiple random integrals and degenerate 
$U$-statistics}.

\medskip
In this subsection I explained how to get good estimates
on non-linear functionals of Gaussian random fields with
the help of Wiener--It\^o integrals. I would remark that
there is another powerful method to deal with such problems.
This is the theory of logarithmic Soboliev inequalities
worked out by E. Nelson and L.~Gross. I would refer to
the paper of L. Gross {\it Logarithmic Soboliev inequalities}
Am. J. Math. {\bf 97} (1061--1083) (1975) about this
subject. The theory of logarithmic Soboliev inequality
is based on a theory completely different from the subject
of these lectures (on the theory of Markov processes), 
so I do not discuss it here.

\Section{Some limit theorems for non-linear functionals 
of Gaussian random fields}

In this section I give the proof of some non-trivial limit
theorems about the limit behaviour of the renormalizations 
$Z_n^N$, $N=1,2,\dots$, $n\in{\mathbb Z}_\nu$, of a 
stationary random field $\xi_n$, $n\in{\mathbb Z}_\nu$, 
defined in formula~(\ref{(1.1)}) if the underlying random 
field $\xi_n$, $n\in{\mathbb Z}_\nu$, has some nice 
properties. In the first subsection I formulate the 
main results, and introduce the notions needed to
formulate them. Here I also explain the main ideas of the 
proofs which also indicate what kind of results we can 
expect. In the next subsection I present the details
of the proof. In that part I copy the original proofs from 
my lecture note with almost no change, only with some 
rearrangement of the text. Finally I discuss the content 
of our results and their relation to some other problems.

\subsection{Formulation of the main results}

To formulate our results first we have to introduce some notions.
First I recall the definition of slowly varying functions, an
important notion also in the theory of limit theorems for sums
of independent variables, and also formulate the Karamata
theorem that gives a useful characterization of them.

\medskip\noindent
{\bf Definition 6A. (Definition of Slowly Varying Functions.)} 
{\rm A function $L(t)$, $t\in[t_0,\infty)$, $t_0>0$, is 
said to be a slowly varying function (at infinity) if
$$
\lim_{t\to\infty}\frac{L(st)}{L(t)}=1 
\quad\textrm{for all \ } s>0.
$$
}

\medskip
We shall apply the following description of slowly 
varying functions.

\medskip\noindent
{\bf Theorem 6A. (Karamata's Theorem.)} 
{\it If a slowly varying function $L(t)$, $t\ge t_0$,
with some $t_0>0$ is bounded on every 
finite interval, then it can be represented in the form
$$
L(t)=a(t)\exp\left\{\int_{t_0}^t 
\frac{\varepsilon(s)}s\,ds\right\},
$$
where $a(t)\to a_0\neq0$, $\varepsilon(t)$ is integrable on any
finite interval $[t_0,t]$, and $\varepsilon(t)\to0$ as 
$t\to\infty$.} 

\medskip
We shall consider a stationary Gaussian random field 
$X_n$, $n\in\,{\mathbb Z}_\nu$, with expectation zero 
and  correlation function
\begin{equation}
r(n)=EX_0X_n=|n|^{-\alpha}a\left(\frac n{|n|}\right)L(|n|),
\quad n\in{\mathbb Z}_\nu, \quad \textrm{if }n\neq(0,\dots,0)), 
\label{(8.1)}
\end{equation}
where $0<\alpha<\nu$, $L(t)$, $t\ge1$, is a slowly varying 
function, bounded in all finite intervals, and 
$a(t)$ is a continuous function on the unit sphere 
${\cal S}_{\nu-1}$, satisfying the symmetry property 
$a(x)=a(-x)$ for all $x\in{\cal S}_{\nu-1}$. We shall
prove limit theorems for such random fields which are
subordinated to this Gaussian random field 
$X_n$, $n\in\,{\mathbb Z}_\nu$,

Let $G$ denote the spectral measure of the field~$X_n$, 
and let us define the measures~$G_N$, $N=1,2,\dots$, 
by the formula
\begin{equation}
G_N(A)=\frac{N^\alpha}{L(N)}G\left(\frac AN\right),
\quad A\in{\cal B}^\nu, \quad N=1,2,\dots. \label{(8.2)}
\end{equation}
To get our results we shall need a result about the
asymptotic behaviour of the rescaled
versions $G_N$ of the spectral measure~$G$. To formulate
this result whose proof I postpone to the next subsection
we have to introduce the notion of vague convergence of 
not necessarily finite measures on a Euclidean space.
This notion is a natural counterpart of the notion of
weak convergence, an important notion in probability
theory if we work with not necessarily finite measures.

\medskip\noindent
{\bf Definition of Vague Convergence of Measures.} 
{\it Let $G_n$, $n=1,2,\dots$, be a sequence of locally 
finite measures over $R^\nu$, i.e. let $G_n(A)<\infty$ 
for all measurable bounded sets~$A$. We say that the 
sequence $G_n$ vaguely converges to a locally finite
measure~$G_0$ on~$R^\nu$ (in notation 
$G_n\stackrel{v}{\rightarrow} G_0$) if
$$
\lim_{n\to\infty}\int f(x)\,G_n(\,dx)=\int f(x)\,G_0(\,dx)
$$
for all continuous functions~$f$ with a bounded support.}

\medskip
I formulate the following

\medskip\noindent
{\bf Lemma 6.1.} {\it Let $G$ be the spectral measure of 
a stationary random field with a correlation function 
$r(n)$ of the form~(\ref{(8.1)}). Then the sequence of 
measures~$G_N$ defined in~(\ref{(8.2)}) tends vaguely to 
a locally finite measure~$G_0$. The measure~$G_0$ has 
the homogeneity property
\begin{equation}
G_0(A)=t^{-\alpha}G_0(tA) \quad \textrm{for all } A\in{\cal B}^\nu
\quad\textrm{and } t>0, \label{(8.3)}
\end{equation}
and it satisfies the identity
\begin{eqnarray}
&&2^\nu\int e^{i(t,x)}
\prod_{j=1}^\nu\frac{1-\cos x^{(j)}}{(x^{(j)})^2}
\,G_0(\,dx) \label{(8.4)} \\
&&\qquad =\int_{[-1,1]^\nu} (1-|x^{(1)}|)\cdots (1-|x^{(\nu)}|)
\frac{a\left(\frac{x+t}{|x+t|}\right)}{|x+t|^\alpha}\,dx, 
\quad\textrm{for all } t\in R^\nu. \nonumber
\end{eqnarray}
}

\medskip
Formula~(\ref{(8.3)}) together with the vague convergence
of $G_n$ to $G_0$ can be  heuristically so interpreted that 
the measure $G$ is asymptotically homogeneous in the 
neighbourhood of zero. On the other hand, it can be proved 
that we get a correlation function of the form~(\ref{(8.1)}) 
by defining it as the Fourier transform of a (positive)  
measure with a density of the form 
$g(u)=|u|^{\alpha-\nu}b(\frac u{|u|})h(|u|)$, $u\in R^\nu$, 
where $b(\cdot)$ is a non-negative smooth function on the 
unit sphere $\{u\colon\; u\in R^\nu,\,|u|=1\}$, and $h(u)$ 
is a non-negative smooth function on $R^1$ which does not 
disappear at the origin, and tends to zero at infinity 
sufficiently fast. The regularizing function $h(|u|)$ is 
needed in this formula to make the function $g(\cdot)$ 
integrable. Results of this type are studied in the 
theory of generalized functions.

The above mentioned result is interesting for us, because
it shows that there are correlation functions of the 
form (\ref{(8.1)}) with appropriate functions $a(\cdot)$
and $L(\cdot)$. The problem with the definition of a
correlation function $r(n)$, $n\in{\mathbb Z}_\nu$, 
satisfying~(\ref{(8.1)}) is that this function  must
be positive definite. We can guarantee this property by 
defining it as the Fourier transform of a measure.

\medskip

I remark that formulae~(\ref{(8.3)}) and~(\ref{(8.4)}) 
imply that the function~$a(t)$ and number~$\alpha$ in the 
definition~(\ref{(8.1)}) of a correlation function~$r(n)$ 
uniquely determine the measure~$G_0$. Indeed, by 
formula~(\ref{(8.4)}) they determine the (finite) measure
$\prod\limits_{j=1}^\nu\frac{1-\cos x^{(j)}}{(x^{(j)})^2}G_0(\,dx)$, 
since they determine its Fourier transform. Hence they also 
determine the measure~$G_0$. (Formula~(\ref{(8.3)}) shows 
that $G_0$ is a locally finite measure). Let us also remark 
that since $G_N(A)=G_N(-A)$ for all $N=1,2,\dots$, and 
$A\in {\cal B}^\nu$, the relation $G_0(A)=G_0(-A)$,
$A\in{\cal B}^\nu$ also holds. These properties of the 
measure~$G_0$ imply that it can be considered as the 
spectral measure of a generalized random field. 

Now I formulate the basic limit theorem of this section.

\medskip\noindent
{\bf Theorem 6.2.} {\it Let $X_n$, $n\in{\mathbb Z}_\nu$, be 
a stationary Gaussian field with a correlation function 
$r(n)$ satisfying relation~(\ref{(8.1)}) and such that 
$r(0)=EX_n^2=1$, $n\in{\mathbb Z}_\nu$. Let us define the 
stationary random field $\xi_j=H_k(X_j)$, $j\in{\mathbb Z}_\nu$, 
with some positive integer~$k$, where $H_k(x)$ denotes the 
$k$-th Hermite polynomial with leading coefficient~1, and 
assume that the parameter~$\alpha$ appearing in~(\ref{(8.1)}) 
satisfies the relation $0<\alpha<\frac\nu k$ with this 
number~$k$. If the random fields $Z^N_n$, $N=1,2,\dots$, 
$n\in{\mathbb Z}_\nu$, are defined by formula~(\ref{(1.1)}) 
with $A_N=N^{\nu-k\alpha/2}L(N)^{k/2}$ and the above defined 
random variables $\xi_j=H_k(X_j)$, then their multi-dimensional
distributions tend to those of the random field~$Z^*_n$,
$$
Z^*_n=\int \tilde\chi_n(x_1+\cdots+x_k)\,
Z_{G_0}(\,dx_1)\dots Z_{G_0}(\,dx_k), 
\quad n\in{\mathbb Z}_\nu.
$$
Here $Z_{G_0}$ is a random spectral measure corresponding to the
spectral measure $G_0$ which appeared in Lemma~6.1. The function
$\tilde\chi_n(\cdot)$, $n=(n^{(1)},\dots,n^{(\nu)})$, is (similarly
to formula~(\ref{(6.1a)}) in Section~5) the Fourier transform of 
the indicator function of the $\nu$-dimensional unit cube
$\prod\limits_{p=1}^\nu[n^{(p)},n^{(p)}+1]$.}

\medskip
I give a heuristic explanation for this result. First I explain 
why the choice of the normalizing constant~$A_N$ in Theorem~6.2 
was natural, then I explain the main ideas of the proof.
I shall work out the details in the next subsection.

There is a fairly well-known result by which if $(\xi,\eta)$ 
is a Gaussian random vector with $E\xi=E\eta=0$ and 
$E\xi^2=E\eta^2=1$, then they satisfy the identity
$EH_k(\xi)H_k(\eta)=k!(E\xi\eta)^k$. I give a short sketch of
a possible proof.

Put $r=E\xi\eta$. Then we can write $\eta=r\xi+(1-r^2)^{1/2}Z$,
where $Z=(1-r^2)^{-1/2}(\eta-r\xi)$ is a standard Gaussian 
random variable, uncorrelated with, hence also independent 
of the random variable~$\xi$. We can express the random 
variable $H_k(\eta)=H_k(r\xi+(1-r^2)^{1/2}Z)$, as a linear 
combination of the random variables $H_j(\xi)H_l(Z)$ with 
indices $j,l$, $0\le j+l\le k$. In this linear combination the term
$H_k(\xi)$ has coefficient $r^k=(E\xi\eta)^k$. This is the only term
in this linear combination which is not orthogonal to the random variable
$H_k(\xi)$. Hence $EH_k(\xi)H_k(\eta)=(E\xi\eta)^k EH^2_k(\xi)$. On 
the other hand, $EH^2_k(\xi)=k!$ that can be seen with the help of 
some calculation. It follows for instance from It\^o's formula 
and the diagram formula, but it can be proved directly from the 
formula by which we defined the Hermite polynomials by applying 
partial integration.

The above identity implies that
$$
E(Z_n^N)^2=\frac{k!}{A_N^2}\sum_{j,\,l\in B_0^N}r(j-l)^k
\sim\frac{k!}{A_N^2}
\sum_{j,\,l\in B_0^N}|j-l|^{-k\alpha}a^k
\left(\frac{j-l}{|j-l|}\right)L(|j-l|)^k,
$$
with the set $B_0^N$ introduced after formula~(\ref{(1.1)}). Some 
calculation with the help of the above formula shows that with 
our choice of~$A_N$ the expectation $E(Z_n^N)^2$ is separated 
both from zero and infinity, therefore this is the natural 
norming factor. In this calculation we have to exploit the 
condition $k\alpha<\nu$, which implies that in the sum 
expressing $E(Z_n^N)^2$ those terms are dominant for which 
$j-l$ is relatively large, more explicitly which are of 
order~$N$. There are $\textrm{const.}\, N^{2\nu}$ such terms.

The field $\xi_n$, $n\in{\mathbb Z}_\nu$, is subordinated 
to the Gaussian field~$X_n$. It is natural to express its
terms via multiple  Wiener--It\^o integrals, and to
write up the canonical representation of the fields $Z_n^N$
defined in Section~5.

It\^o's formula yields the identity
$$
\xi_j=H_k\left(\int e^{i(j,x)}Z_G(\,dx)\right)
=\int e^{i(j,x_1+\cdots+x_k)}Z_G(\,dx_1)\dots Z_G(\,dx_k),
$$
where $Z_G$ is the random spectral measure adapted to the random
field~$X_n$. Then
\begin{eqnarray*}
Z_n^N&=&\frac1{A_N}\sum_{j\in B_n^N}\int e^{i(j,x_1+\cdots+x_k)}
Z_G(\,dx_1)\dots Z_G(\,dx_k)\\
&=&\frac1{A_N}\int e^{i(Nn,x_1+\cdots+x_k)}\prod_{j=1}^\nu
\frac{e^{iN(x_1^{(j)}+\cdots+x_k^{(j)})}-1}
{e^{i(x_1^{(j)}+\cdots+x_k^{(j)})}-1}\,
Z_G(\,dx_1)\dots Z_G(dx_k).
\end{eqnarray*}
Let us make the substitution $y_j=Nx_j$, $j=1,\dots,k$, in the last
formula, and let us rewrite it in a form resembling 
formula~(\ref{($6.5'$)}).
To this end, let us introduce the measures~$G_N$ defined 
in~(\ref{(8.2)}). It is not difficult to see that
$$
Z_n^N\stackrel{\Delta}{=}\int f_N(y_1,\dots,y_k)
\tilde\chi_n(y_1+\cdots+y_k)\,
Z_{G_N}(\,dy_1)\dots Z_{G_N}(dy_k)
$$
with the measure $G_N$ defined in (\ref{(8.2)}) and
\begin{equation}
f_N(y_1,\dots,y_k)=\prod_{j=1}^\nu
\frac{i(y_1^{(j)}+\cdots+y_k^{(j)})}
{\left(\exp\left\{i\frac1N(y_1^{(j)}
+\cdots+y_k^{(j)})\right\}-1\right)N},
\label{(8.5)}
\end{equation}
where $\tilde\chi_n(\cdot)$ is the Fourier transform of the 
indicator function of the unit cube 
$\prod\limits_{j=1}^\nu[n^{(j)},n^{(j)}+1)$.

(In the above calculations we made a small inaccuracy. We calculated
freely with Wiener--It\^o integrals, although we defined them
only in the case when the spectral measure $G$ is non-atomic.
We shall prove a result in Lemma~6B below which states
that if the correlation function satisfies 
relation~(\ref{(8.1)}), then this property holds. Moreover,
we shall prove a stronger statement, namely that all hyperplanes
$x^{(j)}=t$, $1\le j\le\nu$, $t\in R^1$, have zero $G$ measure.
This fact together Fubini's theorem imply that the set 
where the denominator of the the functions $f_N$ defined in 
formula~(\ref{(8.5)}) disappears, i.e. the set, where 
$y_1^{(j)}+\cdots+y_k^{(j)}=2lN\pi$ with some integer
$l\neq0$ and $1\le j\le\nu$ has zero $G_N\times\cdots\times G_N$ 
measure. This means that the functions $f_N$ are well defined.)

The functions~$f_N$ tend to 1 uniformly in all bounded regions, 
and the measures~$G_N$ tend vaguely to~$G_0$ as $N\to\infty$ 
by Lemma~6.1. These relations suggest the following limiting 
procedure. The limit of $Z_n^N$ can be obtained by 
substituting $f_N$ with~1 and the random spectral measure
$Z_{G_N}$ with~$Z_{G_0}$ in the Wiener--It\^o integral
 expressing~$Z_n^N$. This would provide that the large-scale 
limit of the fields $Z^N_n$ equals the random field $Z^*_n$ 
defined in the formulation of Theorem~6.2.  

We have to justify the above formal limiting procedure. 
We shall do this by showing first that the Wiener--It\^o 
integral expressing~$Z_n^N$ is essentially concentrated in 
a large bounded region independent of~$N$. 
The $L_2$-isomorphism of Wiener--It\^o integrals can help us 
to show this. We will justify this argument in Lemma~6.3.
 
I shall discuss Lemma 6.3 in the next subsection, where
the technical details of the proofs is explained. In this 
section I shall still present the proof of Lemma~6B which 
is, as I mentioned before is needed for the justification 
of some formal steps we made before. Then I finish this
subsection with a result formulated in Theorem~$6.2'$
which is a natural continuation of Theorem~6.2.

\medskip\noindent
{\bf Lemma~6B.} {\it Let the correlation function of a 
stationary random field $X_n$, $n\in{\mathbb Z}_\nu$, 
satisfy the relation $r(n)\le A|n|^{-\alpha}$ with some 
$A>0$ and $\alpha>0$ for all $n\in{\mathbb Z}_\nu$, 
$n\neq0$. Then its spectral measure $G$ is non-atomic. 
Moreover, the hyperplanes $x^{(j)}=t$ have zero $G$ measure
for all $1\le j\le \nu$ and $t\in R^1$.}

\medskip\noindent
{\it Proof of Lemma 6B.} Lemma 6B clearly holds if $\alpha>\nu$,
because in this case the spectral measure~$G$ has even a density
 function $g(x)=\sum\limits_{n\in{\mathbb Z}_\nu}e^{-i(n,x)}r(n)$.
On the other hand, the $p$-fold convolution of the spectral measure
$G$ with itself (on the torus $R^\nu/2\pi{\mathbb Z}_\nu$) 
has Fourier transform, $r(n)^p$, $n\in{\mathbb Z}^\nu$, and as 
a consequence in the case $p>\frac\nu\alpha$ this measure is 
non-atomic. Hence it is enough to show that if the convolution 
$G*G$ is a non-atomic measure, then the measure~$G$ is also
non-atomic. But this is obvious, because if there were a point 
$x\in R^\nu/2\pi{\mathbb Z}_\nu$ such that $G(\{x\})>0$, 
then the relation $G*G(\{x+x\})>0$ would hold, and this is 
a contradiction. (Here addition is taken on the torus.) The 
second statement of the lemma can be proved with some small 
modifications of the previous proof, by reducing it to the 
one-dimensional case. 

I finish this subsection with the following

\medskip\noindent
{\bf Theorem~$6.2'$.} {\it Let $X_n$, $n\in{\mathbb Z}_\nu$, 
be a stationary Gaussian field with a correlation function 
$r(n)$ defined in~(\ref{(8.1)}) and such that $r(0)=EX_n^2=1$,
$n\in{\mathbb Z}_\nu$. Let $H(x)$ be a real valued function with 
the properties $EH(X_n)=0$ and $EH(X_n)^2<\infty$. 
Let us consider the orthogonal expansion
\begin{equation}
H(x)=\sum_{j=1}^\infty c_jH_j(x), \quad \sum c_j^2j!<\infty, 
\label{(8.28)}
\end{equation}
of the function~$H(\cdot)$ by the Hermite polynomials~$H_j$ 
(with leading coefficients~1). Let $k$ be the smallest index 
in this expansion such that $c_k\neq0$. If $0<k\alpha<\nu$ 
for the parameter~$\alpha$ in~(\ref{(8.1)}), and the field 
$Z^N_n$ is defined by the field $\xi_n=H(X_n)$, 
$n\in{\mathbb Z}_\nu$, and formula~(\ref{(1.1)}), then 
the multi-dimensional distributions of the fields $Z_n^N$ 
with $A_N=N^{\nu-k\alpha/2}L(N)^{k/2}$ tend to those of the 
fields $c_kZ^*_n$, $n\in{\mathbb Z}_\nu$, where the 
field $Z^*_n$ is the same as in Theorem~6.2.}

\medskip\noindent
{\it Proof of Theorem~$6.2'$ with the help of Theorem~6.2.}\/ 
Define $H'(x)=\sum\limits_{j=k+1}^\infty c_jH_j(x)$ and
$Y_n^N=\frac1{A_N}\sum\limits_{l\in B_n^N}H'(X_l)$.
Because of Theorem~6.2 in order to prove Theorem~$6.2'$ 
it is enough to show that
$$
E(Y_n^N)^2\to0 \quad\textrm{as } N\to\infty.
$$
It can be proved similarly to the identity 
$EH_k(\xi)H_k(\eta)=k!(E\xi\eta)^k EH^2_k(\xi)$
for a Gaussian vector $(\xi,\eta)$ such that
$E\xi=E\eta=0$, $E\xi_k^2=E\eta_k^2=1$, that
also the identity
$EH_j(\xi)H_l(\eta)=j!\delta_{j,l}(E\xi\eta)^j$
holds, where $\delta_{j,l}=0$ if $j\neq l$, 
and $\delta_{j,l}=1$ if $j=l$. 

This means in our case that
$$
EH_j(X_n)H_l(X_m)=\delta_{j,l}j!(EX_nX_m)^j=\delta_{j,l}j!r(n-m)^j.
$$
Hence
$$
E(Y_n^N)^2=\frac1{A_N^2}\sum_{j=k+1}^\infty  c_j^2j!\,
\sum_{s,t\in B_n^N}[r(s-t)]^j.
$$
Some calculation yields with the help of this identity and 
formula~(\ref{(8.1)}) that
$$
E(Y_n^N)^2=\frac1{A_N^2}\left[O(N^{2\nu-(k+1)\alpha}
L(N)^{k+1})+O(N^\nu)\right]\to0.
$$
(Observe that we imposed the condition $\sum c_j^2j!<\infty$
which is equivalent to the condition $E H(X_n)^2<\infty$.)
Theorem~$6.2'$ is proved.

\medskip
The main difference between Theorem~6.2 and Theorem~$6.2'$ 
is that in Theorem~6.2 we considered random variables 
$\xi_n=H_k(X_n)$, while in Theorem~$6.2'$ $\xi_n=H(\xi_n)$ 
with such a function~$H$ in whose expansion with respect 
to Hermite polynomials the Hermite polynomial $H_k(x)$ was 
the term with the smallest index with a non-zero 
coefficient. We saw that in these two cases a very similar 
result holds. In both theorems we imposed the condition 
$0<k\alpha<\nu$ for the parameter~$\alpha$ in~(\ref{(8.1)}). 
One may ask what happens if this condition is violated. 
I mentioned before the proof of Proposition~5.5 that in 
this case the multiple Wiener--It\^o integral defining the 
limiting field $Z_n^*$, $n\in{\mathbb Z}_\nu$, in these 
theorems does not exists. I shall briefly discuss this 
case by pointing out that in this case the central limit 
theorem holds. But because of lack of time I cannot discuss 
the details of the proof.

\subsection{The details of the proofs}

First I prove a lemma that enables us to prove the 
convergence of Wiener--It\^o integrals under some conditions.

\medskip\noindent
{\bf Lemma~6.3.} {\it Let $G_N$, $N=1,2,\dots$, be a sequence
of non-atomic spectral measures on $R^\nu$ tending vaguely 
to a non-atomic spectral measure~$G_0$. Let a sequence of 
measurable functions $K_N=K_N(x_1,\dots,x_k)$, 
$N=0,1,2,\dots$, be given such that 
$K_N\in\bar{{\cal H}}_{G_N}^k$ for $N=1,2,\dots$. Assume 
further that these functions satisfy the following 
properties: For all $\varepsilon>0$ there exist some 
constants $A=A(\varepsilon)>0$ and $N_0=N_0(\varepsilon)>0$ 
such that the conditions~(a) and~(b) formulated below  
are satisfied. 

\medskip
\begin{description}
\item[(a)] The function $K_0$ is continuous on the set
$B=[-A,A]^{k\nu}$, and $K_N\to K_0$ uniformly on 
the set $B$ as $N\to\infty$. Besides, the hyperplanes 
$x_p=\pm A$ have zero $G_0$~measure for all $1\le p\le\nu$.
\item[(b)] $\int_{R^{k\nu}\setminus B}|K_N(x_1,\dots,x_k)|^2
G_N(\,dx_1)\dots G_N(dx_k)<\frac{\varepsilon^2}{k!}$ 
if $N=0$ or $N\ge N_0$, and 
$K_0(-x_1,\dots,-x_k)=\overline{K_0(x_1,\dots,x_k)}$ 
for all $(x_1,\dots,x_k)\in R^{k\nu}$.
\end{description}

\medskip
Then $K_0\in\bar{{\cal H}}_{G_0}^k$, and
$$
\int K_N(x_1,\dots,x_k)\,Z_{G_N}(dx_1)\dots Z_{G_N}(dx_k) 
\stackrel{{\cal D}}{\rightarrow}
\int K_0(x_1,\dots,x_k)\,Z_{G_0}(dx_1)\dots Z_{G_0}(dx_k)
$$
as $N\to\infty$, where $\stackrel{{\cal D}}{\rightarrow}$ denotes
convergence in distribution.}

\medskip\noindent
{\it Remark.}\/ In my Lecture Note a somewhat more general result
is proved that allows to handle also such cases where the function
$K_0$ may have some discontinuities. There are results whose proof 
demands that more general result, but Theorem~6.2 can be proved 
with the help of this simpler result. I shall present a proof 
simpler than in the Lecture Note. The main difference is that 
in the Lecture Note I exploited that the weak convergence of
probability measures is metrizable for instance by means of 
the so-called Prokhorov metric, and I also applied some of its 
properties. Here I use instead that well-known result that a 
sequence of random variables $U_n$ converge weakly to some 
random variable $U_0$ in the Euclidean space $R^p$ if and 
only if their characteristic functions satisfy the relation 
$\lim\limits_{n\to\infty}E^{i(t,U_n)}=E^{i(t,U_0)}$ for all $t\in R^p$. 
I shall give a complete proof of Lemma~6.3 that applies 
Lemma~3.3 whose proof is given in the Appendix. Actually I apply
a slightly stronger version of this result which also formulates
an additional property of the approximation constructed in the
proof of Lemma~3.3 which is mentioned at the end of the proof.

\medskip\noindent
{\it Proof of Lemma~6.3.}\/ First I show that 
$K_0\in\bar{{\cal H}}_{G_0}^k$. Indeed,
Conditions~(a) and~(b) obviously imply that
$$
\int|K_0(x_1,\dots,x_k)|^2\,G_0(\,dx_1)\dots G_0(\,dx_k)<\infty,
$$
hence $K_0\in\bar{{\cal H}}_{G_0}^k$. 

Let us fix an $\varepsilon>0$, and let us choose some 
$A=A(\varepsilon)>0$ and $N_0=N_0(\varepsilon)>0$ for which 
conditions~(a) and~(b) hold with this~$\varepsilon$. Then
\begin{eqnarray}
&&E\left[\int [1-\chi_B(x_1,\dots,x_k)]K_N(x_1,\dots,x_k)\,
Z_{G_N}(\,dx_1)\dots Z_{G_N}(\,dx_k)\right]^2 \nonumber \\
&&\qquad \le k!\int_{R^{k\nu}\setminus B}|K_N(x_1,\dots,x_k)|^2
G_N(\,dx_1)\dots G_N(\,dx_k)<\varepsilon^2
\label{(8.6)}
\end{eqnarray}
for $N=0$ or $N>N_0$, where $\chi_B$ denotes the indicator functions
of the set~$B$ introduced in the formulation of condition~(a).

Since $B=[-A,A]^{k\nu}$, and $G_N\stackrel{v}{\rightarrow} G_0$,
hence $G_N\times\cdots\times G_N(B)<C(A)$ with an appropriate constant
$C(A)<\infty$ for all $N=0,1,\dots$.  Because of this estimate and
the uniform convergence $K_N\to K_0$ on the set~$B$ we have
\begin{eqnarray}
&&\!\!\!   E\left[\int(K_N(x_1,\dots,x_k)-K_0(x_1,\dots,x_k))
\chi_B(x_1,\dots,x_k)\,
Z_{G_N}(\,dx_1)\dots Z_{G_N}(\,dx_k)\right]^2 \nonumber  \\
&&\quad\le k!\int_B
|K_N(x_1,\dots,x_k)-K_0(x_1,\dots,x_k)|^2\,
G_N(\,dx_1)\dots G_N(\,dx_k)<\varepsilon^2 \label{(8.7)}
\end{eqnarray}
for $N>N_1$ with some $N_1=N_1(A,\varepsilon)$.

With the help of relations~(\ref{(8.6)}) and~(\ref{(8.7)}) 
I reduce the proof of Lemma~6.3 to the proof of the relation 
\begin{eqnarray}
&&\int K_0(x_1,\dots,x_k)\chi_B(x_1,\dots,x_k)\,
Z_{G_N}(\,dx_1)\dots Z_{G_N}(\,dx_k) \nonumber \\
&&\qquad \stackrel{{\cal D}}{\rightarrow}
\int K_0(x_1,\dots,x_k)\chi_B(x_1,\dots,x_k)\,
Z_{G_0}(\,dx_1)\dots Z_{G_0}(\,dx_k).
\label{(8.8)}
\end{eqnarray}
For this goal I introduce the random variables
\begin{eqnarray*}
 T_N&=&\int K_N(x_1,\dots,x_k)\,Z_{G_N}(\,dx_1)\dots Z_{G_N}(\,dx_k), \\
 U_N&=&\int K_N(x_1,\dots,x_k)\chi_B(x_1,\dots,x_k)\,
Z_{G_N}(\,dx_1)\dots Z_{G_N}(\,dx_k), \\
V_N&=&\int K_0(x_1,\dots,x_k))\chi_B(x_1,\dots,x_k)\,
Z_{G_N}(\,dx_1)\dots Z_{G_N}(\,dx_k), \\
W&=&\int K_0(x_1,\dots,x_k))\chi_B(x_1,\dots,x_k)\,
Z_{G_0}(\,dx_1)\dots Z_{G_0}(\,dx_k).
\end{eqnarray*}

By inequality (\ref{(8.6)}) we have for all $t\in R^1$ and $N>N_0$
\begin{eqnarray*}
|E(e^{itT_N}-e^{itU_N})|&\le& E|(1-e^{it(U_n-V_N)})|\le E|(t(T_N-U_N)|\\
&\le& |t| (E(T_N-U_N)^2)^{1/2}\le |t|\varepsilon. 
\end{eqnarray*}
Similarly, 
$|E(e^{itU_N}-e^{itV_N})|\le |t| (E(U_N-V_N)^2)^{1/2}\le |t|\varepsilon$
for all $t\in R^1$ and $N>N_0$ by inequality (\ref{(8.7)}).
Finally, $Ee^{itV_N}\to Ee^{itW}$ for all $t\in R^1$ by (\ref{(8.8)}).
These relations together imply that 
$E|e^{itT_N}-Ee^{itW}|\le C(t)|\varepsilon|$ if $N>N_0(t,\varepsilon)$ 
with some numbers $C(t)$ and $N_0(t,\varepsilon)$. Since this
inequality holds for all $\varepsilon>0$, it implies that
$T_N\stackrel{{\cal D}}{\rightarrow}W$ which agrees with the
statement of Lemma~6.3,

\medskip
We shall prove formula~(\ref{(8.8)}) by showing that
$K_0(x_1,\dots,x_k)\chi_B(x_1,\dots,x_k)$ can be well 
approximated by simple functions from 
$\hat{\bar{{\cal H}}}_{G_0}^k$.

More explicitly, I claim that for all 
$\varepsilon>0$ there exists a simple function 
$f_{\varepsilon}\in\hat{\bar{{\cal H}}}^k_{G_0}$ 
such that
\begin{equation}
E\int( K_0(x_1,\dots,x_k)\chi_B(x_1,\dots,x_k)
-f_{\varepsilon}(x_1,\dots,x_k))^2
G_0(\,dx_1)\dots G_0(\,dx_k)\le\frac{{\varepsilon}^2}{k!}
\label{(8.10)}
\end{equation}
and also
\begin{equation}
E\int( K_0(x_1,\dots,x_k)\chi_B(x_1,\dots,x_k)
-f_{\varepsilon}(x_1,\dots,x_k))^2
G_N(\,dx_1)\dots G_N(\,dx_k)\le\frac{{\varepsilon}^2}{k!}
\label{(8.11)}
\end{equation}
if $N\ge N_0$ with some threshold index
$N_0=N_0(\varepsilon,K_0(\cdot)\chi_B(\cdot))$. Moreover, 
this simple function $f_{\varepsilon}$ can be chosen in 
such a way that it is adapted to such a regular system 
${\cal D}=\{\Delta_j,\;j=\pm1,\dots,\pm M\}$ whose 
elements have boundaries of zero $G_0$ measure, i.e. 
$G_0(\partial\Delta_j)=0$ for all $1\le |j|\le M$.

Relation~(\ref{(8.8)}) can be proved with the help
of the estimates (\ref{(8.10)}), (\ref{(8.11)}) and
the relation~(\ref{(8.12)}) formulated below  similarly 
to the reduction of the proof of Lemma~6.3 to~(\ref{(8.8)}). 
Relation (\ref{(8.12)}) states that the function 
$f_{\varepsilon}$ appearing in formulas~(\ref{(8.10)}) and 
(\ref{(8.11)}) satisfies also the relation
\begin{equation}
\int f_{\varepsilon}(x_1,\dots,x_k)\,Z_{G_N}(dx_1)\dots Z_{G_N}(dx_k)
\stackrel{{\cal D}}{\rightarrow}\int f_{\varepsilon}(x_1,\dots,x_k)
\,Z_{G_0}(dx_1)\dots Z_{G_0}(dx_k)
\label{(8.12)}
\end{equation}
as $N\to\infty$. Formula~(\ref{(8.12)}) is a simplified version
of the relation~(\ref{(8.8)}) where the kernel function
$K_0$ in the integral is replaced by a simple function~$f_\varepsilon$.

We can get the proof of~(\ref{(8.8)}) with the help
of the estimates (\ref{(8.10)}), (\ref{(8.11)}) 
and~(\ref{(8.12)}) similarly to the argument leading to
the proof of the limit relation of Lemma~6.3.

In the proof of~(\ref{(8.12)}) we exploit that the function 
$f_\varepsilon\in\hat{\bar{{\cal H}}}_{G_0}^k$ is adapted 
to such a regular system 
${\cal D}=\{\Delta_j,\;j=\pm1,\dots,\pm M\}$ for which all
$\Delta_j$ has the property $G_0(\partial\Delta_j)=0$. In 
the proof of Lemma~3.3 in the Appendix I showed that the 
function $f_{\varepsilon}$ and the regular system $\cal D$ to 
which it is adapted can be chosen in such a way. Besides, 
the spectral measures $G_N$ have the property 
$G_N\stackrel{v}{\rightarrow}G_0$. Hence the (Gaussian) 
random vectors $(Z_{G_N}(\Delta_j),\;j=\pm1,\dots,\pm M)$ 
converge in distribution to the (Gaussian) random vector
$(Z_{G_0}(\Delta_j)$, $j=\pm1,\dots,\pm M)$ as $N\to\infty$.
This implies that if we put the random variables
$(Z_{G_N}(\Delta_j)$, $j=\pm1,\dots,\pm M)$ to the arguments
of a continuous function of $2M$ variables, then these 
random variables converge to the random variable we obtain 
if we put the random variables
$(Z_{G_0}(\Delta_j)$, $j=\pm1,\dots,\pm M)$ to the arguments
of this function. Formula (\ref{(8.12)}) follows from this 
statement because the random vectors we consider in it can 
be written as an appropriate polynomial of these random vectors. 

\medskip
Relation~(\ref{(8.10)}) follows directly from Lemma~3.3
if we apply it to the function $K_0(\cdot)\chi_B(\cdot)$. 
But we need a stronger version of this result, because we 
want to find such a function $f_{\varepsilon}$ which also 
satisfies relations~(\ref{(8.11)}) and~(\ref{(8.12)}). We have
seen that relation~(\ref{(8.12)}) holds if the approximating
function $f_{\varepsilon}$ has the additional property that  
it is adapted to a regular system ${\cal D}$ consisting of
sets with zero $G_0$ measure.


A more careful analysis shows that a function $f_{\varepsilon}$
with this extra property satisfies not only~(\ref{(8.10)}) 
but also~(\ref{(8.11)}) for $N\ge N_0$ with a sufficiently 
large threshold index~$N_0$. We can get another explanation 
of the estimate~(\ref{(8.11)}) by exploiting that the 
function $h_0(x_1,\dots,x_k)$ defined as
$$
h_0(x_1,\dots,x_k)=
K_0(x_1,\dots,x_k)\chi_B(x_1,\dots,x_k)-f_{\varepsilon}(x_1,\dots,x_k)
$$
is almost everywhere continuous with respect to the measure 
$G_0\times\cdots\times G_0$, and it disappears outside the 
compact set~$B$. It can be shown that the vague convergence 
has similar properties as the weak convergence. In particular, 
the above mentioned almost everywhere continuity implies that
$$
\lim_{N\to\infty}\int h_0(x_1,\dots,x_k)G_N(dx_1)\dots G_N(dx_k)
=\int h_0(x_1,\dots,x_k)G_0(dx_1)\dots G_0(dx_k).
$$
In such a way we can reduce the proof of (\ref{(8.10)})
to the proof of~(\ref{(8.11)}). The proof of Lemma~6.3
is finished.

\medskip
Now I show the proof of Theorem~6.2 with the help of Lemma~6.3,
the still unproved Lemma~6.1 and still another result which
will be formulated in Lemma~6.4.

\medskip\noindent
{\it Proof of Theorem~6.2.} We want to prove that for all
positive integers $p$, real numbers $c_1,\dots,c_p$  and
$n_l\in{\mathbb Z}_\nu$, $l=1,\dots,p$,
$$
\sum_{l=1}^p c_l Z^N_{n_l}\stackrel{{\cal D}}{\rightarrow}
\sum_{l=1}^p c_lZ^*_{n_l},
$$
since this relation also implies the convergence of the
multi-dimensional distributions. Applying the same calculation
as in the heuristic justification of Theorem~6.2 we get that
$$
\sum_{l=1}^p c_lZ_{n_l}^N=\frac1{A_N}\sum_{l=1}^p c_l\int
\sum_{j\in B_{n_l}^N}e^{i(j,x_1+\cdots+x_k)}
\,Z_G(\,dx_1)\dots Z_G(\,dx_k),
$$
and by applying the scaling $y_j=Nx_j$ we can show that
$$
\sum_{l=1}^p c_lZ^N_{n_l}\stackrel{\Delta}{=}
\int K_N(x_1,\dots,x_k)
\,Z_{G_N}(\,dx_1)\dots Z_{G_N}(\,dx_k)
$$
with
\begin{eqnarray}
K_N(x_1,\dots,x_k)&=& \frac1{N^\nu}\sum_{l=1}^p c_l\sum_{j\in B_{n_l}^N}
\exp\left\{i\left(\frac jN,x_1+\cdots+x_k\right)\right\} \nonumber  \\
&=&f_N(x_1,\dots,x_k)\sum_{l=1}^p c_l\tilde\chi_{n_l}(x_1+\cdots+x_k).
\label{(8.14)}
\end{eqnarray}
with the function~$f_N$ defined in~(\ref{(8.5)}) and 
the measure~$G_N$ defined in~(\ref{(8.2)}), The function  
$\tilde\chi_n(\cdot)$ denotes again the Fourier 
transform of the indicator function of the unit cube 
$\prod\limits_{j=1}^\nu[n^{(j)},n^{(j)}+1)$, 
$n=(n^{(1)},\dots n^{(\nu)})$.

Let us define the function
$$
K_0(x_1,\dots,x_k)=\sum_{l=1}^p c_l\tilde\chi_{n_l}(x_1+\cdots+x_k)
$$
and the measures~$\mu_N$ on $R^{k\nu}$ by the formula
\begin{eqnarray}
&&\mu_N(A)=\int_A|K_N(x_1,\dots,x_k)|^2 G_N(\,dx_1)\dots G_N(\,dx_k),
\nonumber  \\
&&\qquad \quad A\in{\cal B}^{k\nu} \textrm{ and } N=0,1,\dots.
\label{(8.15)}
\end{eqnarray}
The measure $G_0$ defined with parameter $N=0$ is the vague 
limit of the measures~$G_N$.

We prove Theorem~6.2 by showing that Lemma~6.3 can be 
applied with these spectral measures~$G_N$ and
functions~$K_N$. Since 
$G_N\stackrel{v}{\rightarrow} G_0$, and $K_N\to K_0$ 
uniformly in all bounded regions in $R^{k\nu}$, it is 
enough  to show, beside the proof of Lemma~6.1, that 
the measures $\mu_N$, $N=1,2,\dots$, tend weakly to the 
(necessarily finite) measure~$\mu_0$ which is also defined 
in~(\ref{(8.15)}), (in notation 
$\mu_N\stackrel{w}{\rightarrow}\mu_0$), i.e.
$\int f(x)\mu_N(\,dx)\to\int f(x)\mu_0(\,dx)$ for all 
continuous and bounded functions~$f$ on~$R^{k\nu}$. 
Then this convergence implies condition~(b) in Lemma~6.3. 
Moreover, it is enough to show the slightly weaker statement 
by which there exists some finite measure $\bar\mu_0$ such 
that $\mu_N\stackrel{w}{\rightarrow}\bar\mu_0$, since then 
$\bar\mu_0$ must coincide with $\mu_0$ because of the 
relations $G_N\stackrel{v}{\rightarrow} G_0$ and 
$K_N\to K_0$ uniformly in all bounded regions of $R^{k\nu}$, 
and $K_0$ is a continuous function. This implies that
$\mu_N\stackrel{v}{\rightarrow}\mu_0$, and $\mu_0=\bar\mu_0$.

There is a well-known theorem in probability theory about 
the equivalence between weak convergence of finite measures 
and the convergence of their Fourier transforms. It would be 
natural to apply this theorem for proving
$\mu_N\stackrel{w}{\rightarrow}\bar\mu_0$. But actually we
shall need a version of this result. In this version we
exploit that we have the additional information that the 
measures $\mu_N$, $N=1,2,\dots$, are concentrated in the 
cubes $[-N\pi,N\pi)^{k\nu}$, since the spectral measure~$G$ 
is concentrated in $[-\pi,\pi)^\nu$. On the other hand, 
formula~(\ref{(8.1)}) provides only a restricted information
about the Fourier transform of~$\mu_N$. We have an asymptotic
relation on the function $r(n)$, i.e. on the Fourier transform 
of~$G$ only in the points $n\in{\mathbb Z}_\nu$. As we shall 
see, this implies that we have control on the Fourier 
transform~$\mu_N$ only in the points $\frac nN$, 
$n\in{\mathbb Z}_\nu$. 

Hence it will be more appropriate for us to work with such a version
of the result about the  equivalence of weak convergence of
probability measures which takes into account that we have only
restricted information about the Fourier transform, but on the
other hand we have the additional information that the probability
measures we are working with are concentrated in some 
well-defined cubes. We will formulate such a result in the 
next Lemma~6.4. Its proof will be postponed after the proof 
of Theorem~6.2. 

\medskip\noindent
{\bf Lemma 6.4.} {\it Let $\mu_1,\mu_2,\dots$ be a sequence 
of finite measures on $R^l$ such that 
$\mu_N(R^l\setminus [-C_N\pi,C_N\pi)^l)=0$
for all $N=1,2,\dots$, with some sequence $C_N\to\infty$ as
$N\to\infty$. Define the modified Fourier transform
$$
\varphi_N(t)=\int_{R^l} 
\exp\left\{i\left(\frac{[tC_N]}{C_N},x\right)\right\}
\mu_N(\,dx), \quad t\in R^l,
$$
where $[tC_N]$ is the integer part of the vector $tC_N\in R^l$. 
(For an $x\in R^l$ its integer part $[x]$ is the vector 
$n\in{\mathbb Z}_l$ for which $x^{(p)}-1<n^{(p)}\le x^{(p)}$ 
if $x^{(p)}\ge0$, and $x^{(p)}\le n^{(p)}< x^{(p)}+1$ if 
$x^{(p)}<0$ 
for all $p=1,2,\dots,l$.) If for all $t\in R^l$ 
the sequence $\varphi_N(t)$ tends to a function $\varphi(t)$ 
continuous at the origin, then the measures $\mu_N$ weakly 
tend to a finite measure~$\mu_0$, and $\varphi(t)$ is the 
Fourier transform of~$\mu_0$.}

\medskip
In the proof of Theorem~6.2 we apply Lemma~6.4 with 
$C_N=N$ and $l=k\nu$ for the measures $\mu_N$ defined 
in~(\ref{(8.15)}). Because of the middle term in~(\ref{(8.14)}) 
we can write the modified Fourier transform $\varphi_N$ of the
measure $\mu_N$ as
\begin{equation}
\varphi_N(t_1,\dots,t_k)=\sum_{r=1}^p\sum_{s=1}^p c_rc_s
\psi_N(t_1+n_r-n_s,\dots,t_k+n_r-n_s) \label{(8.19)}
\end{equation}
with
\begin{eqnarray}
&&\psi_N(t_1,\dots,t_r)=\frac1{N^{2\nu}}\int\exp\left\{
i\frac1N((j_1,x_1)+\cdots+(j_k,x_k))\right\}  \nonumber  \\
&&\qquad\qquad\sum_{u\in B^N_0}\sum_{v\in B^N_0}
\exp\left\{i\left(\frac{u-v}N,x_1+\cdots+x_k\right)\right\}
G_N(\,dx_1)\dots G_N(\,dx_k) \nonumber  \\
&&\qquad=\frac1{N^{2\nu-k\alpha}L(N)^k}\sum_{u\in B_0^N}\sum_{v\in B_0^N}
r(u-v+j_1)\cdots r(u-v+j_k), \label{(8.20)}
\end{eqnarray}
where  $j_p=[t_pN]$, $t_p\in R^\nu$, $p=1,\dots,k$. 

The asymptotic behaviour of $\psi_N(t_1,\dots,t_k)$ for 
$N\to\infty$ can be investigated by the help of the last 
relation and formula~(\ref{(8.1)}). Rewriting the last 
double sum in the form of a single sum  by fixing first 
the variable $l=u-v\in [-N,N]^\nu\cap{\mathbb Z}_\nu$, 
and then summing up for~$l$ one gets
$$
\psi_N(t_1,\dots,t_k)=\int_{[-1,1]^\nu} f_N(t_1,\dots,t_k,x)\,dx
$$
with
\begin{eqnarray*}
&&f_N(t_1,\dots,t_k,x) \\
&&\qquad =\left(1-\frac{[|x^{(1)}N|]}N\right)\cdots
\left(1-\frac{[|x^{(\nu)}N|]}N\right)
\frac{r([xN]+j_1)}{N^{-\alpha}L(N)}\cdots
\frac{r([xN]+j_k)}{N^{-\alpha}L(N)}.
\end{eqnarray*}
(In the above calculation we exploited that in the last sum of
formula~(\ref{(8.20)}) the number of pairs $(u,v)$ for which
$u-v=l=(l_1,\dots,l_\nu)$ equals $(N-|l_1|)\cdots(N-|l_\nu|)$.)

Let us fix some vector $(t_1,\dots,t_k)\in R^{k\nu}$. It can 
be seen with the help of formula~(\ref{(8.1)}) that 
for all $\varepsilon>0$ the convergence
\begin{equation}
f_N(t_1,\dots,t_k,x)\to f_0(t_1,\dots,t_k,x) \label{(8.21)}
\end{equation}
holds uniformly with the limit function
\begin{equation}
f_0(t_1,\dots,t_k,x)=(1-|x^{(1)}|)\dots(1-|x^{(\nu)}|)
\frac{a\left(\frac{x+t_1}{|x+t_1|}\right)}{|x+t_1|^\alpha}\dots
\frac{a\left(\frac{x+t_k}{|x+t_k|}\right)}{|x+t_k|^\alpha}
\label{(8.22)}
\end{equation} 
on the set $x\in[-1,1]^\nu\setminus
\bigcup\limits_{p=1}^k\{x\colon\;|x+t_p|>\varepsilon\}$. 

I claim that
$$
\psi_N(t_1,\dots,t_k)\to\psi_0(t_1,\dots,t_k)
=\int_{[-1,1]^\nu}f_0(t_1,\dots,t_k,x)\,dx,
$$
and $\psi_0$ is a continuous function.

This relation implies that $\mu_N\stackrel{w}{\rightarrow}\mu_0$. To
prove it, it is enough to show beside formula~(\ref{(8.21)}) that
\begin{equation}
\left|\int_{|x+t_p|<\varepsilon} 
f_0(t_1,\dots,t_k,x)\,dx\right|<C(\varepsilon),
\quad p=1,\dots,k, \label{(8.23)}
\end{equation}
and
\begin{equation}
\int_{|x+t_p|<\varepsilon} |f_N(t_1,\dots,t_k,x)|\,dx<C(\varepsilon),
\quad p=1,\dots,k, \quad \textrm{and } N=1,2,\dots \label{($8.24'$)}
\end{equation}
with a constant $C(\varepsilon)$ such that 
$C(\varepsilon)\to0$ as $\varepsilon\to0$.

By formula~(\ref{(8.22)}) and H\"older's inequality
\begin{eqnarray*}
\left|\int_{|x+t_p|<\varepsilon} f_0(t_1,\dots,t_k,x)\,dx\right|
&&\le C\prod_{1\le l\le k,\,l\neq p}\left[\int_{x\in[-1,1]^\nu}
|x+t_l|^{-k\alpha}\,dx\right]^{1/k} \\
&&\qquad \left[\int_{|x+t_p|\le\varepsilon}|x
+t_p|^{-k\alpha}\,dx\right]^{1/k}
\le C'\varepsilon^{\nu/k-\alpha}
\end{eqnarray*}
with some appropriate $C>0$ and $C'>0$, since $\nu-k\alpha>0$,
and $a(\cdot)$ is a bounded function. Similarly,
\begin{eqnarray}
\int_{|x+t_p|<\varepsilon} |f_N(t_1,\dots,t_k,x)|\,dx
&&\le \prod_{1\le l\le k,\,l\neq p}\left[\int_{x\in[-1,1]^\nu}
\frac{|r([xN]+j_l)|^k}{N^{-k\alpha}L(N)^k}\,dx\right]^{1/k},
\nonumber \\
&&\qquad \left[\int_{|x+t_p|\le\varepsilon}
\frac{|r([xN]+j_p)|^k}{N^{-k\alpha}L(N)^k}\,dx\right]^{1/k}.
\label{(8.25)}
\end{eqnarray}

It is not difficult to see with the help of Karamata's theorem 
that if $L(t)$, $t\ge1$, is a slowly varying function which is 
bounded in all finite intervals, then for all numbers $\eta>0$ 
and $K>0$ there are some constants $K_1=K_1(\eta,K)>0$, and 
$C=C(\eta,K)>0$ together with a threshold index~$N_0=N_0(\eta,K)$ 
such that
$$
\frac{L(uN)}{L(N)}\le Cu^{-\eta} \quad\textrm{if } uN>K_1, \;\; u\le K,
\textrm{ \ and } N\ge N_0.
$$
Hence formula~(\ref{(8.1)}) implies that
\begin{eqnarray}
|r([xN]+[t_lN])&=&|r([xN]+j_l)|
\le CN^{-\alpha}L(N)|x+t_l|^{-\alpha-\eta} \nonumber \\ 
&&\qquad \textrm{if } |x+t_l|\le K \textrm{ and } N\ge N_0. 
\label{(8.26)}
\end{eqnarray}
Relation~(\ref{(8.26)}) follows from the previous relation 
and~(\ref{(8.1)}) if $|[xN]+[t_lN]|\ge K_1$. It also holds if
$|[xN]+[t_lN]|\le K_1$, since in this case the left-hand side
can be bounded by the inequality $|r([xN]+[t_lN]|\le1$, while 
the right-hand side of~(\ref{(8.26)}) is greater than 1 with 
the choice of a sufficiently large constant~$C$ (depending on 
$\eta$ and $K_1$). This follows from the relation
 $|x+t|^{-\alpha-\eta}
=N^{\alpha+\eta}|N(x+t)|^{-\alpha-\eta}\ge C_1N^{\alpha+\eta}$
if $|[xN]+[t_lN]|\le K_1$, and $L(N)\ge N^{-\eta}$.

We get with the help of~(\ref{(8.26)}) that
\begin{eqnarray*}
&&\int_{|x+t_p|<\varepsilon}\frac{|r([xN]+j_p)|^k}{N^{-k\alpha}L(N)^k}\,dx
\le B\int_{|x+t_p|<\varepsilon}|x+t_p|^{-k(\alpha+\eta)}\,dx
\le B'\varepsilon^{\nu-k(\alpha +\eta)}\\
&&\int_{x\in[-1,1]^\nu}\frac{|r([xN]+j_l)|^k}{N^{-k\alpha}L(N)^k}\,dx
\le B''.
\end{eqnarray*}
for a sufficiently small constant $\eta>0$ with some constants
$B,B',B''<\infty$ depending on $\eta$ and $t_p$, $1\le p\le k$.

Therefore we get from~(\ref{(8.25)}), by choosing an $\eta>0$ 
such that $k(\alpha+\eta)<\nu$, that the inequality
$$
\int_{|x+t_p|<\varepsilon} |f_N(t_1,\dots,t_k,x)|\,dx
\le C\varepsilon^{\nu/k-(\alpha +\eta)}
$$
holds with some $C<\infty$. The right-hand side of this 
inequality tends to zero as $\varepsilon\to0$. Hence 
we proved beside~(\ref{(8.21)}) formulae~(\ref{(8.23)})
and~(\ref{($8.24'$)}), and they have the consequence that
$\psi_N(t_1,\dots,t_k)\to\psi_0(t_1,\dots,t_k)$.
Since $\psi_0(t_1,\dots,t_k)$ is a continuous function
relation~(\ref{(8.19)}) with Lemma~6.4 imply that the
measures $\mu_N$ introduced in~(\ref{(8.18)}) converge
weakly to a probability measure as $N\to\infty$, and as we
saw at the beginning of the proof of Theorem~6.2 this
limit measure must be~$\mu_0$.

Hence we can apply Lemma~6.3 for the spectral measures~$G_N$ 
and functions $K_N(\cdot)$, $N=0,1,2,\dots$, defined in 
Theorem~6.2. The convergence 
$G_N\stackrel{v}{\rightarrow} G_0$ follows from Lemma~6.1. 
Conditions~(a) and~(b) also hold with the choice of a 
sufficiently large number~$A=A(\varepsilon)$. The hard 
point of the proof was to check condition~(b). 
This followed from the relation 
$\mu_N\stackrel{w}{\rightarrow} \mu_0$. Thus we have proved 
Theorem~6.2 with the help of Lemmas~6.1, 6.3 and~6.4. 

\medskip
Now I turn to the proof of Lemma~6.4. Before writing it down I 
make some comments on its conditions. Let us observe
that if the measures~$\mu_N$ or a part of them are shifted with a
vector $2\pi C_N u$ with some $u\in{\mathbb Z}_l$, then their
modified Fourier transforms $\varphi_N(t)$ do not change because of
the periodicity of the trigonometrical functions $e^{i(j/C_N,x)}$,
$j\in{\mathbb Z}_l$. On the other hand, these new measures
which are not concentrated in $[-C_N\pi,C_N\pi)^l$, have no limit.
Lemma~6.4 states that if the measures~$\mu_N$ are concentrated in the
cubes $[-C_N\pi,C_N\pi)^l$, then the convergence of their modified
Fourier transforms defined in Lemma~6.4, which is a weaker condition,
than the convergence of their Fourier transforms, also implies their
convergence to a limit measure.

\medskip\noindent
{\it Proof of Lemma~6.4.}\/ The proof is a natural modification of
the proof about the equivalence of weak convergence of measures and
the convergence of their Fourier transforms. First we show that
for all $\varepsilon>0$ there exits some $K=K(\varepsilon)$ such that
\begin{equation}
\mu_N(x\colon\; x\in R^l,\; |x^{(1)}|>K)<\varepsilon 
\quad \textrm{for all \ }N\ge1.
\label{(8.16)}
\end{equation}

As $\varphi(t)$ is continuous in the origin there is some $\delta>0$
such that
\begin{equation}
|\varphi(0,\dots,0)-\varphi(t,0,\dots,0)|<\frac\varepsilon2 
\quad\textrm{if \ } |t|<\delta. \label{(8.17)}
\end{equation}
We have
\begin{equation}
0\le \textrm{Re}\,[\varphi_N(0,\dots,0)-\varphi_N(t,0,\dots,0)]
\le2\varphi_N(0,\dots,0) \label{(8.18)}
\end{equation}
for all $N=1,2,\dots$. The sequence in the middle term 
of~(\ref{(8.18)}) tends to 
$$
\textrm{Re}\,[\varphi(0,\dots,0)-\varphi(t,0,\dots,0)]
$$ 
as $N\to\infty$. The right-hand side of~(\ref{(8.18)}) is a 
bounded sequence, since it is convergent. Hence the dominated 
convergence theorem can be applied for the functions
$\textrm{Re}\,[\varphi_N(0,\dots,0)-\varphi_N(t,0,\dots,0)]$.
We get because of the condition $C_N\to\infty$ and
relation~(\ref{(8.17)}) that
\begin{eqnarray*}
&&\lim_{N\to\infty} \int_0^{[\delta C_N]/C_N} \frac1\delta\,
\textrm{Re}\,[\varphi_N(0,\dots,0)-\varphi_N(t,0,\dots,0)]\,dt\\
&&\qquad=\int_0^\delta\frac1\delta\,
\textrm{Re}\,[\varphi(0,\dots,0)-\varphi(t,0,\dots,0)]\,dt
<\frac\varepsilon2
\end{eqnarray*}
with the number $\delta>0$ appearing in~(\ref{(8.17)}). Hence
\begin{eqnarray*}
\frac\varepsilon2&>& \lim_{N\to\infty} 
\int_0^{[\delta C_N]/C_N}\frac1\delta\,
\textrm{Re}\,
[\varphi_N(0,\dots,0)-\varphi_N(t,0,\dots,0)]\,dt \\
&=&\lim_{N\to\infty}\int
\left(\frac1\delta\int_0^{[\delta C_N]/C_N}
\textrm{Re}\, [1-e^{i[tC_N]x^{(1)}/C_N}]\,dt\right) 
\mu_N(\,dx)\\
&=&\lim_{N\to\infty}\int
\frac1{\delta C_N}\sum_{j=0}^{[\delta C_N]-1}
\textrm{Re}\,\left[1-e^{ijx^{(1)}/C_N}\right]\mu_N(\,dx)\\
&\ge&\limsup_{N\to\infty} \int_{\{|x^{(1)}|>K\}}
\frac1{\delta C_N} \sum_{j=0}^{[\delta C_N]-1} 
\textrm{Re}\, \left[1-e^{ijx^{(1)}/C_N}\right]\mu_N(\,dx)\\
&=&\limsup_{N\to\infty}\int_{\{|x^{(1)}|>K\}}
\left(1-\frac1{\delta C_N}
\textrm{Re}\,\frac{1-e^{i[\delta C_N]x^{(1)}/C_N}}
{1-e^{ix^{(1)}/C_N}}\right)\mu_N(\,dx)
\end{eqnarray*}
with an arbitrary $K>0$. (In the last but one step of this 
calculation we have exploited that 
$\frac1{\delta C_N}\sum\limits_{j=0}^{[\delta C_N]-1}
\textrm{Re}\,[1-e^{ijx^{(1)}/C_N}]\ge0$ for all 
$x^{(1)}\in R^1$.)

Since the measure $\mu_N$ is concentrated in
$\{x\colon\, x\in R^l,\;|x^{(1)}|\le C_N\pi\}$, and
\begin{eqnarray*}
\textrm{Re}\,\frac{1-e^{i[\delta C_N]x^{(1)}/C_N}}
{1-e^{ix^{(1)}/C_N}}
&=&\frac{\textrm{Re}\,\left(i e^{-ix^{(1)}/2C_N}
\left(1-e^{i[\delta C_N]x^{(1)}/C_N}\right)\right)}
{i(e^{-ix^{(1)}/2CN}-e^{ix^{(1)}/2C_N})}\\
&\le&\frac1{\left|\sin \left(\dfrac{x^{(1)}}{2C_N}\right)\right|}
\le \frac{C_N\pi}{|x^{(1)}|}
\end{eqnarray*}
if $|x^{(1)}|\le C_N\pi$, (here we exploit that
$|\sin u|\ge\frac2\pi|u|$ if $|u|\le\frac\pi2$), hence we have with
the choice $K=\frac{2\pi}{\delta}$
$$
\frac\varepsilon2>\limsup_{N\to\infty}\int_{\{|x^{(1)}|>K\}}
\left(1-\left|\frac\pi{\delta x^{(1)}}\right|\right)\mu_N(\,dx)
\ge\limsup_{N\to\infty}\frac12\mu_N(|x^{(1)}|>K).
$$
As the measures $\mu_N$ are finite the inequality
$\mu_N(|x^{(1)}|>K)<\varepsilon$ holds for each index~$N$ with a
constant~$K=K(N)$ that may depend on~$N$. Hence the above
inequality implies that formula~(\ref{(8.16)}) holds for all $N\ge1$
with a possibly larger index~$K$ that does not depend on~$N$.

Applying the same argument to the other coordinates we 
find that for all $\varepsilon>0$ there exists some 
$C(\varepsilon)<\infty$ such that
$$
\mu_N\left(R^l\setminus[-C(\varepsilon),
C(\varepsilon)]^l\right)<\varepsilon \quad 
\textrm{for all } N=1,2,\dots.
$$

Consider the usual Fourier transforms
$$
\tilde\varphi_N(t)=\int_{R^l}e^{i(t,x)}\mu_N(\,dx), \quad t\in R^l.
$$
Then
\begin{eqnarray*}
|\varphi_N(t)-\tilde\varphi_N(t)|&\le& 2\varepsilon
+\int_{[-C(\varepsilon),C(\varepsilon)]}
\left|e^{i(t,x)}-e^{i([tC_N]/C_N,x)}\right|\mu_N(\,dx)\\
&\le& 2\varepsilon
+\frac{lC(\varepsilon)}{C_N}\mu_N(R^l)
\end{eqnarray*}
for all $\varepsilon>0$. Hence $\tilde\varphi_N(t)-\varphi_N(t)\to0$ as
$N\to\infty$, and $\tilde\varphi_N(t)\to\varphi(t)$. (Observe that
$\mu_N(R^l)=\varphi_N(0)\to\varphi(0)<\infty$ as $N\to\infty$, hence
the measures~$\mu_N(R^l)$ are uniformly bounded, and $C_N\to\infty$
by the conditions of Lemma~6.4.) Then Lemma~6.4 follows from standard
theorems on Fourier transforms.

\medskip
It remained to prove Lemma~6.1.

\medskip\noindent
{\it Proof of Lemma~6.1.}\/ Introduce the notation
$$
K_N(x)=\prod_{j=1}^\nu\frac {e^{ix^{(j)}}-1}{N(e^{ix^{(j)}/N}-1)},
\quad N=1,2,\dots,
$$
and
$$
K_0(x)=\prod_{j=1}^\nu\frac {e^{ix^{(j)}}-1}{ix^{(j)}}.
$$
Let us consider the measures $\mu_N$ defined in 
formula~(\ref{(8.15)}) in the special case $k=1$ with $p=1$, 
$c_1=1$ in the definition of the function $K_N(\cdot)$, i.e. put 
$$
\mu_N(A)=\int_A |K_N(x)|^2\,G_N(\,dx),
\quad A\in {\cal B}^\nu, \quad N=1,2,\dots.
$$
We have already seen in the proof of Theorem~6.2 that
$\mu_N\stackrel{w}{\rightarrow} \mu_0$ with some finite 
measure~$\mu_0$, and the Fourier transform of~$\mu_0$ is
$$
\varphi_0(t)=\int_{[-1,1]^\nu}(1-|x^{(1)}|)\cdots (1-|x^{(\nu)}|)
\frac{a\left(\frac{x+t}{|x+t|}\right)}{|x+t|^\alpha}\,dx.
$$
Moreover, since $|K_N(x)|^2\to |K_0(x)|^2$ uniformly in 
any bounded domain, it is natural to expect that 
$G_N\stackrel{v}{\rightarrow}G_0$ with 
$G_0(\,dx)=\frac1{|K_0(x)|^2}\mu_0(\,dx)$. But $K_0(x)=0$
in some points (if $x^{(j)}=2k\pi$ with some integer 
$k\neq0$ for a coordinate of $x$), and                                       %
the function $K_0(\cdot)^{-2}$ is not continuous here. As 
a consequence, we cannot give a direct proof of the above 
statement. Hence we apply instead a modified version of 
this method. First we prove the following result about 
the behaviour of the restrictions of the measures $G_N$ 
to appropriate cubes: 

For all $T\ge1$ there is a finite measure 
$G_0^T$ concentrated on $(-T\pi,T\pi)^\nu$ such that
\begin{equation}
\lim_{N\to\infty}\int f(x)\,G_N(\,dx)=\int f(x)\,G_0^T(\,dx) 
\label{(8.27)}
\end{equation}
for all continuous functions~$f$ which vanish outside the cube
$(-T\pi,T\pi)^\nu$.

Indeed, let a continuous function $f$ vanish outside the cube 
$(-T\pi, T\pi)^\nu$ with some $T\ge1$. Put $M=[\frac N{2T}]$. 
Then
\begin{eqnarray*}
\int f(x)G_N(\,dx)&=&\frac{N^\alpha}{L(N)}
\cdot \frac{L(M)}{M^\alpha}
\int f\left(\frac NMx\right)G_M(\,dx)\\
&=&\frac{N^\alpha L(M)}{M^\alpha L(N)}\int f
\left(\frac NMx\right) |K_M(x)|^{-2}\mu_M(\,dx)\\
&&\qquad \to (2T)^\alpha
\int f(2Tx)|K_0(x)|^{-2}\mu_0(\,dx)\\
&&\qquad\qquad
=\int f(x)\frac{(2T)^\alpha}{|K_0(\frac x{2T})|^2}\mu_0
\left(\,\frac{dx}{2T}\right)
\quad \textrm{as }N\to\infty,
\end{eqnarray*}
because $f(\frac NMx)|K_M(x)|^{-2}$ vanishes outside the cube
$[-\pi,\pi]^\nu$, the limit relation
$$
f(\frac NMx)|K_M(x)|^{-2}\to f(2Tx)|K_0(x)|^{-2} 
$$ 
holds uniformly, (the function~$K_0(\cdot)^{-2}$ is 
continuous in the cube $[-\pi,\pi]^\nu$), and 
$\mu_M\stackrel{w}{\rightarrow}\mu_0$ as $N\to\infty$. 
Hence relation~(\ref{(8.27)}) holds if we define $G_0^T$ 
as the restriction of the measure
$\frac{(2T)^\alpha}{|K_0(\frac x{2T})|^2}\mu_0
\left(\,\frac{dx}{2T}\right)$ to the cube $(-T\pi,T\pi)^\nu$.
The measures $G_0^T$ appearing in~(\ref{(8.27)}) are 
consistent for different parameters~$T$, i.e.\ $G_0^T$ 
is the restriction of the measure $G_0^{T'}$ to the cube 
$(-T\pi,T\pi)^\nu$ if $T'>T$. This follows from the fact
that $\int f(x)G_0^T(\,dx)=\int f(x)G_0^{T'}(\,dx)$ for
all continuous functions with support in $(-T,T)^\nu$.
We claim that by defining the measure $G_0$ by the relation
$G_0(A)=G_0^T(A)$ for a bounded set $A$ and such number
$T>1$ for which $A\subset (-T\pi,T\pi)^\nu$ we get such a
locally finite measure~$G_0$ for which 
$G_N\stackrel{v}{\rightarrow}G_0$. The above mentioned
vague convergence is a direct consequence of~(\ref{(8.27)})
and the definition of~$G_0$, but to give a complete proof
we have to show that $G_0$ is really a ($\sigma$-additive) 
measure.

Actually it is enough to prove that the restriction of
$G_0$ to the bounded, measurable sets is $\sigma$-additive, 
because it follows then from standard results in measure 
theory that it has a unique $\sigma$-additive extension to 
${\cal B}^\nu$. But this is an almost direct consequence
of the definition of $G_0$. The desired $\sigma$-additivity 
clearly holds, since if $A=\bigcup\limits_{n=1}^\infty A_n$, 
the set $A$ is bounded, and the sets $A_n$, $n=1,2,\dots$, 
are disjoint, then there is a number $T>1$ such that 
$A\subset (-T\pi,T\pi)^\nu$, the same relation holds for
the sets $A_n$, and the $\sigma$-additivity of $G_0^T$ 
implies that $G_0(A)=\sum\limits_{n=1}^\infty G_0(A_n)$.

As $G_N\stackrel{v}{\rightarrow} G_0$, and
$|K_N(x)|^2\to|K_0(x)|^2$ uniformly in all bounded regions, 
the relation $\mu_N\stackrel{v}{\rightarrow}\bar\mu_0$ 
holds with the measure $\bar\mu_0$ defined by the formula
$\bar\mu_0(A)=\int_A|K_0(x)|^2G_0(\,dx)$, $A\in{\cal B}^\nu$. 
Since
$\mu_N\stackrel{w}{\rightarrow}\mu_0$ the measures~$\mu_0$
and~$\bar\mu_0$ must coincide, i.e.
$$
\mu_0(A)=\int_A |K_0(x)|^2\,G_0(\,dx), \quad A\in{\cal B}^\nu.
$$
Relation~(\ref{(8.4)}) expresses the fact that $\varphi_0$ is 
the Fourier transform of~$\mu_0$.
It remained to prove the homogeneity property~(\ref{(8.3)})
of the measure~$G_0$. For this goal let us extend the 
definition of the measures~$G_N$ given in~(\ref{(8.2)}) to
all non-negative real numbers~$u$. It is easy to see that 
the relation $G_u\stackrel{v}{\rightarrow} G_0$ as 
$u\to\infty$ remains valid. Hence we get for all fixed 
$s>0$ and continuous functions~$f$ with compact support that
\begin{eqnarray*}
\int f(x) G_0(\,dx)&=&\lim_{u\to\infty}\int f(x)\,G_u(\,dx)
=\lim_{u\to\infty}\frac{s^\alpha L(\frac us)}{L(u)}
\int f(sx)G_{\frac us}(\,dx)\\
&=&s^\alpha \int f(sx)G_0(\,dx)=\int f(x)s^\alpha 
G_0 \left(\frac{dx}s\right).
\end{eqnarray*}
This identity implies the homogeneity property~(\ref{(8.3)}) 
of~$G_0$. Lemma~6.1 is proved. 

\subsection{A discussion about our results}

Lemma 6.1 is a result about the limit behaviour 
of the spectral distribution of a stationary
random fields if its correlation function satisfies 
formula~(\ref{(8.1)}). It states that under this
condition the appropriately rescaled spectral measure 
has a limit in the vague convergence sense, and
Lemma~6.1 also describes this limit. There is a closer
relation between the behaviour of the correlation
function and spectral measure which  may be worthwhile 
for a more detailed discussion. Moreover, this 
comparison may help us to understand the relation 
between limit theorems about the large scale 
limit of stationary Gaussian random fields and 
such limit theorems about non-linear functionals
of stationary Gaussian fields which are similar to 
Theorem~6.2 of this work.

In the subsequent slightly informal discussion I disregard
the appearance of the slowly varying function $L(t)$ in our
results, I assume simply that $L(t)\sim 1$ as $t\to\infty$.
In this case we can interpret Lemma~6.1 so that if the
correlation function $r(n)$ behaves like $L(n)\sim|n|^{-\alpha}$
in the neighbourhood of the infinity, then the spectral
measure behaves like $G(t)\sim\textrm{const.}\, t^\alpha$ 
as $t\to0$. (Here $G(t)$ denotes the spectral measure of
the ball in with radius $t$ and center point at the origin.) 
One may ask, what can be said in the opposite direction. What 
can we say about the asymptotic behaviour of
the correlation function if we have some information about
the behaviour of the spectral measure? For the
sake of simplicity let us restrict our attention to the
the correlation function of one-dimensional random fields,
i.e. to the case when $\nu=1$. In the following consideration 
I shall apply some heuristic not completely precise argument.

By some results about Fourier analysis we can say that the
smoother is a function the faster tends its Fourier transform
to zero at infinity. On the other hand, if a function has
a singularity, but otherwise it is smooth enough, then the
asymptotic behaviour of its Fourier transform at infinity 
is determined by this singularity. This means in particular 
that if the spectral measure behaves like $G(t)\sim C\cdot t^{\alpha}$, 
$0<\alpha<1$, (or it has a spectral density has the form 
$g(t)\sim|t|^{\alpha-1}$) in the neighbourhood of the origin,
and this is the strongest irregularity of the spectral measure,
then the correlation function of the random field satisfies
condition~(\ref{(8.1)}) of Theorem~6.2.

In the paper P. Major: On renormalizing Gaussian fields. 
{\it Z.~Wahrscheinlichkeitstheorie verw. Gebiete\/} {\bf 59} 
(1982), 515--533. the condition for a limit theorem for
Gaussian fields was was given by means of the spectral measure
while in the result of Theorem~6.2 the condition of a 
limit theorem was given by means of the correlation function
of the underlying Gaussian field. It may be worth comparing 
the conditions of these two results. 

The limit theorem about the existence of the large scale 
limit of a stationary Gaussian field can be interpreted in 
a slightly informal way so that the limit exists if the spectral 
measure has the singularity of the form $G(t)\sim C|t|^{-\alpha}$ 
in a small neighbourhood of the origin. On the other hand, the
condition of Theorem~6.2 was that $r(n)\sim C|n|^{-\alpha}$. As I 
mentioned before this is a stronger condition which implies that
the spectral measure behaves in a small neighbourhood of the
origin similarly to the previous case, but it also implies some
additional restriction. The spectral measure cannot have a
stronger singularity outside zero which would influence too
strongly the behaviour of its Fourier transform at the infinity.
I present an example taken from the fourth section of my paper
with R.~L.~Dobrushin that shows such a picture which
the above considerations suggest.

Take a stationary Gaussian sequence $X_n$, $EX_n=0$, $EX_n^=1$, 
$n=0,\pm1,\pm2,\dots$, with spectral density 
$$
g(x)=C_1|x|^{-\alpha}+C_2(|x-a|^{-\beta}+|x+a|^{-\beta}), \quad -\pi\le x<\pi,
$$
where $0<\alpha<\beta<1$, $\beta>\frac12$, $0<x<\pi$, $C_1,C_2>0$.
We are interested in what kind of limit theorems hold for the
sums $S_n=\frac1{A_n}\sum_{j=1}^nX_j$
and $T_n=\frac1{B_n}\sum_{j=1}^nH_2(X_j)=\frac1{B_n}\sum_{j=1}^n(X_j^2-1)$.
In particular, how do we have to choose the norming constants 
$A_n$ and $B_n$ to get a limit. (In the definition of $T_n$ 
we are working with the Hermite polynomial $H_2(x)=x^2-1$.) 
It can be proved that in the example with this spectral 
density the correlation function has the following form.
\begin{eqnarray*}
EX_kX_{p+k}&=&\int_{-\pi}^\pi e^{ipx}g(x)\,dx \\
&=&K_1p^{\alpha-1}\left(1+O\left(\frac1p\right)\right)
+K_2p^{\beta-1}\cos pa \left(1+O\left(\frac1p\right)\right)
\end{eqnarray*}
with some positive constants $K_1$ and $K_2$.

In the first problem, where we study the limit behaviour of $S_n$
we have a Gaussian limit, and some calculation shows that the
variance of the sum without the normalization is of order 
$n^{1+\alpha}$, which means that we get a limit with the 
norming constant $A_n=n^{(1+\alpha)/2}$. This means that in the 
limit behaviour the singularity $|t|^{-\alpha}$ of the spectral 
density at the origin is important.

In the case of the second limit problem the situation is different.
In this case we can calculate the correlation function of the terms
$H_2(X_j)$ e.g. with the help of the diagram formula, and some 
calculation yields that
$$
EH_2(X_k)H_2(X_{p+k})=K_2^2p^{2\beta-1}(1+\cos 2pa +o(1)).
$$
Further calculation shows that in this case the right norming
for $T_n$ for which the variance of $T_n$ is separated both from
zero and infinity is $B_n=n^\beta$. Some further calculation shows
that all moments of $T_n$ has a limit, moreover these limits
determine the limit distribution, hence there exits limit
theorem in this case. Finally the third moment of the limit
is positive, and this means that the limit is non-Gaussian.

This means that in the second problem the singularity 
$|x\pm a|^\beta$ gives the dominating factor that determines
the limit distribution. A more complete description of the
situation would demand further investigation.

\medskip
In the results of this section we discussed the limit
behaviour of the large scale limit of a random field
$H(X_n)$, $n\in{\mathbb Z}_\nu$, defined with the help
of a function $H(x)$ (square integrable with respect to 
the standard Gaussian measure) and a stationary Gaussian
random sequence $X_n$, $n\in{\mathbb Z}_\nu$, whose correlation
function satisfies relation (\ref{(8.1)}) with some parameter
$\alpha>0$. It was proved that if this parameter $\alpha$ is
not too large (this condition was formulated in a more explicit 
form), then we have a non-Gaussian limit theorem. To get a 
more complete picture one would like to know what can be said 
if this parameter $\alpha$ is relatively large, which means 
some sort of weak dependence. Next I formulate a result in 
this case, but because of lack of time I omit its proofs. 
First I formulate the above problem in more detail.

Let us consider a slightly more general version of the  
problem investigated in Theorem~$6.2'$. Take a stationary 
Gaussian random field  $X_n$, $EX_n=0$, $EX_n^2=1$, 
$n\in{\mathbb Z}_\nu$, with a correlation function 
satisfying relation~(\ref{(8.1)}), and the field 
$\xi_n=H(X_n)$, $n\in{\mathbb Z}_\nu$, subordinated to 
it with a general function $H(x)$ such that $EH(X_n)=0$ and
$EH(X_n)^2<\infty$. We are interested in the large-scale 
limit of such random fields. Take the Hermite 
expansion~(\ref{(8.28)}) of the function~$H(x)$, and let 
$k$ be the smallest such index for which $c_k\neq0$ in the 
expansion~(\ref{(8.28)}). In Theorem~$6.2'$ we solved
this problem if $0<k\alpha<\nu$. We are interested in 
the question what happens in the opposite case when 
$k\alpha>\nu$. Let me remark that in the case 
$k\alpha\ge\nu$ the field $Z^*_n$,
$n\in{\mathbb Z}_\nu$, which appeared in the limit in
Theorem~$6.2'$ does not exist. The Wiener-It\^o integral
defining $Z^*_n$ is meaningless, because the integral 
which should be finite to guarantee the existence of the 
Wiener--It\^o integral is divergent in this case. I 
formulate a general result which contains the answer to 
the above question as a special case.

\medskip\noindent
{\bf Theorem~6.5.} {\it Let us consider a stationary Gaussian 
random field $X_n$, $EX_n=0$, $EX_n^2=1$, 
$n\in{\mathbb Z}_n$, with correlation function 
$r(n)=EX_mX_{m+n}$, $m,n\in{\mathbb Z}_\nu$. Take a 
function $H(x)$ on the real line such that $EH(X_n)=0$ and
$EH(X_n)^2<\infty$. Take the Hermite expansion~(\ref{(8.28)}) 
of the function~$H(x)$, and let $k$ be smallest index in 
this expansion such that $c_k\neq 0$. If
\begin{equation}
\sum_{n\in{\mathbb Z}_\nu}|r(n)|^k<\infty, \label{(8.29)}
\end{equation}
then the limit
$$
\lim_{N\to\infty} EZ^N_n(H_l)^2=\lim_{N\to\infty}
N^{-\nu}\sum_{i\in B^N_n}\sum_{j\in B^N_n}r^l(i-j)=\sigma_l^2l!
$$
exists for all indices $l\ge k$, where $Z^N_n(H_l)$ is 
defined in~(\ref{(1.1)}) with $A_N=N^{\nu/2}$, and 
$\xi_n=H_l(X_n)$ with the $l$-th Hermite polynomial $H_l(x)$ 
with leading coefficient~1. Moreover, also the inequality
$$
\sigma^2=\sum_{l=k}^\infty c_l^2l!\sigma_l^2<\infty
$$
holds.

The finite dimensional distributions of the random field 
$Z^N_n(H)$ defined  in~(\ref{(1.1)}) with $A_N=N^{\nu/2}$ 
and $\xi_n=H(X_n)$ tend to the finite  dimensional 
distributions of a random field $\sigma Z^*_n$ with the  
number~$\sigma$ defined in the previous relation, where 
$Z^*_n$, $n\in{\mathbb Z}_\nu$, are independent,
standard normal random variables.}

\medskip
Theorem 6.5 can be applied if the conditions of Theorem~$6.2'$
hold with the only modification that the condition $k\alpha<\nu$
is replaced by the relation $k\alpha>\nu$. In this case the
relation~(\ref{(8.29)}) holds, and the large-scale limit of the 
random field $Z^N_n$, $n\in{\mathbb Z}_\nu$ with normalization
$A_N=N^{\nu/2}$ is a random field consisting of independent
standard normal random variables multiplied with the
number~$\sigma$. There is a slight generalization of
Theorem~6.5 which also covers the case $k\alpha=\nu$. In this
result we assume instead of the condition~(\ref{(8.29)}) that
$\sum\limits_{n\in \bar B_N} r(n)^k=L(N)$ with a slowly varying
function $L(\cdot)$, where
$\bar B_N=\{(n_1,\dots,n_\nu)\in{\mathbb Z}_\nu\colon\;
-N\le n_j\le N,\; 1\le j\le\nu\}$, and some additional
condition is imposed which states that an appropriately defined
finite number  $\sigma^2=\lim\limits_{N\to\infty}\sigma_N^2$, which
plays the role of the variance of the random variables in the
limiting field, exists. There is a similar large scale limit in
this case as in Theorem~6.5, the only difference is that the
norming constant in this case is $A_N=N^{\nu/2}L(N)^{1/2}$. This
result has the consequence that if the conditions of
Theorem~$6.2'$ hold with the only difference that $k\alpha=\nu$
instead of $k\alpha<\nu$,  then the large scale limit exists
with norming constants $A_N=N^{\nu/2}L(N)$ with an  appropriate
slowly varying function~$L(\cdot)$, and it consists of
independent Gaussian random variables with expectation zero.

The proof of Theorem~6.5 and its generalization that I did 
not formulate here explicitly appeared in in my paper with 
P.~Breuer {\it Central limit theorems for non-linear 
functionals of Gaussian fields} Journal of Multivariate 
Analysis {\bf 13} (1983), 425--441. I omit its proof, I 
only make some short explanation about it.

In the proof we show that all moments of the random variables
$Z^N_n$ converge to the corresponding moments of the Gaussian 
random variables $Z^*_n$ with expectation zero and the right
variance as $N\to\infty$. The moments of the random
variables $Z_n^N$ can be calculated by means of the diagram
formula if we either rewrite them in the form of a Wiener--It\^o
integral or apply a version of it for the moments of Hermite
(or of their generalization, the Wick polynomials) instead of 
Wiener--It\^o integrals. In both cases the moments can be 
expressed explicitly by means of the correlation function of 
the underlying Gaussian random field. The most important step 
of the proof is to show that we can select a special subclass 
of (closed) diagrams, called regular diagrams in my paper 
with P.~Breuer which yield the main contribution to the 
moments $E(Z^N_n)^M$, and their contribution can be simply 
calculated. The contribution of all remaining diagrams 
is~$o(1)$ (after norming), hence it is negligible. For 
the sake of simplicity let us restrict our attention to 
the case $H(x)=H_k(x)$, and let us explain the definition 
of the regular diagrams in this special case. 

If the number of the rows $M$ is an even number, then we 
call a closed diagram regular if there is a pairing of 
the rows, i.e. a partition $\{k_1,k_2\}$,  $\{k_3,k_4\}$,\dots, 
$\{k_{M-1},k_M\}$ of the set $\{1,\dots,M\}$ into subsets
of two elements in such a way that an edge can connect
only vertices in paired rows. If $M$ is an odd number, 
then there is no regular diagram. The main step of the
proof is to show that the contribution of all remaining
closed diagrams is negligibly small.

\Section{Appendix: The proof of Lemma 3.3}

In the Appendix I present a new proof of Lemma~3.3 
which is simpler than its original version that appeared 
as the proof of Lemma~4.1 in my Lecture Note 
{\it Multiple Wiener--It\^o integrals}.           

Our goal is to find for all functions $f\in\bar{{\cal H}}_G^n$
and $\varepsilon>0$ a function $f'\in\hat{\bar{{\cal H}}}_G^n$ 
such that the distance of $f$ and $f'$ is smaller than $\varepsilon$
in the Hilbert space  $\bar{{\cal H}}_G^n$. Then the corresponding
statement about functions in the Hilbert space  ${\cal H}_G^n$
follows from a standard symmetrization procedure.

Let us first observe that if two functions $f_1\in\bar{{\cal H}}_G^n$ 
and $f_2\in\bar{{\cal H}}_G^n$ can be arbitrarily well 
approximated by functions from $\hat{\bar{{\cal H}}}_G^n$ 
in the norm of this space, then the same relation 
holds for any linear combination $c_1f_1+c_2f_2$ with real 
coefficients~$c_1$ and~$c_2$. Indeed, if  the functions $f_i$, 
$i=1,2$, are approximated by some functions 
$g_i\in\hat{\bar{{\cal H}}}_G^n$, then we may assume, by 
applying some refinement of the partitions if it is necessary, 
that the approximating functions $g_1$ and $g_2$ are adapted 
to the same regular partition. Hence also 
$c_1g_1+c_2g_2\in\hat{\bar{{\cal H}}}_G^n$, and it
provides a good approximation of $c_1f_1+c_2f_2$.

The above observation enables us to reduce the proof of
Lemma~3.3 to the proof of a simpler statement formulated
in the following  {\it Statement A}. Here we have to 
approximate simpler functions $f\in\bar{{\cal H}}_G^n$. 
We have to consider two different cases. In the first 
case the function $f$ is the indicator function of some 
set $A\in R^{n\nu}$. In the second case $f$ is a simple function 
taking imaginary values. It takes the value $i=\sqrt{-1}$ 
in a set~$A$, the value~$-i$ in the set $-A$, and otherwise 
it equals zero. Here is the formulation of {\it Statement~A}. 

\medskip\noindent
{\it Statement A.}\/ Let $A\in{\cal B}^{n\nu}$ be a bounded, 
symmetric set, i.e. let $A=-A$. Then for any $\varepsilon>0$ 
there is a function $g\in\hat{\bar{{\cal H}}}_G^n$ such that 
$g=\chi_B$ with some set $B\in{\cal B}^{n\nu}$, i.e. $g$~is 
the indicator function of a set~$B$ such that the inequality 
$\|g-\chi_A\|<\varepsilon$ holds with the norm of the 
space $\bar{{\cal H}}_G^n$. (Here $\chi_A$ denotes the 
indicator function of the set~$A$, and we have
$\chi_A\in\bar{{\cal H}}_G^n$.) 

If $\chi_A\in\bar{{\cal H}}_G^n$ is a bounded set, and there 
is such a set $A_1$ for which the set $A$  can be written 
in the form $A=A_1\cup (-A_1)$, and the sets $A_1$ and $-A_1$ 
have a positive distance from each other, i.e. 
$\rho(A_1,-A_1)=\inf\limits_{x\in A_1,\,y\in -A_1}\rho(x,y)>\delta$,
with some $\delta>0$, where $\rho$ denotes the Euclidean 
distance in $R^{n\nu}$, then a good approximation of 
$\chi_A$ can be given with such a function 
$g=\chi_{B\cup(-B)}\in\hat{\bar{{\cal H}}}_G^n$ for which
the sets $B$ and $-B$ are disjoint, and the set $B$ is 
close to $A_1$. More explicitly, for all $\varepsilon>0$ 
there is a set $B\in{\cal B}^{n\nu}$ such that
$B\subset A_1^{\delta/2}=\{x\colon\;\rho(x,A_1)\le\frac\delta2\}$,
$g=\chi_{B\cup(-B)}\in\hat{\bar{{\cal H}}}_G^n$, where $\delta>0$ 
may depend on $\varepsilon>0$, and
$G^n(A_1\,\Delta\,B)<\frac\varepsilon2$. Here $A\Delta B$ 
denotes the symmetric difference of the sets $A$ and $B$, 
and $G^n$ is the $n$-fold direct product of the spectral 
measure~$G$ on the space $R^{n\nu}$. (The above 
properties of the set~$B$ imply that the function
$g=\chi_{B\cup(-B)}\in\hat{\bar{{\cal H}}}_G^n$ satisfies 
the relation $\|g-\chi_A\|<\varepsilon$.)

\medskip
The reduction of Lemma~3.3 to {\it Statement~A} is relatively
simple. Given a function $f\in\bar{{\cal H}}_G^n$ 
we can write $f=f_1+if_2$ with $f_1=\textrm{Re}\,f$,
$if_2=\textrm{Im}\,f$, and both $f_1\in\bar{{\cal H}}_G^n$ 
and $if_2\in\bar{{\cal H}}_G^n$. Hence it is enough to prove
the arbitrarily good approximability of 
$\textrm{Re}\, f\in\bar{{\cal H}}_G^n$ and 
$\textrm{Im}\, f\in\bar{{\cal H}}_G^n$ by a function in 
$\hat{\bar{{\cal H}}}_G^n$. 

Moreover, the real part and imaginary part of the 
function~$f$ can be arbitrarily well approximated by such 
real or imaginary valued functions from the space 
$\bar{{\cal H}}_G^n$ which take only finitely many values, 
and which take a non-zero value only on a bounded set. 
Since we know that the if some functions $f_1,\dots,f_m$ 
from $\bar{{\cal H}}_G^n$ can be approximated arbitrary well
by a function from  $\hat{\bar{{\cal H}}}_G^n$,
then the same relation holds for their linear combination
$\sum_{j=1}^m c_jf_j$ with real coefficients $c_j$, the
good approximability of $\textrm{Re}\,f$ follows from
the first part of {\it Statement~A}. 

The proof of the good approximability of $\textrm{Im}\,f$
is similar, but it demands an additional argument.
We can reduce the statement we want to prove to the good 
approximability of such a function $f$ for which $f(x)=i$ 
on a bounded set $A_0$, $f(x)=-i$ on the set $-A_0$, and 
$f(x)=0$ otherwise. Naturally the sets $A_0$ and $-A_0$ are 
disjoint, but their distance may be zero. Let us observe 
that for any $\varepsilon>0$ there is such a compact set 
$A_1\subset A_0$ for which $G^n(A_1\setminus A_0)<\varepsilon$. 
Then $\rho(A_1,-A_1)>\delta$ with some $\delta>0$, and we can 
reduce the statement about the good approximability of the 
function $\textrm{Im}\,f$ to the good approximability of the
function $g$ which is defined as $g(x)=i$ on the set $A_1$,
$g(x)=-i$ on the set $-A_1$, and it equals zero otherwise.
But the latter statement follows from the second part of
{\it Statement~A} if it is applied for $A=A_1\cup(-A_1)$.

\medskip
To prove {\it Statement A} first I make the following
observation.

For all numbers $M>0$ and $\varepsilon>0$ there is a number 
$\delta=\delta(\varepsilon,M)>0$ such that the set 
\begin{eqnarray*}
K(\delta)&=&
\left\{x=(x^{(1)},\dots, x^{(n)})\colon\; |x_j\pm x_k|<\delta
\textrm{ for a pair } (j,k),\; 1\le j<k\le n\right\} \\
&&\qquad\cap \{x\colon\; x\in R^{n\nu},\,|x|\le M\}
\end{eqnarray*}
satisfies the inequality $G^n(K(\delta))<\varepsilon$.

Similarly, for all $\varepsilon>0$ and $M>0$ 
there is a number $\eta=\eta(\varepsilon,M)>0$ such that
$$
G^n(L(\eta))<\varepsilon \quad \textrm{with }  
L(\eta)=\left(\bigcup_{j=1}^n L_j(\eta)\right)\cap
\{x\colon\; x\in R^{n\nu},\,|x|\le M\},
$$ 
where
$L_j(\eta)=\{(x_1,\dots,x_n)\colon\; x_l\in R^\nu,\;l=1,\dots,n,\;
\rho(x_j,0)\le\eta\}$. 

Indeed, because of condition that $G(\{x\})=0$ for all $x\in R^\nu$
we get by means of the Fubini theorem that for all $j\neq k$, 
$1\le j,k\le n$, $G^n(\{x_j\neq\pm x_k\})=0$. The first statement
follows from this relation, since it implies that the intersection 
of the sets $K(\delta)$ with $\delta\to0$ is contained in a set with
zero $G^n$ measure.

The second statement follows similarly from the relation
$G(\{0\})=0$, since it implies that for all $1\le j\le n$
$G^n(x=(x_1,\dots,x_n)\colon\; x_j=0)=0$.

Since the set $A$ considered in {\it Statement A} is bounded,
the above relations enable us to replace the set $A$ by the
set $A'=A\setminus(K(\delta)\cup L(\eta))$ with a sufficiently small
$\delta>0$ and $\eta>0$ in the formulation of {\it Statement A}.
Let us consider first the first part of {\it Statement A}. Observe 
that the property $A'=-A'$ is preserved. 

We can choose some open rectangles 
\begin{eqnarray*}
D(j)&=&
(a_{(1,1)}(j),b_{(1,1)}(j))\times\cdots\times
(a_{(1,\nu)}(j),b_{(1,\nu)}(j))\times\cdots\\
&&\qquad\cdots \times
(a_{(n,1)}(j),b_{(n,1)}(j))\times\cdots\times
(a_{(n,\nu)}(j),b_{(n,\nu)}(j)),
\end{eqnarray*}
$j=1,\dots,M$ with some number $M>0$ which
satisfy the following relations: \hfill\break
$G^n(x_{(p-1)\nu+s}=\pm a_{(p,s)}(j))=0$ and  
$G^n(x_{(p-1)\nu+s}=\pm b_{(p,s)}(j))=0$ for all $1\le p\le n$
and $1\le s\le\nu$, and also the inequality
$G^n\left(\left(\bigcup_{j=1}^M D(j)\right)\Delta A'\right)
<\frac\varepsilon2$ holds. Let us define the rectangles
$D(-j)=-D(j)$ for all $1\le j\le n$. Since $A'=-A'$ the last
inequality also implies that
$G^n\left(\left(\bigcup_{j=-M}^M D(j)\right)\Delta A'\right)
<\varepsilon$. 

We split up the set $\bigcup_{j=-M}^M D(j)$ (by omitting
some hyperplanes with zero $G_N$ measure) to disjoint open 
rectangles in the following way. First we choose some disjoint
intervals $S_l=(a'_l,b'_l)$, $-P'\le l\le P$ with some $P>0$
in such a way that it satisfies the following properties. 
$S_l=-S_{-l}$, for all $-P\le l\le P$, (in particular, 
$S_0=-S_{0}$). The relations 
$G(x_k=a_l')=G(x_k=b_l')=0$ hold for all $-P\le l\le P$ and 
$1\le k\le \nu$. Besides,
$b'_l-a'_l\le\min(\frac\delta{2n\nu},\frac\eta{2n\nu})$
for all $-P\le l\le P$ with the parameters $\delta$ and 
$\eta$  of those sets $K(\delta)$ and $L(\eta)$ which we 
chose in the definition of the set $A'$, and all edges
$(a_{(p,s)}(j),b_{(p,s)}(j))$ of the rectangles $D(j)$,
$1\le s\le n$, $1\le p\le\nu$, $-M\le j\le M$, (except
finitely many points of the form $a'_l$ or $b'_l$) can be
presented as the union of some elements from the set of
intervals $(a'_l,b'_l)$, $-P\le l\le P$. 

Then we take all those rectangles $D'(k)$ of the form
\begin{eqnarray*}
D'(k) &=&
(a'_{u(1,1,k)},b'_{u(1,1,k)})\times\cdots\times
(a'_{u(1,\nu,k)},b'_{u(1,\nu,k)})\times\cdots\\
&&\qquad \cdots \times
(a'_{u(n,1,k)},b'_{u(n,1,k)})\times\cdots\times
(a'_{u(n,\nu.k)},b'_{u(n,\nu,k)})
\end{eqnarray*}
for which $D'(k)\subset D(j)$ with some $-M\le j\le M$.
The union of these rectangles equals $\bigcup_{j=-M}^M D(j)$ 
minus finitely many hyperplanes of dimension $n\nu-1$ with
zero $G^n$ measure.

In the next step of our construction we preserve those elements
from the set of these rectangles whose intersection with
the set $A'$ is non-empty. Let us reindex the set of these
preserved rectangles $D'(k)$ by the numbers $1\le k\le  M'$ 
with some number $M'$. Clearly we have
$G^n\left(\left(\bigcup_{k=1}^{M'} D'(k)\right)\Delta A'\right)
<\varepsilon$. Let us still define set of rectangles $\cal D$
that consists of those rectangles $\Delta_s\in{\cal B}^\nu$, 
$1\le s\le P$, which are a side of one of the above defined 
rectangles $D'(k)$, $1\le k\le M'$. More precisely $\cal D$ 
consists of those rectangles in $R^\nu$ which can be written 
in the form $\Delta_p=(a'_{u(p,1,k)},b'_{u(p,1,k)})\times\cdots\times
(a'_{u(p,\nu,k)},b'_{u(p,\nu,k)})$, $1\le p\le n$, $1\le k\le M'$,
where the intervals $(a'_{u(p,l,k)},b'_{u(p,l,k)})$ appear in the 
representation of one of the rectangles $D'(k)$, $1\le k\le M'$.

I claim that the class of sets $\cal D$ (with an appropriate
indexation) is a regular system, and if we define the function
$g(x)$ as $g(x)=1$ if $x\in D'(k)$ with some $1\le k\le M'$,
and $g(x)=0$ otherwise, then $g(x)$ is a simple function adapted 
to the regular system~$\cal D$. This fact together with the above
mentioned inequality imply the first part of {\it Statement~A}.

It is clear that $\cal D$ consists of disjoint sets, and if
$\Delta_l\in{\cal D}$, then also $-\Delta_l\in{\cal D}$. We
still have to show that $-\Delta_l\neq\Delta_l$ for all sets 
$\Delta_l\in{\cal D}$. To prove this let us first observe 
that $D'(k)\cap K(\frac\eta2)=\emptyset$ for all $1\le k\le M'$.
Indeed, there is some point $x\in D'(k)\setminus K(\eta)$, because
$D'(k)\cap A'$ is non-empty. As the diameter of $D'(k)$ is less
than $\frac\eta2$ this implies that $D'(k)\cap K(\frac\eta2)=\emptyset$.
Since this relation holds for all sets $D'(k)$, $1\le k\le M'$, 
the definition of the set $K(\frac\eta2)$ and of the class of 
sets $\cal D$ imply that $-\Delta_l\neq\Delta_l$ for all sets 
$\Delta_l\in{\cal D}$. 

To prove that $g(x)$ is a simple function adapted to $\cal D$
we still have to show that for all rectangles
$D'(k)=\Delta_{k_1}\times\cdots\times\Delta_{k_n}$, $1\le k\le M'$,
the relation $k_l\neq\pm k_{l'}$ holds if $l\neq l'$, $1\le l,l'\le n$. 
To prove this statement observe that 
$D'(k)\cap K(\frac \delta2)=\emptyset$ for all $0\le k\le M'$. 
Indeed, there is some point $x\in D'(k)\setminus K(\delta)$, since
$D'(k)\cap A'\neq\emptyset$. Since the length of all edges of
$D'(k)$ is less than $\frac\delta{2n\nu}$, this implies this
statement. Finally this statement together with the definition
of the set $K(\frac\delta2)$ imply the desired property.

The proof of the second part of {\it Statement A} can be proved
with some small modifications of the previous argument. The main
difference is that in this case we start our construction with
a good approximation of the set $A_1$ (and not of $A$) with the
union of some rectangles. Then we take these rectangles $D(j)$
together with their reflection $-D(j)$, and apply the same
procedure as before to get the proof of the second part of
{\it Statement A.} There is still a small additional 
modification in this construction. We choose the rectangles 
$D'(k)$ in our construction with such a little diameter that
guarantees that if one of these rectangles intersects the 
set $A_1$, another one intersects the set $-A_1$, then they 
are disjoint.

\medskip
Let me finally remark that we got such an approximation of
a function $f\in\bar{\cal H}^n_G$ with elementary functions
which are adapted to such a regular system $\cal D$, whose 
elements satisfy the property $G(\partial\Delta_j)=0$ for 
all $\Delta_j\in{\cal D}$, where $\partial\Delta$ denotes
the boundary of the set $\Delta$. I made this remark, because
this means that we have such an approximation in Lemma~3.3
which also satisfies the extra property needed in the proof
of Lemma~6.3. 

\end{document}